\documentclass[12pt,a4paper]{amsart}
\usepackage{amsfonts}
\usepackage{amsmath}
\usepackage{amssymb}

\def\gG{\mathfrak{G}}
\def\gH{\mathfrak{H}}

\def\gM{\mathfrak{M}}
\def\gN{\mathfrak{N}}

\begin{document}

\title[Continuous fields of $\mathbb{C}^*$-algebras and Quantum Groupoids] 
{Continuous fields of $\mathbb{C}^*$-algebras and \\ Lesieur's Quantum Groupoids}
\author{Michel Enock}
\address{Institut de Math\'ematiques de Jussieu, Unit\'{e} Mixte Paris 6 / Paris 7 /
CNRS de Recherche 7586 \\175, rue du Chevaleret, Plateau 7E, F-75013 Paris}
 \email{enock@math.jussieu.fr}
\date{september 06}

\begin{abstract}
In two articles ([L2], [L3]), Franck Lesieur had introduced a notion of quantum groupoid, in the setting of von Neumann algebras, using intensively the notion of pseudo-multiplicative unitary, which had been introduced in a previous article of the author, in collaboration with Jean-Michel Vallin [EV]. We are here studying the analog of Lesieur's construction in a $\mathbb{C}^*$-framework, when the basis of the quantum groupoid is central; in this case, the $\mathbb{C}^*$ structure obtained can be described using the tools of continuous fields of $\mathbb{C}^*$-algebras. This allows us to re-use as quantum groupoids some examples introduced by Blanchard.  \end{abstract}

\maketitle
\newpage
\section{Introduction}
\label{intro}
\subsection{}
 In two articles ([Val2], [Val3]), J.-M. Vallin has introduced two notions (pseudo-multiplicative
unitary, Hopf-bimodule), in order to generalize, up to the groupoid
case, the classical notions of multiplicative unitary [BS] and of Hopf-von Neumann algebras [ES]
which were introduced to describe and explain duality of groups, and leaded to appropriate notions
of quantum groups ([ES], [W1], [W2], [BS], [MN], [W3], [KV1], [KV2], [MNW]). 
\\ In another article [EV], J.-M. Vallin and the author have constructed, from a depth 2 inclusion of
von Neumann algebras $M_0\subset M_1$, with an operator-valued weight $T_1$ verifying a regularity
condition, a pseudo-multiplicative unitary, which leaded to two structures of Hopf bimodules, dual
to each other. Moreover, we have then
construct an action of one of these structures on the algebra $M_1$ such that $M_0$
is the fixed point subalgebra, the algebra $M_2$ given by the basic construction being then
isomorphic to the crossed-product. We construct on $M_2$ an action of the other structure, which
can be considered as the dual action.
\\  If the inclusion
$M_0\subset M_1$ is irreducible, we recovered quantum groups, as proved and studied in former papers
([EN], [E1]).
\\ Therefore, this construction leads to a notion of "quantum groupo\"{\i}d", and a construction of a
duality within "quantum groupo\"{\i}ds". 
\subsection{}
In a finite-dimensional setting, this construction can be
mostly simplified, and is studied in [NV1], [BSz1],
[BSz2], [Sz],[val4], [Val5], and examples are described. In [NV2], the link between these "finite quantum
groupo\"{\i}ds" and depth 2 inclusions of
$II_1$ factors is given. 
\subsection{}
Franck Lesieur introduced in his thesis [L1] a notion of "measured quantum groupoids", in which a modular hypothesis on the basis is required. Mimicking in a wider setting the technics of Kustermans and Vaes [KV], he obtained then a pseudo-multiplicative unitary, which, as in the quantum group case, "contains" all the information of the object (the von Neuman algebra, the coproduct, the antipod, the co-inverse). This construction will be published in an article [L2]. Unfortunately, the axioms chosen by Lesieur don't fit perfectely with the duality (namely, the dual object doesnot fit the modular condition on the basis chosen in [L1] and [L2]), and, for this purpose, Lesieur had to define a wider class, called the "generalized measured quantum groupoids", whose axioms could be described as the analog of [MNW], in which a duality is defined and studied. In [E3] had been shown that, with suitable conditions, the objects constructed in [EV] from depth 2 inclusions, are "generalized measured quantum groupoids". 
\subsection{}
All these constructions had been made in a von Neumann setting, which was natural, once we are dealing with (or thinking of) depth 2 inclusions of von Neumann algebras. But, as for quantum groups, a 
 $\mathbb{C}^*$-version of this theory is to be done, at least to obtain quantum objects similar to locally compact groupoids. Many difficulties exist on that direction : how to define a relative $\mathbb{C}^*$-tensor product ? how to define the analog of operator-valued weights at the $\mathbb{C}^*$ level ?
 \newline
 This article is just a first step in that direction, and is devoted only to the special case when the basis of the measured quantum groupoid is central. 
 \newline
 In that case, we get closed links with the theory of continuous fields of $\mathbb{C}^*$-algebras, as studied by E. Blanchard; using this theory and formalism, we give then a definition of $\mathbb{C}^*$-quantum groupoid with central basis. A procedure is given to associate to such an object a measured quantum groupoid (with a central basis). Conversely, any measured quantum groupoid, with a central basis, contains as a sub-$\mathbb{C}^*$-algebra a $\mathbb{C}^*$-quantum groupoid with central basis; but this last construction is far from being unique; as an example, one must think that finding a dense sub-$\mathbb{C}^*$-algebra of an abelian von Neumann algebra is just a particular case of the problem. 
 \subsection{}
The paper is organized as follows : in chapter \ref{pr}, we give all the preliminaries needed for that theory, mostly Connes-Sauvageot relative tensor product, weights and operator-valued weights on $\mathbb{C}^*$-algebras, and Hopf-bimodules; in chapter \ref{pmu} is recalled the notion of a pseudo-multiplicative unitary, and the Hopf-bimodules associated. In  chapter \ref{quantum} is recalled the notion of measured quantum groupoid (and of generalized measured quantum groupoid), as defined by Lesieur. We show that, when the basis is central, the norm closed algebra $A_n(W)$, which had been defined in [E2], is a $\mathbb{C}^*$-algebra invariant by all the data of the generalized quantum groupoid. 
 \newline
 In chapter \ref{CXC*}, we introduce a notion of "measured continuous field of $\mathbb{C}^*$-algebras", which occurs when the basis is central, and this case is studied in details in chapter \ref{chapcentral}. In particular, we look after a measured quantum groupoid, whose underlying von Neuman algebra itself is abelian; it is then proved that we obtain, in that case, locally compact groupoids. Using the constructions made in chapter \ref{CXC*}, we give, in chapter \ref{C*central} a defintion of a $\mathbb{C}^*$-quantum groupoid, in the case of a central basis, and we show how to construct then a measured quantum groupoid. We finish this article (chapter \ref{fieldqg}) by applying this construction to define a notion of continuous fields of locally compact quantum groups, which was underlying in Blanchard's work [B2]; these are exactly the measured quantum groupoids with central basis, and with a dual which has also a central basis. Blanchard's examples are recalled.  
 \subsection{}
The author is mostly indebted to E. Blanchard, F. Lesieur, S. Vaes, 
L. Va\u{\i}nerman and J.-M. Vallin for many fruitful conversations.

\section{Preliminaries}
\label{pr}
In this chapter are mainly recalled definitions and notations about Connes' spatial
theory (\ref{spatial}, \ref{rel}) and the fiber product construction (\ref{fiber}, \ref{slice})
which are the main technical tools of the theory of measured quantum theory. The definition of Hopf-bimodules is given
(\ref{Hbimod}). 
 \newline
In \ref{C*} are recalled classical results about weights on $\mathbb{C}^*$-algebras, and a standard procedure is defined to go from $\mathbb{C}^*$-algebra weight theory to von Neumann weight theory and vice versa. 

\subsection{Spatial theory [C1], [S2], [T]}
 \label{spatial}
 Let $N$ be a von Neumann algebra, 
and let $\psi$ be a faithful semi-finite normal weight on $N$; let $\gN_{\psi}$, 
$\gM_{\psi}$, $H_{\psi}$, $\pi_{\psi}$, $\Lambda_{\psi}$,$J_{\psi}$, 
$\Delta_{\psi}$,... be the canonical objects of the Tomita-Takesaki construction 
associated to the weight $\psi$. Let $\alpha$  be a non-degenerate normal representation of $N$ on a
Hilbert space
$\mathcal{H}$. We may as well consider $\mathcal{H}$ as a left $N$-module, and write it then
$_\alpha\mathcal{H}$. Following ([C1], definition 1), we define the set of 
$\psi$-bounded elements of $_\alpha\mathcal{H}$ as :
\[D(_\alpha\mathcal{H}, \psi)= \lbrace \xi \in \mathcal{H};\exists C < \infty ,\| \alpha (y) \xi\|
\leq C \| \Lambda_{\psi}(y)\|,\forall y\in \gN_{\psi}\rbrace\]
Then, for any $\xi$ in $D(_\alpha\mathcal{H}, \psi)$, there exists a bounded operator
$R^{\alpha,\psi}(\xi)$ from $H_\psi$ to $\mathcal{H}$,  defined, for all $y$ in $\gN_\psi$ by :
\[R^{\alpha,\psi}(\xi)\Lambda_\psi (y) = \alpha (y)\xi\]
If there is no ambiguity about the representation $\alpha$, we shall write
$R^{\psi}(\xi)$ instead of $R^{\alpha,\psi}(\xi)$. This operator belongs to $Hom_N (H_\psi , \mathcal{H})$;
therefore, for any
$\xi$, $\eta$ in
$D(_\alpha\mathcal{H}, \psi)$, the operator :
\[\theta^{\alpha,\psi} (\xi ,\eta ) =  R^{\alpha,\psi}(\xi)R^{\alpha,\psi}(\eta)^*\]
belongs to $\alpha (N)'$; moreover, $D(_\alpha\mathcal{H}, \psi)$ is dense ([C1], lemma 2), stable
under $\alpha (N)'$, and the linear span generated by the operators $\theta^{\alpha,\psi} (\xi ,\eta
)$ is a weakly dense ideal in $\alpha (N)'$. We shall write $\mathcal{K}_{\alpha, \psi}$ the norm closure of this ideal, which is also weakly dense in $\alpha (N)'$. 
 \newline
 With the same hypothesis, the operator :
\[<\xi,\eta>_{\alpha,\psi} = R^{\alpha,\psi}(\eta)^* R^{\alpha,\psi}(\xi)\]
belongs to $\pi_{\psi}(N)'$. Using Tomita-Takesaki's theory, this last algebra is equal to
$J_\psi
\pi_{\psi}(N)J_\psi$, and therefore anti-isomorphic to $N$ (or isomorphic to the opposite von
Neumann algebra $N^o$). We shall consider now $<\xi,\eta>_{\alpha,\psi}$ as an element of $N^o$, and
the linear span generated by these operators is a dense algebra in $N^o$. More precisely ([C], lemma 4, and [S1], lemme 1.5), we get that $<\xi, \eta>_{\alpha, \psi}^o$ belongs to $\gM_\psi$, and that :
 \[\Lambda_{\psi}(<\xi, \eta>_{\alpha, \psi}^o)=J_\psi R^{\alpha, \psi}(\xi)^*\eta\]
 \newline
 If $y$ in $N$ is analytical with respect to $\psi$, and if $\xi\in D(_\alpha\mathcal{H}, \psi)$, then we get that $\alpha(y)\xi$ belongs to $D(_\alpha\mathcal{H}, \psi)$ and that :
 \[R^{\alpha,\psi}(\alpha(y)\xi)=R^{\alpha,\psi}(\xi)J_\psi\sigma_{-i/2}^{\psi}(y^*)J_\psi\]
  So, if $\eta$ is another $\psi$-bounded element of $_\alpha\mathcal H$, we get :
 \[<\alpha(y)\xi, \eta>_{\alpha, \psi}^o=\sigma_{i/2}^\psi(y)<\xi, \eta>_{\alpha, \psi}^o\]
There exists ([C], prop.3) a family  $(e_i)_{i\in I}$ of 
$\psi$-bounded elements of $_\alpha\mathcal H$, such that
\[\sum_i\theta^{\alpha, \psi} (e_i ,e_i )=1\]
Such a family will be called an $(\alpha,\psi )$-basis of $\mathcal H$. 
 \newline
 In that situation, let us consider, for all $n\in \mathbb{N}$ and finite $J\subset I$ with $|J|=n$, the $(1,n)$ matrix $R_J=(R^{\alpha, \nu}(e_i))_{i\in J}$. As $R_JR_J^*\leq1$, we get that $\|R_J\|\leq 1$, and that the $(n,n)$ matrix $(<e_i, e_j>_{\alpha, \nu})_{i,j\in J}\in M_n(N^o)$ is less than the unit matrix. 
 \newline
 It is possible ([EN] 2.2) to construct 
an $(\alpha,\psi )$-basis of $\mathcal H$, $(e_i)_{i\in I}$, such that the
operators $R^{\alpha, \psi}(e_i)$ are partial isometries with final supports 
$\theta^{\alpha, \psi}(e_i ,e_i )$ 2 by 2 orthogonal, and such that, if $i\neq j$, then 
$<e_i ,e_j>_{\alpha, \psi}=0$. Such a family will be called an $(\alpha, \psi)$-orthogonal basis of $\mathcal H$. 
\newline
 We have, then :
 \[R^{\alpha, \psi}(\xi)=\sum_i\theta^{\alpha, \psi}(e_i, e_i)R^{\alpha, \psi}(\xi)=\sum_iR^{\alpha, \psi}(e_i)<\xi, e_i>_{\alpha, \psi}\]
 \[<\xi, \eta>_{\alpha, \psi}=\sum_i<\eta, e_i>_{\alpha, \psi}^*<\xi, e_i>_{\alpha, \psi}\]
 \[\xi=\sum_iR^{\alpha, \psi}(e_i)J_\psi\Lambda_\psi(<\xi, e_i>^o_{\alpha, \psi})\]
 the sums being weakly convergent. Moreover, we get that, for all $n$ in $N$, $\theta^{\alpha, \psi}(e_i, e_i)\alpha(n)e_i=\alpha(n)e_i$, and $\theta^{\alpha, \psi}(e_i, e_i)$ is the orthogonal projection on the closure of the subspace $\{\alpha(n)e_i, n\in N\}$. 
 \newline
Let $\beta$ be a normal non-degenerate
anti-representation
 of
$N$ on
$\mathcal{H}$. We may then as well consider $\mathcal{H}$ as a right $N$-module, and write it $\mathcal{H}_\beta$, or
consider
$\beta$ as a normal non-degenerate representation of the opposite von Neumann algebra $N^o$, and
consider
$\mathcal{H}$ as a left $N^o$-module. 
\newline
We can then define on $N^o$ the opposite faithful
 semi-finite normal weight $\psi ^o$; we have $\gN_{\psi ^o}=\gN_\psi^*$, and 
 the Hilbert space $H_{\psi ^o}$ will be, as usual, identified with $H_\psi$, 
by the identification, for all $x$ in $\gN_\psi$, of 
$\Lambda_{\psi ^o}(x^*)$ with $J_\psi \Lambda_\psi (x)$.
\newline
From these remarks, we infer that the set of 
$\psi ^o$-bounded elements of
$\mathcal{H}_\beta$ is :
\[D(\mathcal{H}_\beta, \psi^o) = \lbrace\xi\in\mathcal{H} ;\exists C < \infty ,
\|\beta (y^*)\xi\|
\le C \| \Lambda_{\psi}(y)\|,\forall y\in \gN_{\psi}\rbrace\]
and, for any $\xi$ in $D(\mathcal{H}_\beta, \psi^o)$ and $y$ in $\gN_\psi$,
 the bounded operator $R^{\beta,\psi^{o}}(\xi )$ is given by the formula :
\[R^{\beta,\psi{^o}}(\xi)J_{\psi}\Lambda_\psi (y) = \beta (y^*)\xi\]
This operator belongs to $Hom_{N^{o}}(H_\psi ,\mathcal{H} )$.
Moreover, $D(\mathcal{H}_\beta, \psi^o)$ is dense, stable under $\beta( N)'=P$, 
and, for all $y$ in $P$, we have :
\[R^{\beta,\psi{^o}}(y\xi)= yR^{\beta,\psi{^o}}(\xi)\]
Then, for any $\xi$, $\eta$ in $D(\mathcal{H}_\beta, \psi^o)$, the operator 
\[\theta^{\beta, \psi^{o}}(\xi ,\eta )=R^{\beta,\psi{^o}}(\xi)R^{\beta,\psi{^o}}(\eta)^*\] 
belongs to $P$, and the linear span generated by these operators is a dense ideal in $P$;
moreover, the operator-valued product 
$<\xi ,\eta >_{\beta,\psi^o}= R^{\beta,\psi{^o}}(\eta)^*R^{\beta,\psi{^o}}(\xi)$
belongs to $\pi_\psi (N)$; we shall consider now, for simplification, that $<\xi ,\eta
>_{\beta,\psi^o}$ belongs to $N$, and the linear span generated by these operators is a dense algebra
in $N$, stable under multiplication by analytic elements with respect to $\psi$. More precisely, $<\xi ,\eta
>_{\beta,\psi^o}$ belongs to $\gM_\psi$ ([C], lemma 4) and we have ([S1], lemme 1.5) 
\[\Lambda_\psi(<\xi, \eta>_{\beta,\psi^o})=R^{\beta,\psi{^o}}(\eta)^*\xi\]
A $(\beta,\psi^o)$-basis of $\mathcal H$ is a family 
$(e_i )_{i\in I}$ of $\psi^o$-bounded elements of $\mathcal H_\beta$, such that 
\[\sum_i\theta^{\beta, \psi^{o}}(e_i ,e_i )=1\]
We have then, for all $\xi$ in $D(\mathcal H_\beta)$ :
\[\xi=\sum_i R^{\beta,\psi^o}(e_i)\Lambda_\psi(<\xi, e_i>_{\beta,\psi^o})\]
It is possible to choose the $(e_i )_{i\in I}$ such that
the $R^{\beta, \psi{^o}}(e_i)$ are partial isometries, with final supports
$\theta^{\beta, \psi^{o}}(e_i ,e_i )$ 2 by 2 orthogonal, and $<e_i, e_j>_{\beta, \psi^o}=0$ if $i\neq j$; such
a family will be then called a $(\beta, \psi^o)$-orthogonal basis of $\mathcal H$. We have then 
\[R^{\beta, \psi{^o}}(e_i)=\theta^{\beta, \psi^{o}}(e_i ,e_i )R^{\beta, \psi{^o}}(e_i)=
R^{\beta, \psi{^o}}(e_i)<e_i, e_i>_{\beta, \psi^o}\]
 Moreover, we get that, for all $n$ in $N$, and for all $i$, we have :
\[\theta^{\beta, \psi^{o}}(e_i ,e_i )\beta (n)e_i=\beta (n)e_i\]
and that  $\theta^{\beta, \psi^{o}}(e_i ,e_i )$ is the orthogonal projection on the closure of the subspace $\{\beta (n)e_i, n\in N\}$.

\subsection{Jones' basic construction and operator-valued weights}
\label{basic}
Let $M_0\subset M_1$ be an inclusion of von Neumann algebras  (for simplification, these algebras will be supposed to be $\sigma$-finite), equipped with a normal faithful semi-finite operator-valued weight $T_1$ from $M_1$ to $M_0$ (to be more precise, from $M_1^{+}$ to the extended positive elements of $M_0$ (cf. [T] IX.4.12)). Let $\psi_0$ be a normal faithful semi-finite weight on $M_0$, and $\psi_1=\psi_0\circ T_1$; for $i=0,1$, let $H_i=H_{\psi_i}$, $J_i=J_{\psi_i}$, $\Delta_i=\Delta_{\psi_i}$ be the usual objects constructed by the Tomita-Takesaki theory associated to these weights. Following ([J], 3.1.5(i)), the von Neumann algebra $M_2=J_1M'_0J_1$ defined on the Hilbert space $H_1$ will be called the basic construction made from the inclusion $M_0\subset M_1$. We have $M_1\subset M_2$, and we shall say that the inclusion $M_0\subset M_1\subset M_2$ is standard.  
 \newline
Following ([EN] 10.6), for $x$ in $\gN_{T_1}$, we shall define $\Lambda_{T_1}(x)$ by the following formula, for all $z$ in $\gN_{\psi_{0}}$ :
\[\Lambda_{T_1}(x)\Lambda_{\psi_{0}}(z)=\Lambda_{\psi_1}(xz)\]
Then, $\Lambda_{T_1}(x)$ belongs to $Hom_{M_{0}^o}(H_{0}, H_1)$; if $x$, $y$ belong to $\gN_{T_1}$, then $\Lambda_{T_1}(x)^*\Lambda_{T_1}(y)=T_1(x^*y)$, and $\Lambda_{T_1}(x)\Lambda_{T_1}(y)^*$ belongs to $M_{2}$.
\newline
 Using then Haagerup's construction ([T], IX.4.24), it is possible to construct a normal semi-finite faithful operator-valued weight $T_2$ from $M_2$ to $M_1$ ([EN], 10.7), which will be called the basic construction made from $T_1$. If $x$, $y$ belong to $\gN_{T_1}$, then $\Lambda_{T_1}(x)\Lambda_{T_1}(y)^*$
belongs to $\gM_{T_{2}}$, and $T_{2}(\Lambda_{T_1}(x)\Lambda_{T_1}(y)^*)=xy^*$. 
\newline
By Tomita-Takesaki theory, the Hilbert space $H_1$ bears a natural structure of $M_1-M_1^o$-bimodule, and, therefore, by restriction, of $M_0-M_0^o$-bimodule. Let us write $r$ for the canonical representation of $M_0$ on $H_1$, and $s$ for the canonical antirepresentation given, for all $x$ in $M_0$, by $s(x)=J_1r(x)^*J_1$. Let us have now a closer look to the subspaces $D(H_{1s}, \psi_0^o)$ and $D(_rH_1, \psi_0)$. If $x$ belongs to $\gN_{T_1}\cap\gN_{\psi_1}$, we easily get that $J_1\Lambda_{\psi_1}(x)$ belongs to $D(_rH_1, \psi_0)$, with :
\[R^{r, \psi_0}(J_1\Lambda_{\psi_1}(x))=J_1\Lambda_{T_1}(x)J_0\]
and $\Lambda_{\psi_1}(x)$ belongs to $D(H_{1s}, \psi_0)$, with :
\[R^{s, \psi_0^o}(\Lambda_{\psi_1}(x))=\Lambda_{T_1}(x)\]
In ([E3], 2.3) was proved that the subspace $D(H_{1s}, \psi_0^o)\cap D(_rH_1, \psi_0)$ is dense in $H_1$; let us write down and precise this result :

\subsubsection{{\bf Proposition}}
\label{propbasic}
{\it  Let us keep on the notations of this paragraph; let $A_{\psi_1, T_1}$ be the algebra made of elements $x$ in $\gN_{\psi_1}\cap\gN_{T_1}\cap\gN_{\psi_1}^*\cap\gN_{T_1}^*$, analytical with respect to $\psi_1$, and such that, for all $z$ in $\mathbb{C}$, $\sigma^{\psi_1}_z(x_n)$ belongs to $\gN_{\psi_1}\cap\gN_{T_1}\cap\gN_{\psi_1}^*\cap\gN_{T_1}^*$. Then :
\newline
(i) for any $x$ in  $A_{\psi_1, T_1}$, $\Lambda_{\psi_1}(x)$ belongs to $D(H_{1s}, \psi_0)\cap D(_rH_1, \psi_0)$;
\newline
(ii) for any $\xi$ in $D(H_{1s}, \psi_0^o))$, there exists a sequence $x_n$ in $A_{\psi_1, T_1}$ such that $\Lambda_{T_1}(x_n)=R^{s, \psi_0^o}(\Lambda_{\psi_1}(x))$ is weakly converging to $R^{s, \psi_0^o}(\xi)$ and $\Lambda_{\psi_1}(x_n)$ is converging to $\xi$;  
\newline
(iii) if $M_0$ is abelian, then we get that $\Lambda_{T_1}(x_n)=R^{s, \psi_0^o}(\Lambda_{\psi_1}(x_n))$ is norm converging to $R^{s, \psi_0^o}(\xi)$.}
\begin{proof}
We get in ([E3], 2.3) an increasing sequence of projections $p_n$ in $M_1$, converging to $1$, and elements $x_n$ in $A_{\psi_1, T_1}$ such that $\Lambda_{\psi_1}(x_n)=p_n\xi$. So, (i) and (ii) were obtained in ([E3], 2.3) from this construction. More precisely, we get that :
\begin{eqnarray*}
T_1(x_n^*x_n)&=&<R^{s, \psi^o}(\Lambda_{\psi_1}(x_n)), R^{s, \psi_0^o}(\Lambda_{\psi_1}(x_n))>_{s, \psi_0^o}\\
&=&<p_n\xi, p_n\xi>_{s, \psi_0^o}\\
&=&R^{s, \psi^o}(\xi)^*p_nR^{s, \psi^o}(\xi)
\end{eqnarray*}
which is increasing and weakly converging to $<\xi, \xi>_{s, \psi_0^o}$. Let us suppose now that $M_0$ is abelian, and let $X$ be the spectrum of the $\mathbb{C}^*$-algebra generated by all elements $<\eta_1, \eta_2>_{s, \psi_0^o}$, when $\eta_1$ and $\eta_2$ are in $D(H_{1s}, \psi_0)$. Using Dini's theorem in $C_0(X)$, we get that $T_1(x_n^*x_n)$ is norm converging to $<\xi, \xi>_{s, \psi_0^o}$, and that :
\[\|\Lambda_{T_1}(x_n)-R^{s, \psi_0^o}(\xi)\|^2=\|T_1(x_n^*x_n)-<\xi, \xi>_{s, \psi_0^o}\|\]
is converging to $0$. \end{proof}


\subsection{Relative tensor product [C1], [S2], [T]}
\label{rel}
Using the notations of \ref{spatial}, let now $\mathcal{K}$ be another Hilbert space on which there exists
a non-degenerate representation
$\gamma$ of
$N$. Following J.-L. Sauvageot ([S2], 2.1), we define
the relative tensor product $\mathcal{H}\underset{\psi}{_\beta\otimes_\gamma}\mathcal{K}$ as the
Hilbert space obtained from the algebraic tensor product $D(\mathcal{H}_\beta ,\psi^o )\odot
\mathcal{K} $ equipped with the scalar product defined, for $\xi_1$, $\xi_2$ in $D(\mathcal{H}_\beta
,\psi^o )$,
$\eta_1$, $\eta_2$ in $\mathcal{K}$, by 
\[(\xi_1\odot\eta_1 |\xi_2\odot\eta_2 )=(\gamma(<\xi_1 ,\xi_2 >_{\beta,\psi^o})\eta_1 |\eta_2 )\]
where we have identified $N$ with $\pi_\psi (N)$ to simplifly the notations.
\newline
The image of $\xi\odot\eta$ in $\mathcal{H}\underset{\psi}{_\beta\otimes_\gamma}\mathcal{K}$ will be
denoted by
$\xi\underset{\psi}{_\beta\otimes_\gamma}\eta$. We shall use intensively this construction; one
should bear in mind that, if we start from another faithful semi-finite normal weight $\psi '$, we
get another Hilbert space $\mathcal{H}\underset{\psi'}{_\beta\otimes_\gamma}\mathcal{K}$; there exists an isomorphism $U^{\psi, \psi'}_{\beta, \gamma}$ from $\mathcal{H}\underset{\psi}{_\beta\otimes_\gamma}\mathcal{K}$ to $\mathcal{H}\underset{\psi'}{_\beta\otimes_\gamma}\mathcal{K}$, which is unique up to some functorial property ([S2], 2.6) (but this isomorphism
does not send  $\xi\underset{\psi}{_\beta\otimes_\gamma}\eta$ on
$\xi\underset{\psi'}{_\beta\otimes_\gamma}\eta$ !). 
\newline
When no confusion is possible about the representation and the anti-representation, we shall write
$\mathcal{H}\otimes_{\psi}\mathcal{K}$ instead of
$\mathcal{H}\underset{\psi}{_\beta\otimes_\gamma}\mathcal{K}$, and $\xi\otimes_\psi\eta$ instead
of
$\xi\underset{\psi}{_\beta\otimes_\gamma}\eta$.
\newline
For any $\xi$ in $D(\mathcal{H}_\beta,
\psi^o)$, we define the bounded linear application $\lambda_\xi^{\beta, \gamma}$ from $\mathcal
K$ to
$\mathcal{H}\underset{\psi}{_\beta\otimes_\gamma}\mathcal{K}$ by, for all $\eta$ in $\mathcal K$,
$\lambda_\xi^{\beta, \gamma} (\eta)=\xi\underset{\psi}{_\beta\otimes_\gamma}\eta$. We shall write
$\lambda_\xi$ if no confusion is possible. We get ([EN], 3.10) :
\[\lambda_\xi^{\beta, \gamma}=R^{\beta, \psi^o}(\xi)\otimes_\psi 1_\mathcal K\]
where we recall the canonical identification (as left $N$-modules) of 
$L^2(N)\otimes_\psi\mathcal K$ with $\mathcal K$. We have :
\[(\lambda_\xi^{\beta, \gamma})^*\lambda_\xi^{\beta, \gamma}=\gamma(<\xi, \xi>_{\beta, \psi^o})\]
\newline
In ([S1] 2.1), the relative tensor product
$\mathcal{H}\underset{\psi}{_\beta\otimes_\gamma}\mathcal{K}$ is defined also, if
$\xi_1$, $\xi_2$ are in $\mathcal{H}$, $\eta_1$, $\eta_2$ are in $D(_\gamma\mathcal{K},\psi)$, by the
following formula :
\[(\xi_1\odot\eta_1 |\xi_2\odot\eta_2 )= (\beta(<\eta_1, \eta_2>_{\gamma,\psi})\xi_1 |\xi_2)\]
which leads to the the definition of a relative flip $\sigma_\psi$ which will be an isomorphism from
$\mathcal{H}\underset{\psi}{_\beta\otimes_\gamma}\mathcal{K}$ onto
$\mathcal{K}\underset{\psi^o}{_\gamma\otimes _\beta}\mathcal{H}$, defined, for any 
$\xi$ in $D(\mathcal{H}_\beta ,\psi^o )$, $\eta$ in $D(_\gamma \mathcal{K},\psi)$, by :
\[\sigma_\psi (\xi\otimes_\psi\eta)=\eta\otimes_{\psi^o}\xi\]
This allows us to define a relative flip $\varsigma_\psi$ from
$\mathcal{L}(\mathcal{H}\underset{\psi}{_\beta\otimes_\gamma}\mathcal{K})$ to $\mathcal{L}
(\mathcal{K}\underset{\psi^o}{_\gamma\otimes_\beta}\mathcal{H})$ which sends
$X$ in 
$\mathcal{L}(\mathcal{H}\underset{\psi}{_\beta\otimes_\gamma}\mathcal{K})$ onto
$\varsigma_\psi(X)=\sigma_\psi X\sigma_\psi^*$. Starting from another faithful semi-finite normal
weight $\psi'$, we get a von Neumann algebra
$\mathcal{L}(\mathcal{H}\underset{\psi'}{_\beta\otimes_\gamma}\mathcal{K})$ which is isomorphic to
$\mathcal{L}(\mathcal{H}\underset{\psi}{_\beta\otimes_\gamma}\mathcal{K})$, and a von Neumann
algebra $\mathcal{L} (\mathcal{K}\underset{\psi^{'o}}{_\gamma\otimes_\beta}\mathcal{H})$ which is
isomorphic to
$\mathcal{L} (\mathcal{K}\underset{\psi^o}{_\gamma\otimes_\beta}\mathcal{H})$; as we get that :
 \[\sigma_{\psi'}\circ U^{\psi, \psi'}_{\beta, \gamma}=U^{\psi^o, \psi'^o}_{\gamma, \beta}\]
 we see that these isomorphisms exchange $\varsigma_\psi$ and
$\varsigma_{\psi '}$. Therefore, the homomorphism $\varsigma_{\psi}$ can be denoted $\varsigma_N$
without any reference to a specific weight.
\newline
We may define, for any $\eta$ in $D(_\gamma\mathcal{K}, \psi)$, an application
$\rho_\eta^{\beta, \gamma}$ from $\mathcal H$ to
$\mathcal{H}\underset{\psi}{_\beta\otimes_\gamma}\mathcal{K}$ by
$\rho_\eta^{\beta, \gamma} (\xi)=\xi\underset{\psi}{_\beta\otimes_\gamma}\eta$. We shall write
$\rho_\eta$ if no confusion is possible. We get that :
\[(\rho_\eta^{\beta, \gamma})^*\rho_\eta^{\beta, \gamma}=\beta(<\eta, \eta>_{\gamma, \psi})\]
We recall, following
([S2], 2.2b) that, for all $\xi$ in $\mathcal{H}$, $\eta$ in $D(_\gamma\mathcal{K},\psi)$, $y$ in
$N$, analytic with respect to $\psi$, we have :
\[\beta (y)\xi \otimes_\psi\eta=\xi\otimes_\psi\gamma(\sigma^\psi_{-i/2}(y))\eta\]
Let $x$ be an element
of
$\mathcal{L}(\mathcal{H})$, commuting with the right action of $N$ on $\mathcal{H}_\beta$ (i.e. $x\in\beta(N)'$). It
is possible to define an operator $x\underset{\psi}{_\beta\otimes_\gamma} 1_{\mathcal{K}}$ on
$\mathcal{H}\underset{\psi}{_\beta\otimes_\gamma}
\mathcal{K}$. We can easily evaluate $\|x\underset{\psi}{_\beta\otimes_\gamma} 1_{\mathcal{K}}\|$, for any finite $J\subset I$, for any $\eta_i$ in $\mathcal K$, we have :
 \begin{multline*}
 ((x^*x\underset{\psi}{_\beta\otimes_\gamma} 1_{\mathcal{K}})(\Sigma_{i\in J}e_i\underset{\psi}{_\beta\otimes_\gamma}\eta_i)|(\Sigma_{i\in J}e_i\underset{\psi}{_\beta\otimes_\gamma}\eta_i))=\\
 =\Sigma_{i\in J}(\gamma(<xe_i, xe_i>_{\beta, \psi^o})\eta_i|\eta_i)\\
 \leq \|x\|^2\Sigma_{i\in J}(\gamma(<e_i, e_i>_{\beta, \psi^o})\eta_i|\eta_i)\
 =\|x\|^2\|\Sigma_{i\in J}e_i\underset{\psi}{_\beta\otimes_\gamma}\eta_i\|
 \end{multline*}
 from which we get $\|x\underset{\psi}{_\beta\otimes_\gamma} 1_{\mathcal{K}}\|\leq\|x\|$. 
 \newline
 By the same way, if $y$ commutes with the left action of $N$ on
$_\gamma\mathcal{K}$ (i.e. $y\in\gamma(N)'$), it is possible to define
$1_{\mathcal{H}}\underset{\psi}{_\beta\otimes_\gamma}y$ on
$\mathcal{H}\underset{\psi}{_\beta\otimes_\gamma} \mathcal{K}$, and by composition, it is possible
to define then
$x\underset{\psi}{_\beta\otimes_\gamma} y$. If we start from another faithful semi-finite normal
weight $\psi '$, the canonical isomorphism $U^{\psi, \psi'}_{\beta, \gamma}$ from $\mathcal{H}\underset{\psi}{_\beta\otimes_\gamma}
\mathcal{K}$ to $\mathcal{H}\underset{\psi'}{_\beta\otimes_\gamma} \mathcal{K}$ sends
$x\underset{\psi}{_\beta\otimes_\gamma} y$ on $x\underset{\psi'}{_\beta\otimes_\gamma} y$ ([S2],
2.3 and 2.6); therefore, this operator can be denoted $x\underset{N}{_\beta\otimes_\gamma} y$
without any reference to a specific weight, and we get $\|x\underset{N}{_\beta\otimes_\gamma} y\|\leq\|x\|\|y\|$. 
 \newline
 With the notations of \ref{spatial}, let $(e_i)_{i\in I}$ a $(\beta, \psi^o)$-orthogonal basis of $\mathcal H$; 
 let us remark that, for all $\eta$ in $\mathcal K$, we have :
 \[e_i\underset{\psi}{_\beta\otimes_\gamma}\eta=e_i\underset{\psi}{_\beta\otimes_\gamma}\gamma(<e_i, e_i>_{\beta, \psi^o})\eta\]
On the other hand, $\theta^{\beta, \psi^o}(e_i, e_i)$ is an orthogonal projection, and so is $\theta^{\beta, \psi^o}(e_i, e_i)\underset{N}{_\beta\otimes_\gamma}1$; this last operator is the projection on the subspace $e_i\underset{\psi}{_\beta\otimes_\gamma}\gamma(<e_i, e_i>_{\beta, \psi^o})\mathcal K$ ([E2], 2.3) and, therefore, we get that $\mathcal H\underset{\psi}{_\beta\otimes_\gamma}\mathcal K$ is the orthogonal sum of the subspaces $e_i\underset{\psi}{_\beta\otimes_\gamma}\gamma(<e_i, e_i>_{\beta, \psi^o})\mathcal K$;  for any $\Xi$ in $\mathcal{H}\underset{\psi}{_\beta\otimes_\gamma}
\mathcal{K}$, there exist $\xi_i$ in $\mathcal K$, such that $\gamma(<e_i, e_i>_{\beta, \psi^o})\xi_i=\xi_i$ and $\Xi=\sum_i e_i\underset{\psi}{_\beta\otimes_\gamma}\xi_i$, from which we get that $\sum_i\|\xi_i\|^2=\|\Xi\|^2$. 
 \newline
 Let $A_1$ a non-degenerate $\mathbb{C}^*$-subalgebra of $\beta (N)'$, $A_2$ a non-degenerate $\mathbb{C}^*$-subalgebra of $\gamma (N)'$; we shall denote by $A_1\underset{N}{_\beta\otimes_\gamma}A_2$ the $\mathbb{C}^*$-algebra generated by all the operators $x\underset{N}{_\beta\otimes_\gamma} y$, where $x\in A_1$ and $y\in A_2$\vspace{2mm}.  
\newline
Let us suppose now that $\mathcal{K}$ is a $N-P$ bimodule; that means that there exists a von
Neumann algebra $P$, and a non-degenerate normal anti-representation $\epsilon$ of $P$ on
$\mathcal{K}$, such that
$\epsilon (P)\subset\gamma (N)'$. We shall write then $_\gamma\mathcal{K}_\epsilon$. If $y$ is in $P$, we
have seen that it is possible to define then the operator
$1_{\mathcal{H}}\underset{\psi}{_\beta\otimes_\gamma}\epsilon (y)$ on
$\mathcal{H}\underset{\psi}{_\beta\otimes_\gamma}\mathcal{K}$, and we define this way a
non-degenerate normal antirepresentation of $P$ on
$\mathcal{H}\underset{\psi}{_\beta\otimes_\gamma}\mathcal{K}$, we shall call again $\epsilon$ for
simplification. If $\mathcal H$ is a $Q-N$ bimodule, then
$\mathcal{H}\underset{\psi}{_\beta\otimes_\gamma}\mathcal{K}$ becomes a $Q-P$ bimodule (Connes'
fusion of bimodules).
\newline
Taking a faithful semi-finite normal weight
$\nu$  on $P$, and a left $P$-module $_{\zeta}\mathcal{L}$ (i.e. a Hilbert space $\mathcal{L}$ and a normal
non-degenerate representation $\zeta$ of $P$ on $\mathcal{L}$), it is possible then to define
$(\mathcal{H}\underset{\psi}{_\beta\otimes_\gamma}\mathcal{K})\underset{\nu}{_\epsilon\otimes_\zeta}\mathcal{L}$.
Of course, it is possible also to consider the Hilbert space
$\mathcal{H}\underset{\psi}{_\beta\otimes_\gamma}(\mathcal{K}\underset{\nu}{_\epsilon\otimes_\zeta}\mathcal{L})$.
It can be shown that these two Hilbert spaces are isomorphics as $\beta (N)'-\zeta
(P)^{'o}$-bimodules. (In ([V1] 2.1.3), the proof, given for $N=P$ abelian can be used, without
modification, in that wider hypothesis). We shall write then
$\mathcal{H}\underset{\psi}{_\beta\otimes_\gamma}\mathcal{K}\underset{\nu}{_\epsilon\otimes_\zeta}\mathcal{L}$
without parenthesis, to emphazise this coassociativity property of the relative tensor
product.
 \newline
 If $\pi$ denotes the canonical left representation of $N$ on the Hilbert space $L^2(N)$, then it is straightforward to verify that the application which sends, for all $\xi$ in $\mathcal {H}$, $\chi$ normal faithful semi-finite weight on $N$, and $x$ in $\gN_\chi$, the vector $\xi{}_\beta\underset{\chi}{\otimes}{}_\pi J_\chi\Lambda_\chi (x)$ on $\beta(x^*)\xi$, gives an isomorphism of $\mathcal{H}{}_\beta\underset{\chi}{\otimes}{}_\pi L^2(N)$ on $\mathcal{H}$, which will send the antirepresentation of $N$ given by $n\mapsto 1_{\mathcal{H}}{}_\beta\underset{\chi}{\otimes}{}_\pi J_\chi n^*J_\chi$ on $\beta$
 \newline
 If $\mathcal K$ is a Hilbert space on which there exists a non-degenerate representation $\gamma$ of $N$, then $\mathcal K$ is a $N-\gamma(N)'^o$ bi-module, and the conjugate Hilbert space $\overline{\mathcal K}$ is a $\gamma(N)'-N^o$ bimodule, and, ([S2]), for any normal faithful semi-finite weight $\phi$ on $\gamma(N)'$, the fusion $_\gamma\mathcal K\underset{\phi^o}{\otimes}\overline{\mathcal K}_\gamma$ is isomorphic to the standard space $L^2(N)$, equipped with its standard left and right representation. 
 \newline
 Using that remark, one gets for any $x\in \beta(N)'$ :
 \[\|x\underset{N}{_\beta\otimes_\gamma}1_{\mathcal K}\|\leq\|x\underset{N}{_\beta\otimes_\gamma}1_{\mathcal K}\underset{\gamma(N)'^o}{\otimes}1_{\overline{\mathcal K}}\|=\|x\underset{N}{_\beta\otimes}1_{L^2(N)}\|=\|x\|\]
from which we have $\|x\underset{N}{_\beta\otimes_\gamma}1_{\mathcal K}\|=\|x\|$ \vspace{2mm}. 
\newline
If $\mathcal H$ and $\mathcal K$ are finite-dimensional Hilbert spaces, the relative tensor product 
$\mathcal{H}\underset{\psi}{_\beta\otimes_\gamma}
\mathcal{K}$ can be identified with a subspace of the tensor Hilbert space 
$\mathcal{H}\otimes
\mathcal{K}$ ([EV] 2.4), the projection on which belonging to $\beta (N)\otimes\gamma (N)$.

\subsection{Fiber product [V1], [EV]} 
\label{fiber}
Let us follow the notations of \ref{rel}; let now
$M_1$ be a von Neumann algebra on $\mathcal{H}$, such that $\beta (N)\subset
M_1$, and $M_2$ be a von Neumann algebra on $\mathcal{K}$, such that $\gamma (N)\subset
M_2$. The von Neumann algebra generated by all elements $x\underset{N}{_\beta\otimes_\gamma} y$,
where
$x$ belongs to $M'_1$, and $y$ belongs $M'_2$ will be denoted
$M'_1\underset{N}{_\beta\otimes_\gamma} M'_2$ (or $M'_1\otimes_N M'_2$ if no confusion if
possible), and will be called the relative tensor product of
$M'_1$ and $M'_2$ over $N$. The commutant of this algebra will be denoted 
$M_1\underset{N}{_\beta *_\gamma} M_2$ (or $M_1*_N M_2$ if no confusion is possible) and called the
fiber product of $M_1$ and
$M_2$, over
$N$. It is straightforward to verify that, if $P_1$ and $P_2$ are two other von Neumann
algebras satisfying the same relations with $N$, we have 
\[M_1*_N M_2\cap P_1*_N P_2=(M_1\cap P_1)*_N (M_2\cap P_2)\]
Moreover, we get that $\varsigma_N (M_1\underset{N}{_\beta *_\gamma}
M_2)=M_2\underset{N^o}{_\gamma *_\beta}M_1$.
\newline
In particular, we have :
\[(M_1\cap \beta (N)')\underset{N}{_\beta\otimes_\gamma} (M_2\cap \gamma (N)')\subset
M_1\underset{N}{_\beta *_\gamma} M_2\] and :
\[M_1\underset{N}{_\beta *_\gamma} \gamma(N)=(M_1\cap\beta (N)')\underset{N}{_\beta\otimes_\gamma} 1\]
More generally, if
$\beta$ is a non-degenerate normal involutive antihomomorphism from
$N$ into a von Neumann algebra
$M_1$, and
$\gamma$ a non-degenerate normal involutive homomorphism from $N$ into a von Neumann
algebra
$M_2$, it is possible
to define, without any reference to a specific Hilbert space, a von Neumann algebra
$M_1\underset{N}{_\beta *_ \gamma}M_2$. 
\newline
Moreover, if now $\beta '$ is a non-degenerate normal involutive antihomomorphism from $N$ into
another von Neumann algebra
$P_1$,
$\gamma '$ a non-degenerate normal involutive homomorphism from $N$ into another
von Neumann algebra $P_2$, $\Phi$ a normal involutive homomorphism from $M_1$ into $P_1$ such that
$\Phi\circ\beta =\beta '$, and $\Psi$ a normal involutive homomorphism from $M_2$ into $P_2$ such that
$\Psi\circ\gamma=\gamma'$, it is possible then to define a normal involutive homomorphism (the proof
given in ([S1] 1.2.4) in the case when $N$ is abelian can be extended without modification in the
general case) :
\[\Phi\underset{N}{_\beta *_\gamma}\Psi 
: M_1\underset{N}{_\beta
*_\gamma}M_2\mapsto P_1\underset{N}{_{\beta '}*_{\gamma '}}P_2\]
In the case when $_\gamma\mathcal{K}_\epsilon$ is a $N-P^o$ bimodule as explained in \ref{rel} and
$_\zeta\mathcal{L}$ a $P$-module, if
$\gamma (N)\subset M_2$ and $\epsilon (P)\subset M_2$, and if $\zeta (P)\subset M_3$, where $M_3$ is
a von Neumann algebra on $\mathcal{L}$, it is possible to consider then $(M_1\underset{N}{_\beta
*_\gamma}M_2)\underset{P}{_\epsilon *_\zeta}M_3$ and $M_1\underset{N}{_\beta
*_\gamma}(M_2\underset{P}{_\epsilon *_\zeta}M_3)$. The coassociativity property for relative tensor
products leads then to the isomorphism of these von Neumann algebra we shall write now 
$M_1\underset{N}{_\beta
*_\gamma}M_2\underset{P}{_\epsilon *_\zeta}M_3$ without parenthesis.
\newline
If $M_1$ and $M_2$ are finite-dimensional, the fiber product $M_1\underset{N}{_\beta *_ \gamma}M_2$
can be identified to a reduced algebra of $M_1\otimes M_2$ (reduced by a projector which belongs to
$\beta (N)\otimes \gamma (N)$). ([EV] 2.4)

\subsection{Hopf-bimodules}
\label{Hbimod}
A quadruplet $(N, M, \alpha, \beta, \Gamma)$ will be called a Hopf-bimodule, following ([Val2], [EV] 6.5), if
$N$,
$M$ are von Neumann algebras, $\alpha$ a faithful non-degenerate representation of $N$ into $M$, $\beta$ a
faithful non-degenerate anti-representation of
$N$ into $M$, with commuting ranges, and $\Gamma$ an injective involutive homomorphism from $M$
into
$M\underset{N}{_\beta *_\alpha}M$ such that, for all $X$ in $N$ :
\newline
(i) $\Gamma (\beta(X))=1\underset{N}{_\beta\otimes_\alpha}\beta(X)$
\newline
(ii) $\Gamma (\alpha(X))=\alpha(X)\underset{N}{_\beta\otimes_\alpha}1$ 
\newline
(iii) $\Gamma$ satisfies the co-assosiativity relation :
\[(\Gamma \underset{N}{_\beta *_\alpha}id)\Gamma =(id \underset{N}{_\beta *_\alpha}\Gamma)\Gamma\]
This last formula makes sense, thanks to the two preceeding ones and
\ref{fiber}\vspace{5mm}.\newline
If $(N, M, \alpha, \beta, \Gamma)$ is a Hopf-bimodule, it is clear that
$(N^o, M, \beta, \alpha,
\varsigma_N\circ\Gamma)$ is another Hopf-bimodule, we shall call the symmetrized of the first
one. (Recall that $\varsigma_N\circ\Gamma$ is a homomorphism from $M$ to
$M\underset{N^o}{_r*_s}M$).
\newline
If $N$ is abelian, $\alpha=\beta$, $\Gamma=\varsigma_N\circ\Gamma$, then the quadruplet $(N, M, \alpha, \alpha,
\Gamma)$ is equal to its symmetrized Hopf-bimodule, and we shall say that it is a symmetric
Hopf-bimodule\vspace{5mm}.\newline
Let $\mathcal G$ be a groupo\"{\i}d, with $\mathcal G^{(0)}$ as its set of units, and let us denote
by $r$ and $s$ the range and source applications from $\mathcal G$ to $\mathcal G^{(0)}$, given by
$xx^{-1}=r(x)$ and $x^{-1}x=s(x)$. As usual, we shall denote by $\mathcal G^{(2)}$ (or $\mathcal
G^{(2)}_{s,r}$) the set of composable elements, i.e. 
\[\mathcal G^{(2)}=\{(x,y)\in \mathcal G^2; s(x)=r(y)\}\]
In [Y1] and [Val2] were associated to a measured groupo\"{\i}d $\mathcal G$, equipped with a Haar system $(\lambda^u)_{u\in \mathcal G ^{(0)}}$ and a quasi-invariant measure $\mu$ on $\mathcal G ^{(0)}$ (see [R1],
[R2], [C2] II.5 and [AR] for more details, precise definitions and examples of groupo\"{\i}ds) two
Hopf-bimodules : 
\newline
The first one is $(L^\infty (\mathcal G^{(0)}, \mu), L^\infty (\mathcal G, \nu), r_{\mathcal G}, s_{\mathcal G}, \Gamma_{\mathcal
G})$, where $\nu$ is the measure constructed on $\mathcal G $ using $\mu$ and the Haar system $(\lambda^u)_{u\in \mathcal G ^{(0)}}$, where we define $r_{\mathcal G}$ and $s_{\mathcal G}$ by writing , for $g$ in $L^\infty (\mathcal G^{(0)})$ :
\[r_{\mathcal G}(g)=g\circ r\]
\[s_{\mathcal G}(g)=g\circ s\]
 and where
$\Gamma_{\mathcal G}(f)$, for $f$ in $L^\infty (\mathcal G)$, is the function defined on $\mathcal G^{(2)}$ by $(s,t)\mapsto f(st)$;
$\Gamma_{\mathcal G}$ is then an involutive homomorphism from $L^\infty (\mathcal G)$ into $L^\infty
(\mathcal G^2_{s,r})$ (which can be identified to
$L^\infty (\mathcal G){_s*_r}L^\infty (\mathcal G)$).
\newline
The second one is symmetric; it is $(L^\infty (\mathcal G^{(0)}), \mathcal L(\mathcal G), r_{\mathcal G}, r_{\mathcal G},
\widehat{\Gamma_{\mathcal G}})$, where
$\mathcal L(\mathcal G)$ is the von Neumann algebra generated by the convolution algebra associated to the
groupo\"{\i}d
$\mathcal G$, and $\widehat{\Gamma_{\mathcal G}}$ has been defined in [Y1] and
[Val2]\vspace{5mm}.\newline
If $(N,M,r,s,\Gamma)$ be a Hopf-bimodule with a finite-dimensional algebra $M$, then, the
identification of $M\underset{N}{_\beta*_\alpha}M$ with a reduced algebra $(M\otimes M)_e$ (\ref{fiber})
leads to an injective homomorphism $\widetilde{\Gamma}$ from $M$ to $M\otimes M$ such that
$\widetilde{\Gamma}(1)=e\not= 1$ and $(\widetilde{\Gamma}\otimes
id)\widetilde{\Gamma}=(id\otimes\widetilde{\Gamma})\widetilde{\Gamma}$ ([EV] 6.5). Then $(M,
\widetilde{\Gamma})$ is a weak Hopf
$\mathbb{C}^*$-algebra in the sense of ([BSz1], [BSz2], [Sz]).

\subsection{Slice maps [E3]}
\label{slice}
Let $A$ be in $M_1\underset{N}{_\beta *_\gamma}M_2$, $\psi$ a normal faithful semi-finite weight on $N$, $\mathcal{H}$ an Hilbert space on which $M_1$ is acting, $\mathcal{K}$ an Hilbert space on which $M_2$ is acting, and let $\xi_1$, $\xi_2$ be in
$D(\mathcal{H}_\beta,\psi^o)$; let us define :
\[(\omega_{\xi_1, \xi_2}\underset{\psi}{_\beta*_\gamma}id)(A)=(\lambda^{\beta, \gamma}_{\xi_2})^*A\lambda^{\beta, \gamma}_{\xi_1}\]
We define this way a $(\omega_{\xi_1, \xi_2}\underset{\psi}{_\beta*_\gamma}id)(A)$ as a bounded operator on $\mathcal{K}$,
which belongs to $M_2$, such that :
\[((\omega_{\xi_1, \xi_2}\underset{\psi}{_\beta*_\gamma}id)(A)\eta_1|\eta_2)=
(A(\xi_1\underset{\psi}{_\beta\otimes_\gamma}\eta_1)|
\xi_2\underset{\psi}{_\beta\otimes_\gamma}\eta_2)\]
One should note that $(\omega_{\xi_1, \xi_2}\underset{\psi}{_\beta*_\gamma}id)(1)=\gamma (<\xi_1, \xi_2 >_{\beta, \psi^o})$. 
\newline
Let us define the same way, for any $\eta_1$, $\eta_2$ in
$D(_\gamma\mathcal{K}, \psi)$:
\[(id\underset{\psi}{_\beta*_\gamma}\omega_{\eta_1, \eta_2})(A)=(\rho^{\beta, \gamma}_{\eta_2})^*A\rho^{\beta, \gamma}_{\eta_1}\]
which belongs to $M_1$. 
\newline
We therefore have a Fubini formula for these slice maps : for any $\xi_1$, $\xi_2$ in
$D(\mathcal{H}_\beta,\psi^o)$, $\eta_1$, $\eta_2$ in $D(_\gamma\mathcal{K}, \psi)$, we have :
\[<(\omega_{\xi_1, \xi_2}\underset{\psi}{_\beta*_\gamma}id)(A), \omega_{\eta_1, \eta_2}>=<(id\underset{\psi}{_\beta*_\gamma}\omega_{\eta_1,
\eta_2})(A),\omega_{\xi_1, \xi_2}>\]
Let $\phi_1$ be a normal semi-finite weight on
$M_1^+$, and $A$ be a positive element of the fiber product
$M_1\underset{N}{_\beta*_\gamma}M_2$, then we may define an element of the extended positive part
of $M_2$, denoted
$(\phi_1\underset{\psi}{_\beta*_\gamma}id)(A)$, such that, for all $\eta$ in $D(_\gamma L^2(M_2), \psi)$, we have :
\[\|(\phi_1\underset{\psi}{_\beta*_\gamma}id)(A)^{1/2}\eta\|^2=\phi_1(id\underset{\psi}{_\beta*_\gamma}\omega_\eta)(A)\]
Moreover, then, if $\phi_2$ is a normal semi-finite weight on $M_2^+$, we have :
\[\phi_2(\phi_1\underset{\psi}{_\beta*_\gamma}id)(A)=\phi_1(id\underset{\psi}{_\beta*_\gamma}\phi_2)(A)\]
and if $\omega_i$ be in $M_{1*}$ such that $\phi_1=sup_i\omega_i$, we
have $(\phi_1\underset{\psi}{_\beta*_\gamma}id)(A)=sup_i(\omega_i\underset{\psi}{_\beta*_\gamma}id)(A)$.
\newline
Let now $P_1$ be a von Neuman algebra such that :
\[\beta(N)\subset P_1\subset M_1\]
and let $\Phi_i$ ($i=1,2$)
be a normal faithful semi-finite operator valued weight from $M_i$ to $P_i$; for any positive
operator $A$ in the fiber product
$M_1\underset{N}{_\beta*_\gamma}M_2$, there exists an element $(\Phi_1\underset{\psi}{_\beta*_\gamma}id)(A)$
of the extended positive part
of $P_1\underset{N}{_\beta*_\gamma}M_2$, such that ([E3], 3.5), for all
$\eta$ in $D(_\gamma L^2(M_2), \psi)$, and $\xi$ in $D(L^2(P_1)_\beta, \psi^o)$, we have :
\[\|(\Phi_1\underset{\psi}{_\beta*_\gamma}id)(A)^{1/2}(\xi\underset{\psi}{_\beta\otimes_\gamma}\eta)\|^2=
\|\Phi_1(id\underset{\psi}{_\beta*_\gamma}\omega_\eta)(A)^{1/2}\xi\|^2\]
If $\phi$ is a normal semi-finite weight on $P$, we have :
\[(\phi\circ\Phi_1\underset{\psi}{_\beta*_\gamma}id)(A)=(\phi\underset{\psi}{_\beta*_\gamma}id)
(\Phi_1\underset{\psi}{_\beta*_\gamma}id)(A)\]
We define the same way an element $(id\underset{\psi}{_\beta*_\gamma}\Phi_2)(A)$ of the extended positive part
of
$M_1\underset{N}{_\gamma*_\beta}P_2$, and we have :
\[(id\underset{\psi}{_\beta*_\gamma}\Phi_2)((\Phi_1\underset{\psi}{_\beta*_\gamma}id)(A))=
(\Phi_1\underset{\psi}{_\beta*_\gamma}id)((id\underset{\psi}{_\beta*_\gamma}\Phi_2)(A))\]
 Considering now an element $x$ of $M_1{}_\beta\underset{\psi}{*}{}_\pi \pi(N)$, which can be identified 
(\ref{fiber}) to $M_1\cap\beta(N)'$, we get that, for $e$ in $\gN_\psi$, we have  \[(id_\beta\underset{\psi}{*}{}_\pi\omega_{J_\psi \Lambda_{\psi}(e)})(x)=\beta(ee^*)x\]
 Therefore, by increasing limits, we get that $(id_\beta\underset{\psi}{*}{}_\pi\psi)$ is the injection of $M_1\cap\beta(N)'$ into $M_1$.  More precisely, if $x$ belongs to $M_1\cap\beta(N)'$, we have :
 \[(id_\beta\underset{\psi}{*}{}_\pi\psi)(x{}_\beta\underset{\psi}{\otimes}{}_\pi 1)=x\]
 \newline
 Therefore, if $\Phi_2$ is a normal faithful semi-finite operator-valued weight from $M_2$ onto $\gamma(N)$, we get that, for all $A$ positive in $M_1\underset{N}{_\beta*_\gamma}M_2$, we have :
 \[(id_\beta\underset{\psi}{*}{}_\gamma\psi\circ\Phi_2)(A){}_\beta\underset{\psi}{\otimes}{}_\gamma 1=
 (id_\beta\underset{\psi}{*}{}_\gamma\Phi_2)(A)\]

With the notations of \ref{spatial}, let $(e_i)_{i\in I}$ be a $(\beta, \psi^o)$-orthogonal basis of $\mathcal H$; using the fact (\ref{rel}) that, for all $\eta$ in $\mathcal K$, we have :\[e_i\underset{\psi}{_\beta\otimes_\gamma}\eta=e_i\underset{\psi}{_\beta\otimes_\gamma}\gamma(<e_i, e_i>_{\beta, \psi^o})\eta\]
we get that, for all $X$ in $M_1\underset{N}{_\beta*_\gamma}M_2$, $\xi$ in $D(\mathcal H_\beta, \psi^o)$, we have 
\[(\omega_{\xi, e_i}\underset{\psi}{_\beta *_\gamma}id)(X)=\gamma(<e_i, e_i>_{\beta, \psi^o})(\omega_{\xi, e_i}\underset{\psi}{_\beta *_\gamma}id)(X)\]

\subsection{Notations} 
\label{aut}
Let $M$ be a von Neuman algebra, and $\alpha$ an action from a locally compact group $G$ on $M$, i.e. a homomorphism from $G$ into $Aut M$, such that, for all $x\in M$, the function $g\mapsto \alpha_g(x)$ is $\sigma$-weakly continuous. Let us denote by $\mathbb{C}^*(\alpha)$ the set of elements $x$ of $M$, such that this function $t\mapsto \alpha_g(x)$ is norm continuous. It is ([P], 7.5.1) a sub-$\mathbb{C}^*$-algebra of $M$, invariant under the $\alpha_g$, generated by the elements ($x\in N$,$f\in L^1(G)$) :
\[\alpha_f(x)=\int_{\mathbb{R}}f(s)\alpha_s(x)ds\]
More precisely, we get that, for any $x$ in $M$, $\alpha_f(x)$ is $\sigma$-weakly converging to $x$ when $f$ goes in an approximate unit of $L^1(\mathbb{G})$, which proves that $\mathbb{C}^*(\alpha)$ is $\sigma$-weakly dense in $M$, and that $x\in M$ belongs to $\mathbb{C}^*(\alpha)$ if and only if this file is norm converging. 
\newline
If $\alpha_t$ and $\gamma_s$ are two one-parameter automorphism groups of $M$, such that, for all $s$, $t$ in $\mathbb{R}$, we have $\alpha_t\circ\gamma_s=\gamma_s\circ\alpha_t$, by considering the action of $\mathbb{R}^2$ given by $(s,t)\mapsto \gamma_s\circ\alpha_t$, we obtain a dense sub-$\mathbb{C}^*$-algebra of $M$, on which both $\alpha$ and $\gamma$ are norm continuous, we shall denote $\mathbb{C}^*(\alpha, \gamma)$.

\subsection{Weights on $\mathbb{C}^*$-algebras}
\label{C*}
Let $A$ be a $\mathbb{C}^*$-algebra, and $\varphi$ a lower semi-continuous, densely defined non zero weight on $A$ ([Co1]). We shall use all classical notations, and, in particular, we shall denote $(H_\varphi, \Lambda_\varphi, \pi_\varphi)$ the GNS construction for $\varphi$; $\pi_\varphi$ is faithful if, and only if, $\varphi$ is faithful; let us denote $M=\pi_\varphi(A)''$ and $\overline{\varphi}$ the semi-finite normal weight on $M^+$, constructed by  ([B], cor. 9), which verify $\overline{\varphi}\circ\pi_\varphi = \varphi$. 
\newline
Let us recall that if the $\mathbb{C}^*$-algebra $A$ is unital, any densely defined weight $\varphi$ is everywhere defined, and therefore finite. 
\newline
Following [Co1], we shall say that $\varphi$ is KMS if there exists a norm-continuous one parameter group of automorphisms $\sigma_t$ of $A$ such that, for all $t\in \mathbb{R}$, $\varphi=\varphi\circ\sigma_t$, and such that $\varphi$ verifies the KMS conditions with respect to $\sigma$. (For an equivalent definition of these conditions, see [KV1], 1.3). One can find in ([KV1], 1.35) the proof that every KMS weight extends to a faithful extension $\overline{\varphi}$ on $M^+$, and that we have then $\pi_\varphi\circ\sigma_t=\sigma^{\overline{\varphi}}_t\circ\pi_\varphi$, where $\sigma^{\overline{\varphi}}_t$ is the modular automorphism group of $M$ given by the Tomita-Takesaki theory of the faithful semi-finite weight $\overline{\varphi}$ on $M$. This leads easily to the uniqueness of the one-parameter group $\sigma_t$, which we shall emphasize by writing it $\sigma_t^\varphi$.
\newline
Moreover, it is well known that the set of elements $x$ in $A$ such that the function $t\mapsto\sigma_t^\varphi(x)$ extends to an analytic function in $A$ is a dense involutive subalgebra of $A$ (see for instance [Val1] 0.3.2 and 0.3.4). 
\newline
Let $\varphi$ be a lower semi-continuous, densely defined, faithful weight on $A$; let $\alpha$ be a representation of $A$ on a Hilbert space $_\alpha H$. A vector $\xi$ in $_\alpha H$ will be said to be bounded with respect to $\varphi$ if there exists a positive $C$ such that, for all $x\in\gN_\varphi$, we have :
\[\|\alpha (x)\xi\|^2\leq C\varphi(x^*x)\]
Following ([Co2], 1.7), we shall say that $\alpha$ is square-integrable with respect to $\varphi$ if the set of bounded vectors is total in $_\alpha H$. Then ([Co2], 1.8), for any representation $\alpha$ of $A$, which is square-integrable with respect to $\varphi$, there exists a unique normal representation $\overline{\alpha}$ of $M$ on $_\alpha H$ such that $\overline{\alpha}\circ\pi_\varphi = \alpha$. 
\newline
It is proved in ([KV1], 1.34) that the representation $\pi_\varphi$ itself is square-integrable with respect to $\varphi$ if and only if the weight $\overline{\varphi}$ is faithful. 
\newline

Let now $\psi$ be a normal semi-finite faithful weight on $N$. We shall write $\mathbb{C}^*(\psi)$ the sub-$\mathbb{C}^*$-algebra of $\mathbb{C}^*(\sigma^\psi)$ generated by elements $\sigma^\psi_f(x)$, with $f\in L^1(\mathbb{R})$ and $x\in\gM_\psi$. The weak closure of $\mathbb{C}^*(\psi)$ contains $\gM_\psi$, and, therefore, $\mathbb{C}^*(\psi)$ is weakly dense in $N$; moreover, it is straightforward to see that the restriction of $\psi$ to that $\mathbb{C}^*$-algebra is densely defined, lower semi-continuous and KMS. If $1\in\mathbb{C}^*(\psi)$, then the restriction of $\psi$ to  $\mathbb{C}^*(\psi)$ is finite, so $\psi(1)<\infty$ and $\mathbb{C}^*(\psi)=\mathbb{C}^*(\sigma^\psi)$. If $\psi$ is a trace, then $\mathbb{C}^*(\psi)$ is the norm closure of $\gM_\psi$, and $M(\mathbb{C}^*(\psi))=N$. 
\newline
If $\gamma_t$ is one-parameter group of $N$, such that $\psi\circ \gamma_t=\psi$, for all $t\in \mathbb{R}$, we may as well define the $\mathbb{C}^*$-algebra $\mathbb{C}^*(\psi, \gamma)$ generated by all elements :
\[\int_{\mathbb{R}^2}f(s)g(t)\sigma_s^\psi\circ\gamma_t(x)dsdt\]
where $f,g$ belong to $L^1(\mathbb{R})$, and $x$ belongs to $\gM_\psi$; this $\mathbb{C}^*$-algebra $\mathbb{C}^*(\psi, \gamma)$ is weakly dense in $N$, invariant under $\gamma$, the restriction of $\psi$ to this $\mathbb{C}^*$-algebra is densely defined, lower semi-continuous and KMS, and the restriction of $\gamma$ to this $\mathbb{C}^*$-algebra is norm-continuous. If $\psi$ is a trace, we have $M(\mathbb{C}^*(\psi, \gamma))=\mathbb{C}^*(\gamma)$. 
\newline
Let us recall another result borrowed from Baaj ([B], th. 8) ; if $\varphi$ is a lower semi-continuous weight defined on a hereditary cone $F$ of a $\mathbb{C}^*$-algebra $A$, such that the subcone $F'$ made of elements $X\in F$ such that $\varphi (X)<\infty$ is norm dense in $A$, then $\varphi$ has an extension $\tilde{\varphi}$ on $A^+$ which is lower semi-continuous and semi-finite.


\section{Pseudo-multiplicative unitary}
\label{pmu}
In this chapter, we recall (\ref{defmult}) the definition of a pseudo-multiplicative unitary, give the fundamental example given by groupo\"{\i}ds (\ref{gd}), and construct the algebras and the Hopf-bimodules "generated by the left (resp. right) leg" of a pseudo-multiplicative unitary (\ref{AW}). In \ref{manageable}, we recall the definition of a "manageable" pseudo-multiplicative unitary, and we obtain in \ref{propmanag} a result about the norm closed algebra generated by the left leg of a manageable pseudo-multiplicative unitary. In particular, if the basis is abelian, this algebra is a $\mathbb{C}^*$-algebra. (\ref{propmanag}(iv)).
\subsection{Definition}
\label{defmult}
Let $N$ be a von Neumann algebra; let
$\gH$ be a Hilbert space on which $N$ has a non-degenerate normal representation $\alpha$ and two
non-degenerate normal anti-representations $\hat{\beta}$ and $\beta$. These 3 applications
are supposed to be injective, and to commute two by two.  Let $\nu$ be a normal semi-finite faithful weight on
$N$; we can therefore construct the Hilbert spaces
$\gH\underset{\nu}{_\beta\otimes_\alpha}\gH$ and
$\gH\underset{\nu^o}{_\alpha\otimes_{\hat{\beta}}}\gH$. A unitary $W$ from
$\gH\underset{\nu}{_\beta\otimes_\alpha}\gH$ onto
$\gH\underset{\nu^o}{_\alpha\otimes_{\hat{\beta}}}\gH$.
will be called a pseudo-multiplicative unitary over the basis $N$, with respect to the
representation $\alpha$, and the anti-representations $\hat{\beta}$ and $\beta$ (we shall write it is an $(\alpha, \hat{\beta}, \beta)$-pseudo-multiplicative unitary), if :
\newline
(i) $W$ intertwines $\alpha$, $\hat{\beta}$, $\beta$  in the following way :
\[W(\alpha
(X)\underset{N}{_\beta\otimes_\alpha}1)=
(1\underset{N^o}{_\alpha\otimes_{\hat{\beta}}}\alpha(X))W\]
\[W(1\underset{N}{_\beta\otimes_\alpha}\beta
(X))=(1\underset{N^o}{_\alpha\otimes_{\hat{\beta}}}\beta (X))W\]
\[W(\hat{\beta}(X) \underset{N}{_\beta\otimes_\alpha}1)=
(\hat{\beta}(X)\underset{N^o}{_\alpha\otimes_{\hat{\beta}}}1)W\]
\[W(1\underset{N}{_\beta\otimes_\alpha}\hat{\beta}(X))=
(\beta(X)\underset{N^o}{_\alpha\otimes_{\hat{\beta}}}1)W\]
(ii) The operator satisfies :
\begin{multline*}
(1_{\gH}\underset{N^o}{_\alpha\otimes_{\hat{\beta}}}W)
(W\underset{N}{_\beta\otimes_\alpha}1_{\gH})=\\
=(W\underset{N^o}{_\alpha\otimes_{\hat{\beta}}}1_{\gH})
(\sigma_{\nu^o}\underset{N^o}{_\alpha\otimes_{\hat{\beta}}}1_{\gH})
(1_{\gH}\underset{N^o}{_\alpha\otimes_{\hat{\beta}}}W)
\sigma_{2\nu}
(1_{\gH}\underset{N}{_\beta\otimes_\alpha}\sigma_{\nu^o})
(1_{\gH}\underset{N}{_\beta\otimes_\alpha}W)
\end{multline*}
In that formula, the first $\sigma_{\nu^o}$ is the relative flip defined in \ref{rel} from $\gH\underset{\nu^o}{_\alpha\otimes_\beta}\gH$
to $\gH\underset{\nu}{_\beta\otimes_\alpha}\gH$, and the second is the relative flip from
$\gH\underset{\nu^o}{_\alpha\otimes_{\hat{\beta}}}\gH$ to $\gH\underset{\nu}{_{\hat{\beta}}\otimes_\alpha}\gH$; while $\sigma_{2\nu}$ is the relative flip
from $\gH\underset{\nu}{_\beta\otimes_\alpha}\gH\underset{\nu}{_{\hat{\beta}}\otimes_\alpha}\gH$ to
$\gH\underset{\nu^o}{_\alpha\otimes_{\hat{\beta}}}(\gH\underset{\nu}{_\beta\otimes_\alpha}\gH)$. The index $2$ is written to recall that the
flip "turns" around the second relative tensor product, and, in such a formula, the parenthesis are written to recall that, in such a situation,
associativity rules does not occur because the anti-representation $\hat{\beta}$ is here acting in the second leg of
$\gH\underset{\nu}{_\beta\otimes_\alpha}\gH$. 
\newline
All the properties supposed in (i) allow us to write such a formula, which will be called the
"pentagonal relation". 
\newline
If we start from another normal semi-finite faithful weight $\nu'$ on $N$, we may define, using \ref{rel}, another unitary $W^{\nu'}=U^{\nu^o,
\nu^{'o}}_{\alpha, {\hat{\beta}}}WU^{\nu', \nu}_{\beta, \alpha}$ from $\gH\underset{\nu'}{_\beta\otimes_\alpha}\gH$ onto
$\gH\underset{\nu^{'o}}{_\alpha\otimes_{\hat{\beta}}}\gH$. The formulae which link these isomorphims between relative product Hilbert spaces and the
relative flips allow us to check that this operator $W^{\nu'}$ is also pseudo-multiplicative; which can be resumed in saying that a
pseudo-multiplicative unitary does not depend on the choice of the weight on $N$. 
\newline
If $W$ is an $(\alpha, \hat{\beta}, \beta)$-pseudo-multiplicative unitary, then the unitary $\sigma_\nu W^*\sigma_\nu$ from $\gH\underset{\nu}{_{\hat{\beta}}\otimes_\alpha}\gH$ to $\gH\underset{\nu^o}{_\alpha\otimes_\beta}\gH$ is an $(\alpha, \beta, \hat{\beta})$-pseudo-multiplicative unitary, called the dual of $W$.

\subsection{Algebras and Hopf-bimodules associated to a pseudo-mul-tiplicative unitary}
\label{AW}
For $\xi_2$ in $D(_\alpha\gH, \nu)$, $\eta_2$ in $D(\gH_{\hat{\beta}}, \nu^o)$, the operator $(\rho_{\eta_2}^{\alpha,
\hat{\beta}})^*W\rho_{\xi_2}^{\beta, \alpha}$ will be written $(id*\omega_{\xi_2, \eta_2})(W)$; we have, therefore, for all
$\xi_1$, $\eta_1$ in $\gH$ :
\[((id*\omega_{\xi_2, \eta_2})(W)\xi_1|\eta_1)=(W(\xi_1\underset{\nu}{_\beta\otimes_\alpha}\xi_2)|
\eta_1\underset{\nu^o}{_\alpha\otimes_{\hat{\beta}}}\eta_2)\]
and, using the intertwining property of $W$ with $\hat{\beta}$, we easily get that $(id*\omega_{\xi_2, \eta_2})(W)$ belongs
to $\hat{\beta} (N)'$. 
\newline
If $x$ belongs to $N$, we have $(id*\omega_{\xi_2, \eta_2})(W)\alpha (x)=(id*\omega_{\xi_2, \alpha(x^*)\eta_2})(W)$, and $\beta(x)(id*\omega_{\xi_2, \eta_2})(W)=(id*\omega_{\hat{\beta}(x)\xi_2, \eta_2})(W)$. 
\newline
If $y$ belongs to $N$ and is analytical with respect to $\nu$, then, we have 
\[\alpha(y)(id*\omega_{\xi_2, \eta_2})(W)=(id*\omega_{\xi_2,\hat{\beta}(\sigma^\nu_{i/2}(y^*)) \eta_2})(W)\]
\[(id*\omega_{\xi_2, \eta_2})(W)\beta(y)=(id*\omega_{\alpha(\sigma^\nu_{i/2}(y))\xi_2, \eta_2})(W)\]
\newline
If $\xi$ belongs to $D(_\alpha\gH, \nu)\cap D(\gH_{\hat{\beta}}, \nu^o)$, we shall write $(id*\omega_\xi)(W)$ instead of
$(id*\omega_{\xi, \xi})(W)$.
\newline
We shall write $A_n(W)$ (resp. $A_w(W)$ the norm (resp. weak) closure of the linear span of these operators, which are right $\alpha(N)$-modules and left $\beta(N)$-modules. If we define $B$ as the norm closure of the subalgebra of elements in $N$ which are analytical with respect to $\nu$, we get also that $A_n(W)$ is a left $\alpha (B)$-module, and a right $\beta(B)$-module. Applying ([E2] 3.6), we get that $A_n(W)^*$, $A_n(W)$, $A_w(W)^*$ and $A_w(W)$ are non-degenerate algebras. One should note that the notations of ([E2]) had been changed in order to fit with Lesieur's notations. 
We shall write $\mathcal A(W)$ the von Neumann algebra generated by  $A_n(W)$ (or $A_w(W)$) .
We then have $\mathcal A(W)\subset\hat{\beta}(N)'$.
\newline
For $\xi_1$ in $D(\gH_\beta,\nu^o)$, $\eta_1$ in $D(_\alpha\gH, \nu)$, the operator $(\lambda_{\eta_1}^{\alpha,
\hat{\beta}})^*W\lambda_{\xi_1}^{\beta, \alpha}$ will be written $(\omega_{\xi_1,
\eta_1}*id)(W)$ for ; we have,
therefore, for all
$\xi_2$,
$\eta_2$ in
$\gH$ :
\[((\omega_{\xi_1,\eta_1}*id)(W)\xi_2|\eta_2)=(W(\xi_1\underset{\nu}{_\beta\otimes_\alpha}\xi_2)|
\eta_1\underset{\nu^o}{_\alpha\otimes_{\hat{\beta}}}\eta_2)\]
and, using the intertwining property of $W$ with $\beta$, we easily get that $(\omega_{\xi_1,
\eta_1}*id)(W)$ belongs to $\beta(N)'$. If $\xi$ belongs to $D(\gH_\beta, \nu^o)\cap D(_\alpha\gH, \nu)$, we shall write
$(\omega_\xi*id)(W)$ instead of $(\omega_{\xi, \xi}*id)(W)$. 
\newline
We shall write $\widehat{A_n(W)}$ the norm closure of the linear span of these operators. Then, we had obtained ([E2], 3.6) that this norm closed subspace is a non degenerate algebra (and so is $\widehat{A_n(W)}^*$ also). 
\newline
We shall write $\widehat{A_w(W)}$ the weak closure of the linear span of these operators. It is clear that this weakly closed subspace is a non degenarate algebra; following ([EV] 6.1 and 6.5), we shall write $\widehat{\mathcal A(W})$ the von Neumann algebra generated by  $\widehat{A_n(W)}$ or $\widehat{A_w(W)}$. We then have $\widehat{\mathcal A(W)}\subset\beta(N)'$. 
\newline
In ([EV] 6.3 and 6.5), using the pentagonal equation, we got
that
$(N,\mathcal A(W),\alpha,\beta,\Gamma)$, and
$(N^o,\widehat{\mathcal A(W)}, \hat{\beta}, \alpha, \widehat{\Gamma})$ are Hopf-bimodules, where $\Gamma$ and
$\widehat{\Gamma}$ are defined, for any $x$ in $\mathcal A(W)$ and $y$ in $\widehat{\mathcal
A(W)}$, by :
\[\Gamma(x)=W^*(1\underset{N^o}{_\alpha\otimes_{\hat{\beta}}}x)W\]
\[\widehat{\Gamma}(y)=W(y\underset{N}{_\beta\otimes_\alpha}1)W^*\]
In ([EV] 6.1(iv)), we had obtained that $x$ in $\mathcal L(\gH)$ belongs to $\mathcal A(W)'$ if and only if $x$ belongs to $\alpha(N)'\cap
\beta(N)'$ and verifies \[(x\underset{N^o}{_\alpha\otimes_{\hat{\beta}}}1)W=W(x\underset{N}{_\beta\otimes_\alpha}1)\]
We obtain the same way
that $y$ in $\mathcal L(\gH)$ belongs to $\widehat{\mathcal A(W)}'$ if and only if $y$ belongs to $\alpha(N)'\cap
\hat{\beta}(N)'$ and verify $(1\underset{N^o}{_\alpha\otimes_{\hat{\beta}}}y)W=W(1\underset{N}{_\beta\otimes_\alpha}y)$.  \newline
Moreover, we get that $\alpha(N)\subset\mathcal A\cap\widehat{\mathcal A}$, $\beta(N)\subset\mathcal A$,
$\hat{\beta}(N)\subset\widehat{\mathcal A}$, and, for all $x$ in $N$ :
\[\Gamma (\alpha (x))=\alpha (x)\underset{N}{_\beta\otimes_\alpha}1\]
\[\Gamma (\beta (x))=1\underset{N}{_\beta\otimes_\alpha}\beta (x)\]
\[\widehat{\Gamma}(\alpha(x))=1\underset{N^o}{_\alpha\otimes_{\hat{\beta}}}\alpha (x)\]
\[\widehat{\Gamma}(\hat{\beta}(x))=\hat{\beta}(x)\underset{N^o}{_\alpha\otimes_{\hat{\beta}}}1\]
Following ([E2], 3.7) If $\eta_1$, $\xi_2$ are in $D(_\alpha\gH, \nu)$, let us write $(id*\omega_{\xi_2, \eta_1})(\sigma_{\nu^o}W)$ for $(\lambda_{\eta_1}^{\alpha, \hat{\beta}})^*W\rho_{\xi_2}^{\beta, \alpha}$; we have, therefore, for all $\xi_1$ and $\eta_2$ in $\gH$ :
\[(id*\omega_{\xi_2, \eta_1})(\sigma_{\nu^o}W)\xi_1|\eta_2)=(W(\xi_1\underset{\nu}{_\beta\otimes_\alpha}\xi_2)|\eta_1\underset{\nu^o}{_\alpha\otimes_{\hat{\beta}}}\eta_2)\]
Using the intertwining property of $W$ with $\alpha$, we get that it belongs to $\alpha(N)'$; we write $C_w(W)$ for the weak closure of the linear span of these operators, and we have $C_w(W)\subset \alpha(N)'$. It had been proved in ([E2], 3.10) that $C_w(W)$ is a non degenerate algebra; following ([E2] 4.1), we shall say that $W$ is weakly regular if $C_w(W)=\alpha (N)'$.

\subsection{Fundamental example}
\label{gd}
Let $\mathcal G$ be a measured groupo\"{\i}d, with $\mathcal G^{(0)}$ as space
of units, and $r$ and $s$ the range and source functions from $\mathcal G$ to $\mathcal G^{(0)}$, with a Haar system $(\lambda^u)_{u\in \mathcal G^{(0)}}$ and a quasi-invariant measure $\mu$ on $\mathcal G^{(0)}$. Let us write $\nu$ the associated measure on $\mathcal G$. Let us
note :
\[\mathcal G^2_{r,r}=\{(x,y)\in \mathcal G^2, r(x)=r(y)\}\]
Then, it has been shown [Val2] that the formula $W_{\mathcal G}f(x,y)=f(x,x^{-1}y)$, where $x$, $y$ are
in
$\mathcal G$, such that $r(y)=r(x)$, and $f$ belongs to $L^2(\mathcal G^{(2)})$ (with respect to an
appropriate measure, constructed from $\lambda^u$ and $\mu$), is a unitary from $L^2(\mathcal G^{(2)})$ to $L^2(\mathcal G^2_{r,r})$ (with respect also to another
appropriate measure, constructed from $\lambda^u$ and $\mu$). 
\newline
Let us define $r_{\mathcal G}$ and $s_{\mathcal G}$ from
$L^\infty (\mathcal G^{(0)})$ to $L^\infty (\mathcal G)$ (and then considered as representations on $\mathcal L(L^2(\mathcal
G))$, for any
$f$ in
$L^\infty (\mathcal G^{(0)})$, by
$r_{\mathcal G}(f)=f\circ r$ and $s_{\mathcal G}(f)=f\circ s$.
\newline
We shall identify ([Y1], 3.2.2) the Hilbert space $L^2(\mathcal G^{(2)})$ with the relative Hilbert tensor product $L^2(\mathcal G, \nu)\underset{L^{\infty}(\mathcal G^{(0)}, \mu)}{_{s_{\mathcal G}}\otimes_{r_{\mathcal G}}}L^2(\mathcal G, \nu)$, and the Hilbert space $L^2(\mathcal G^2_{r,r})$ with $L^2(\mathcal G, \nu)\underset{L^{\infty}(\mathcal G^{(0)}, \mu)}{_{r_{\mathcal G}}\otimes_{r_{\mathcal G}}}L^2(\mathcal G, \nu)$. Moreover, the unitary $W_{\mathcal G}$ can be then interpreted [Val3] as a pseudo-multiplicative unitary over the basis
$L^\infty (\mathcal G^{(0)})$, with respect to the representation $r_{\mathcal G}$, and anti-representations
$s_{\mathcal G}$ and
$r_{\mathcal G}$ (as here the basis is abelian, the notions of representation and anti-representations are the same, and the commutation property is fulfilled). So, we get that $W_{\mathcal G}$ is a $(r_\mathcal G, s_\mathcal G, r_\mathcal G)$ pseudo-multiplicative unitary. 
\newline
Let us take the notations of \ref{AW}; the von Neumann algebra $\mathcal A(W_{\mathcal G})$ is equal to the von Neumann algebra $L^{\infty}(\mathcal G, \nu)$ ([Val3], 3.2.6 and 3.2.7); using ([Val3]
3.1.1), we get that the Hopf-bimodule homomorphism
$\widehat{\Gamma}$ defined on
$L^{\infty}(\mathcal G, \nu)$ by $W_{\mathcal G}$ is equal to the usual Hopf-bimodule homomorphism $\Gamma_{\mathcal G}$ studied in [Val2], and
recalled in
\ref{Hbimod}.
Moreover, the von Neumann algebra $\widehat{\mathcal A(W_{\mathcal G})}$ is equal to the von Neumann algebra $\mathcal
L(\mathcal G)$ ([Val3], 3.2.6 and 3.2.7); using ([Val3] 3.1.1), we get that the Hopf-bimodule homomorphism $\Gamma$ defined on $\mathcal L(\mathcal
G)$ by
$W_{\mathcal G}$ is the usual Hopf-bimodule homomorphism $\widehat{\Gamma_{\mathcal G}}$ studied in [Y1] and [Val2]. 
\newline
More precisely, we easily get ([Y1], 2.1) that $D(_{r_{\mathcal G}}L^2(\mathcal G, \nu))$ is the set of (equivalent classes of) functions $f$ in $L^2(\mathcal G, \nu)$, such that the function :
\[u\mapsto \int_{G^u}|f(x)|^2d\lambda^u(x)\]
belongs to $L^{\infty}(\mathcal G^{(0)}, \mu)$ (and then $\|R^{r_{\mathcal G}, \mu}(f)\|$ is equal to the square root of the essential supremum of this function). This vector space, with this norm, had been denoted $L^\infty (\mathcal G^{(0)}, L^2(\mathcal G, \lambda))$ in ([AR], 1.2.2). If $f$, $g$ are two functions in $D(_{r_{\mathcal G}}L^2(\mathcal G, \nu))$, we get that $(id*\omega_{f,g})(W_{\mathcal G})$ is the function on $\mathcal G$ defined by :
\[x\mapsto\int_{G^{r(x)}}f(x^{-1}y)\overline{g(y)}d\lambda^{r(x)}(y)\]
So we get that $A_n(W_{\mathcal G})$ is invariant under the adjoint operation, and is therefore a $\mathbb{C}^*$-algebra. It is also, using $r_\mathcal G$ and $s_\mathcal G$, a bimodule over $L^{\infty}(\mathcal G^{(0)})$. 
\newline
Let us suppose now that the groupoid $\mathcal G$ is locally compact in the sense of [R1]; let $\mathcal K(\mathcal G)$ denote the continuous functions on $\mathcal G$, with compact support, and let us suppose that, for all $f$ in $\mathcal K(\mathcal G)$, the function $u\mapsto\int_{\mathcal G}fd\lambda^u$ is continuous. Then $\mathcal K(\mathcal G)\subset D(_{r_{\mathcal G}}L^2(\mathcal G, \nu), \mu)$ ([Val2], 3.2.1); its closure for the norm $f\mapsto\|R^{r_{\mathcal G}, \mu}(f)\|$ had been studied in ([AR], 1.1.1) and shall be denoted $C_0(\mathcal G^{(0)}, L^2(\mathcal G, \lambda))$. If $f$, $g$ belongs to $\mathcal K(\mathcal G)$, it is clear from the formula above that $(id*\omega_{f,g})(W_{\mathcal G})$ belongs also to $\mathcal K(\mathcal G)$; by continuity, for any $f$, $g$ in $C_0(\mathcal G^{(0)}, L^2(\mathcal G, \lambda))$, we get that $(id*\omega_{f,g})(W_{\mathcal G})$ belongs to the algebra $C_0(\mathcal G)$ of continuous functions going to $0$ at infinity. Finally, using Stone-Weirstrass theorem, we shall get that the sub-$\mathbb{C}^*$-algebra generated by these elements is equal to $C_0(\mathcal G)$, and, therefore, that $A_n(W_{\mathcal G})$ contains $C_0(\mathcal G)$ as a sub-$\mathbb{C}^*$-algebra. As $A_n(W_{\mathcal G})$ is a bimodule over $L^{\infty}(\mathcal G^{(0)})$, there is no chance in general that $A_n(W_{\mathcal G})$ should be equal to $C_0(\mathcal G)$. 
\newline
It has been proved in ([E2] 4.8) that $W_\mathcal G$ is weakly regular (in fact was proved a much stronger condition, namely the norm regularity).

\subsection{Lemma}
\label{lem1}
{\it Let $W$ be an $(\alpha, \hat{\beta}, \beta)$-pseudo-multiplicative unitary, $\xi_1$ in $D(\gH_\beta, \nu^o)$, $\xi_2$ in $D(_\alpha \gH, \nu)$, $\eta$ in $\gH$; let $\zeta_i$ in $D(\gH_\beta, \nu^o)$ and $\zeta'_i$ in $\gH$ such that $W^*(\xi_2\underset{\nu^o}{_\alpha\otimes_{\hat{\beta}}}\eta)=\sum_i \zeta_i\underset{\nu}{_\beta\otimes_\alpha}\zeta'_i$; then we have :}
\[\sum_i\alpha(<\zeta_i, \xi_1>_{\beta, \nu^o})\zeta'_i=(\omega_{\xi_1, \xi_2}*id)(W)^*\eta\]

\begin{proof}
Let $\theta$ in $\gH$; we have :
\begin{eqnarray*}
((\omega_{\xi_1, \xi_2}*id)(W)^*\eta|\theta)
&=&(W^*(\xi_2\underset{\nu^o}{_\alpha\otimes_{\hat{\beta}}}\eta)|\xi_1\underset{\nu}{_\beta\otimes_\alpha}\theta)\\
&=&(\sum_i \zeta_i\underset{\nu}{_\beta\otimes_\alpha}\zeta'_i|\xi_1\underset{\nu}{_\beta\otimes_\alpha}\theta)\\
&=&(\sum_i\alpha(<\zeta_i, \xi_1>_{\beta, \nu^o})\zeta'_i|\theta)
\end{eqnarray*}
from which we get the result.  \end{proof}

\subsection{Lemma}
\label{lem2}
{\it Let $W$ be an $(\alpha, \hat{\beta}, \beta)$-pseudo-multiplicative unitary, $\xi_1$, $\zeta_1$ in $D(\gH_\beta, \nu^o)$, $\xi$ in $D(_\alpha\gH, \nu)$ and $\eta_1$, $\eta_2$ in $\gH$. Then, we have :}
\begin{multline*}
(\sigma_{2\nu}
(1_{\gH}\underset{N}{_\beta\otimes_\alpha}\sigma_{\nu^o})
(1_{\gH}\underset{N}{_\beta\otimes_\alpha}W)(\xi_1\underset{\nu}{_\beta\otimes_\alpha}\eta_1\underset{\nu}{_\beta\otimes_\alpha}\xi)|\eta_2\underset{\nu^o}{_\alpha\otimes_{\hat{\beta}}}(\zeta_1\underset{\nu}{_\beta\otimes_\alpha}\zeta_2))=\\
(W(\eta_1\underset{\nu}{_\beta\otimes_\alpha}\xi)|\eta_2\underset{\nu^o}{_\alpha\otimes_{\hat{\beta}}}\alpha(<\zeta_1, \xi_1>_{\beta, \nu^o})\zeta_2)
\end{multline*}
\begin{proof}
The scalar product 
\[(\sigma_{2\nu}
(1_{\gH}\underset{N}{_\beta\otimes_\alpha}\sigma_{\nu^o})
(1_{\gH}\underset{N}{_\beta\otimes_\alpha}W)(\xi_1\underset{\nu}{_\beta\otimes_\alpha}\eta_1\underset{\nu}{_\beta\otimes_\alpha}\xi)|\eta_2\underset{\nu^o}{_\alpha\otimes_{\hat{\beta}}}(\zeta_1\underset{\nu}{_\beta\otimes_\alpha}\zeta_2))\]
is equal to :
\[((1_{\gH}\underset{N}{_\beta\otimes_\alpha}\sigma_{\nu^o})
(1_{\gH}\underset{N}{_\beta\otimes_\alpha}W)(\xi_1\underset{\nu}{_\beta\otimes_\alpha}\eta_1\underset{\nu}{_\beta\otimes_\alpha}\xi)|\zeta_1\underset{\nu}{_\beta\otimes_\alpha}\zeta_2\underset{\nu}{_{\hat{\beta}}\otimes_\alpha}\eta_2)\]
and to :
\[(\xi_1\underset{\nu}{_\beta\otimes_\alpha}W(\eta_1\underset{\nu}{_\beta\otimes_\alpha}\xi)|
\zeta_1\underset{\nu}{_\beta\otimes_\alpha}(\eta_2\underset{\nu^o}{_\alpha\otimes_{\hat{\beta}}}\zeta_2))\]
from which we get the result. \end{proof}

\subsection{Proposition}
\label{prop2Gamma}
{\it Let $W$ be an $(\alpha, \hat{\beta}, \beta)$-pseudo-multiplicative unitary, $\Gamma$ the coproduct constructed in \ref{AW}, $\xi$ in $D(_\alpha\gH, \nu)$, $\eta$ in $D(\gH_{\hat{\beta}}, \nu^o)$. Let $\xi_1$, $\eta_1$ in $D(\gH_\beta, \nu^o)$, $\xi_2$, $\eta_2$ in $D(_\alpha\gH, \nu)$; then, we have :}
\begin{multline*}
(\Gamma((id*\omega_{\xi, \eta})(W))(\xi_1\underset{\nu}{_\beta\otimes_\alpha}\eta_1)|
\xi_2\underset{\nu}{_\beta\otimes_\alpha}\eta_2)=\\
((\omega_{\xi_1, \xi_2}*id)(W)(\omega_{\eta_1, \eta_2}*id)(W)\xi|\eta)
\end{multline*}
\begin{proof}
The scalar product 
\[(\Gamma((id*\omega_{\xi, \eta})(W))(\xi_1\underset{\nu}{_\beta\otimes_\alpha}\eta_1)|
\xi_1\underset{\nu}{_\beta\otimes_\alpha}\eta_2)\]
is equal to, using the definition of $\Gamma$ (\ref{AW}) :
\[((1\underset{\nu^o}{_\alpha\otimes_{\hat{\beta}}}(id*\omega_{\xi, \eta})(W))W(\xi_1\underset{\nu}{_\beta\otimes_\alpha}\eta_1)|W(\xi_2\underset{\nu}{_\beta\otimes_\alpha}\eta_2))\]
or, to :
\[((1\underset{N^o}{_\alpha\otimes_{\hat{\beta}}}W)(W\underset{N}{_\beta\otimes_\alpha}1)(\xi_1\underset{\nu}{_\beta\otimes_\alpha}\eta_1\underset{\nu}{_\beta\otimes_\alpha}\xi)|(W\underset{N^o}{_\alpha\otimes_{\hat{\beta}}}1)((\xi_2\underset{\nu}{_\beta\otimes_\alpha}\eta_2)\underset{\nu^o}{_\alpha\otimes_{\hat{\beta}}}\eta)\]
which, using the pentagonal equation (\ref{defmult}), is equal to the scalar product of the vector :
\[(\sigma_{\nu^o}\underset{N^o}{_\alpha\otimes_{\hat{\beta}}}1_{\gH})
(1_{\gH}\underset{N^o}{_\alpha\otimes_{\hat{\beta}}}W)
\sigma_{2\nu}
(1_{\gH}\underset{N}{_\beta\otimes_\alpha}\sigma_{\nu^o})
(1_{\gH}\underset{N}{_\beta\otimes_\alpha}W)(\xi_1\underset{\nu}{_\beta\otimes_\alpha}\eta_1\underset{\nu}{_\beta\otimes_\alpha}\xi)\]
with the vector $(\xi_2\underset{\nu}{_\beta\otimes_\alpha}\eta_2)\underset{\nu^o}{_\alpha\otimes_{\hat{\beta}}}\eta$, 
which is equal to :
\[(\sigma_{2\nu}
(1_{\gH}\underset{N}{_\beta\otimes_\alpha}\sigma_{\nu^o})
(1_{\gH}\underset{N}{_\beta\otimes_\alpha}W)(\xi_1\underset{\nu}{_\beta\otimes_\alpha}\eta_1\underset{\nu}{_\beta\otimes_\alpha}\xi)|\eta_2\underset{\nu^o}{_\alpha\otimes_{\hat{\beta}}}(W^*(\xi_2\underset{\nu}{_\alpha\otimes_{\hat{\beta}}}\eta)))\]
Defining now $\zeta_i$, $\zeta'_i$ as in \ref{lem1}, we get, using \ref{lem2}, that it is equal to :
\[(W(\eta_1\underset{\nu}{_\beta\otimes_\alpha}\xi)|\eta_2\underset{\nu^o}{_\alpha\otimes_{\hat{\beta}}}\sum_i\alpha(<\zeta_i, \xi_1>_{\beta, \nu^o})\zeta'_i)\]
which, thanks to \ref{lem1}, is equal to :
\[(W(\eta_1\underset{\nu}{_\beta\otimes_\alpha}\xi)|\eta_2\underset{\nu^o}{_\alpha\otimes_{\hat{\beta}}}(\omega_{\xi_1, \xi_2}*id)(W)^*\eta)\]
and, therefore, to 
\[((\omega_{\eta_1, \eta_2}*id)(W)\xi|(\omega_{\xi_1, \xi_2}*id)(W)^*\eta)\]
which finishes the proof.  \end{proof}

\subsection{Manageable pseudo-multiplicative unitaries}
\label{manageable}
Following S.L. Woronowicz ([W3]), F. Lesieur had introduced the notion of a "manageable" pseudo-mutiplicative unitary. 

{\bf Definition} Let $N$ be a von Neumann algebra; let $\gH$ be a Hilbert space on which $N$ has a non degenerate normal representation $\alpha$ and a non-degenerate normal anti-representation $\beta$. Let us suppose that there exists an antilinear bijective unitary $J$ on $\gH$, such that $J^2=1$, and, let us define, for all $n$ in $N$ :
\[\hat{\beta}(n)=J\alpha(n^*)J\]
which gives another non-degenerate normal anti-representation $\hat{\beta}$ of $N$ on $\gH$. We then easily get that $D(\gH_{\hat{\beta}})=JD(_\alpha\gH, \nu)$, and, for all $\xi$, $\eta$ in $D(_\alpha\gH, \nu)$, we have $R^{\hat{\beta}}(J\xi)=JR^\alpha(\xi)J_\nu$, and, therefore :
\[<J\xi, J\eta>_{\hat{\beta}, \nu^o}=<\eta, \xi>_{\alpha, \nu}^o\]
Let now $W$ be a pseudo-multiplicative unitary over the basis $N$, with respect to the representation $\alpha$, and the anti-representations $\beta$ and $\hat{\beta}$, in the sense of  \ref{defmult}; let us suppose that there exists a positive self-adjoint operator $P$ on $\gH$, such that, for all $t$ in $\mathbb{R}$ and $n$ in $N$, we have :
\[P^{it}\alpha(n)P^{-it}=\alpha(\sigma_t^{\nu}(n))\]
\[P^{it}\beta(n)P^{-it}=\beta(\sigma_t^{\nu}(n))\]
\[P^{it}\hat{\beta}(n)P^{-it}=\hat{\beta}(\sigma_t^{\nu}(n))\]
which allows us to define a one-parameter group of unitaries $P^{it}\underset{\nu}{_\beta\otimes_\alpha}P^{it}$ on $\gH\underset{\nu}{_\beta\otimes_\alpha}\gH$, and also $P^{it}\underset{\nu^o}{_\alpha\otimes_{\hat{\beta}}}P^{it}$ on $\gH\underset{\nu^o}{_\alpha\otimes_{\hat{\beta}}}\gH$. Then, $W$ will be said manageable (with managing operator $P$), if we have :
\[W(P^{it}\underset{\nu}{_\beta\otimes_\alpha}P^{it})=(P^{it}\underset{\nu^o}{_\alpha\otimes_{\hat{\beta}}}P^{it})W\]
and if, for all $v$ in $\mathcal D(P^{-1/2})$, $w$ in $\mathcal D(P^{1/2})$, and $p$, $q$ in $D(_\alpha\gH, \nu)\cap D(\gH_{\hat{\beta}}, \nu^o)$, we have :
\[(W^* (v\underset{\nu^o}{_\alpha\otimes_{\hat{\beta}}}q)|w\underset{\nu}{_\beta\otimes_\alpha}p)=
(W(P^{-1/2}v\underset{\nu}{_\beta\otimes_\alpha}Jp)|P^{1/2}w\underset{\nu^o}{_\alpha\otimes_{\hat{\beta}}}Jq)\]

\subsection{Lemma}
\label{lemmanageable}
{\it Let us take the hypothesis of \ref{manageable}. Then, if $p$ belongs to  $D(_\alpha\gH, \nu)\cap D(\gH_{\hat{\beta}}, \nu^o)\cap\mathcal D(P^{1/2})$ such that $P^{1/2}p$ belongs to $D(_\alpha\gH, \nu)$, and $q$ belongs to $D(_\alpha\gH, \nu)\cap D(\gH_{\hat{\beta}}, \nu^o)\cap\mathcal D(P^{-1/2})$ such that $P^{-1/2}q$ belongs to $D(\gH_{\hat{\beta}}, \nu^o)$, then we have :}
\[(id*\omega_{Jp, Jq})(W)^*=(id*\omega_{P^{1/2}p, P^{-1/2}q})(W)\]

\begin{proof}
Let us take $v$ in $\mathcal D(P^{-1/2})$, $w$ in $\mathcal D(P^{1/2})$; then, we have, using \ref{manageable} :
\begin{eqnarray*}
((id*\omega_{Jp, Jq})(W)^*v|w)&=&(v|(id*\omega_{Jp, Jq})(W)w)\\
&=&(v\underset{\nu^o}{_\alpha\otimes_{\hat{\beta}}}Jq|W(w\underset{\nu}{_\beta\otimes_\alpha}Jp))\\
&=&
(W(P^{-1/2}v\underset{\nu}{_\beta\otimes_\alpha}p)|P^{1/2}w\underset{\nu^o}{_\alpha\otimes_{\hat{\beta}}}q)\\
&=&(W(v\underset{\nu}{_\beta\otimes_\alpha}P^{1/2}p)|w\underset{\nu^o}{_\alpha\otimes_{\hat{\beta}}}P^{-1/2}q)\\
&=&((id*\omega_{P^{1/2}p, P^{-1/2}q})(W)v|w)
\end{eqnarray*}
which, by density, gives the result. \end{proof}

\subsection{Lemma}
\label{lemP}
{\it Let us take the hypothesis of \ref{manageable}. For any $p$ in $D(_\alpha\gH, \nu)\cap D(\gH_{\hat{\beta}}, \nu^o)$, there exists a sequence $p_n$ in $D(_\alpha\gH, \nu)\cap D(\gH_{\hat{\beta}}, \nu^o)\cap\mathcal D(P^{1/2})\cap\mathcal D(P^{-1/2})$, such that $P^{1/2}p_n$ belongs to $D(_\alpha\gH, \nu)$, $P^{-1/2}p_n$ belongs to $D(\gH_{\hat{\beta}}, \nu^o)$, and such that $R^{\alpha, \nu}(p_n)$ is weakly converging to $R^{\alpha, \nu}(p)$ and $R^{\hat{\beta}, \nu^o}(p_n)$ is weakly converging to $R^{\hat{\beta}, \nu^o}(p)$. }

\begin{proof}
Let us write :
\[p_n= \frac{\sqrt{n}}{\pi}\int_{-\infty}^{\infty}e^{-nt^2}P^{it}pdt\]
It is a usual calculation to prove that $p_n$ belongs to $\mathcal D(P^{1/2})\cap \mathcal D(P^{-1/2})$; moreover, we get, for any $a$ in $\gN_\nu$ :
\begin{eqnarray*}
\alpha(a)p_n&=&\frac{\sqrt{n}}{\pi}\int_{-\infty}^{\infty}e^{-nt^2}\alpha(a)P^{it}pdt\\
&=&\frac{\sqrt{n}}{\pi}\int_{-\infty}^{\infty}e^{-nt^2}P^{it}\alpha(\sigma^{\nu}_{-t}(a))pdt\\
&=&\frac{\sqrt{n}}{\pi}\int_{-\infty}^{\infty}e^{-nt^2}P^{it}R^{\alpha, \nu}(p)\Delta_{\nu}^{-it}\Lambda_{\nu}(a)dt
\end{eqnarray*}
from which we get that :
\[\|\alpha(a)p_n\|\leq\frac{\sqrt{n}}{\pi}\int_{-\infty}^{\infty}e^{-nt^2}\|R^{\alpha, \nu}(p)\|\|\Lambda_{\nu}(a)\|dt\]
which proves that $p_n$ belongs to $D(_\alpha\gH, \nu)$ and that :
\[\|R^{\alpha, \nu}(p_n)\|\leq\|R^{\alpha, \nu}(p)\|\]
Moreover, we have, going on the same calculation :
\[R^{\alpha, \nu}(p_n)\Lambda_\nu (a)=\frac{1}{\pi}\int_{-\infty}^{\infty}e^{-t^2}P^{\frac{it}{\sqrt{n}}}R^{\alpha, \nu}(p)\Delta_{\nu}^{\frac{-it}{\sqrt{n}}}\Lambda_{\nu}(a)dt\]
which, using Lebesgue's theorem, is converging to $R^{\alpha, \nu}(p)\Lambda_{\nu}(a)$. With the norm majoration, we get this way the weak convergence of $R^{\alpha, \nu}(p_n)$ to $R^{\alpha, \nu}(p)$. The fact that  $R^{\hat{\beta}, \nu^o}(p_n)$ is weakly converging to $R^{\hat{\beta}, \nu^o}(p)$ is obtained the very same way.  \end{proof}

\subsection{Proposition}
\label{propmanag}
{\it Let $N$ be a von Neumann algebra; let $\gH$ be a Hilbert space on which $N$ has a non degenerate normal representation $\alpha$ and a non-degenerate normal anti-representation $\beta$. Let us suppose that there exists an antilinear bijective unitary $J$ on $\gH$, such that $J^2=1$, and, let us define, for all $n$ in $N$ :
\[\hat{\beta}(n)=J\alpha(n^*)J\]
which gives another non-degenerate normal anti-representation $\hat{\beta}$ of $N$ on $\gH$. Let us suppose that, for each $\xi$ in $D(_\alpha\gH, \nu)$, there exists a sequence $\xi_n$ in $D(_\alpha\gH, \nu)\cap D(\gH_{\hat{\beta}}, \nu^o)$ such that $R^{\alpha, \nu}(\xi_n)$ is weakly converging to $R^{\alpha, \nu}(\xi)$. 
\newline
Let now $W$ be a pseudo-multiplicative unitary over the basis $N$, with respect to the representation $\alpha$, and the anti-representations $\beta$ and $\hat{\beta}$, in the sense of  \ref{defmult}; let us suppose that there exists a positive self-adjoint operator $P$ on $\gH$ such that $W$ is manageable, with managing operator $P$, in the sense of \ref{manageable}. Then :
\newline
(i) $W$ is weakly regular; 
\newline
(ii) the weakly closed algebra $A_w(W)$ is closed under the $*$ operation, and is therefore equal to the von Neumann algebra $\mathcal A(W)$;
\newline
(iii) $A_n(W)\cap A_n(W)^*$ is a non degenerate $\mathbb{C}^*$-algebra, which is weakly dense in $\mathcal  A(W)$. Moreover, if $y\in N$ is analytical with respect to $\nu$, then $\alpha(y)$ and $\beta(y)$ belong to the multipliers of this $\mathbb{C}^*$-algebra.
\newline
(iv) if $N$ is abelian, then $A_n(W)$ is a non degenerate $\mathbb{C}^*$-algebra, which is weakly dense in $\mathcal A(W)$; moreover, we have $\alpha(N)\subset M(A_n(W))$ and $\beta(N)\subset M(A_n(W))$. }

\begin{proof}
Let's take $v$ in $D(_\alpha\gH, \nu)\cap\mathcal D(P^{-1/2})$ such that such that $P^{-1/2}v$ belongs to $D(_\alpha\gH, \nu)$, and $w$ in $D(_\alpha\gH, \nu)\cap\mathcal D(P^{1/2})$ such that $P^{1/2}w$ belongs to $D(_\alpha\gH, \nu)$; we have then, for any $p$ in $D(_\alpha\gH, \nu)$ and $q$ in $\gH$ :
\begin{eqnarray*}
(R^{\alpha, \nu}(v)R^{\alpha, \nu}(p)^*q|w)
&=&(R^{\alpha, \nu}(v)J_\nu R^{\hat{\beta}, \nu^o}(Jp)^*q|w)\\
&=&(R^{\alpha, \nu}(v)J_\nu\Lambda_\nu(<Jq, J p>_{\hat{\beta}, \nu^o})|w)\\
&=&(P^{-1/2}R^{\alpha, \nu}(v)J_\nu\Lambda_\nu(<J q, J p>_{\hat{\beta}, \nu^o})|P^{1/2}w)\\
&=&(R^{\alpha, \nu}(P^{-1/2}v)J_\nu\Lambda_\nu(<J q, J p>_{\hat{\beta}, \nu^o})|P^{1/2}w)\\
&=&(\alpha(<J q, J p>_{\hat{\beta}, \nu^o})P^{-1/2}v|P^{1/2}w)\\
&=&(J p\underset{\nu}{_{\hat{\beta}}\otimes_\alpha}P^{-1/2}v|J q\underset{\nu}{_{\hat{\beta}}\otimes_\alpha}P^{1/2}w)
\end{eqnarray*}
There exists $\Xi\in \gH\underset{\nu}{_\alpha\otimes_\beta}\gH$ such that \[\sigma_{\nu^o}W\sigma_\nu\Xi=J p\underset{\nu}{_{\hat{\beta}}\otimes_\alpha}P^{-1/2}v\]
Using \ref{lemP}, there exists elements $p^i_k$ in $D(_\alpha\gH, \nu)\cap D(\gH_{\hat{\beta}}, \nu^o)$, and $v^i_k$ in $D(_\alpha\gH, \nu)\cap\mathcal D(P^{-1/2})$ such that $P^{-1/2}v^i_k$ belongs to $D(\gH_{\alpha}, \nu)$, such that \[\Xi=lim_i\sum_{k=1}^{k=n(i)}Jp^i_k\underset{\nu}{_\alpha\otimes_\beta}P^{-1/2}v^i_k\]
from which we get, using definition \ref{manageable} that $(R^{\alpha, \nu}(v)R^{\alpha, \nu}(p)^*q|w)$ is the limit of :
\begin{eqnarray*}
\sum_k(W(P^{-1/2}v^i_k\underset{\nu^o}{_\beta\otimes_\alpha}J p^i_k)|P^{1/2}w\underset{\nu}{_\alpha\otimes_{\hat{\beta}}}Jq)
&=& (W^*(v^i_k\underset{\nu^o}{_\alpha\otimes_{\hat{\beta}}}q)|w\underset{\nu}{_\beta\otimes_\alpha}p^i_k)\\
&=&((id*\omega_{p^i_k, v^i_k}(\sigma_{\nu^o}W)^*q|w)
\end{eqnarray*}
from which we get that $\theta^{\alpha, \nu}(p,v)=lim_i\sum_k(id*\omega_{p^i_k, v^i_k})(\sigma_{\nu^o}W)$, and that $\theta^{\alpha, \nu}(p,v)$ belongs to $C_w(W)$; using again \ref{lemP}, we get (i). One should note that this proof is mostly borrowed from ([L3], 6.5); the proof that (i) implies (ii) is in ([E2], 3.12(ii)), but we shall prove directly (ii) and (iii) by the following considerations. 
\newline
Thanks to \ref{lemmanageable}, if $p$ belongs to  $D(_\alpha\gH, \nu)\cap D(\gH_{\hat{\beta}}, \nu^o)\cap\mathcal D(P^{1/2})$ such that $P^{1/2}p$ belongs to $D(_\alpha\gH, \nu)$, and if $q$ belongs to $D(_\alpha\gH, \nu)\cap D(\gH_{\hat{\beta}}, \nu^o)\cap\mathcal D(P^{-1/2})$ such that $P^{-1/2}q$ belongs to $D(\gH_{\hat{\beta}}, \nu^o)$, then we have :
\[(id*\omega_{Jp, Jq})(W)^*=(id*\omega_{P^{1/2}p, P^{-1/2}q})(W)\]
and, therefore, $(id*\omega_{Jp, Jq})(W)$ belongs to $A_n(W)\cap A_n(W)^*$. 
\newline
Moreover, if $\xi$ and $\eta$ belong to $D(_\alpha\gH, \nu)\cap D(\gH_{\hat{\beta}}, \nu^o)$, then, using \ref{lemP}, it is possible to construct sequences $p_n$ and $q_n$ verifying such conditions, such that $R^{\hat{\beta}, \nu^o}(p_n)$ is weakly converging to $R^{\hat{\beta}, \nu^o}(\xi)$ (or, equivalently, $R^{\alpha, \nu}(Jp_n)$ is weakly converging to $R^{\alpha, \nu}(J\xi)$) and $R^{\alpha, \nu}(q_n)$ is weakly converging to $R^{\alpha, \nu}(\eta)$ (or, equivalently, $R^{\hat{\beta}, \nu^o}(Jq_n)$ is weakly converging to $R^{\hat{\beta}, \nu^o}(J\eta)$). 
\newline
So, the element $(id*\omega_{J\xi, J\eta})(W)$ is the weak limit of the sequence of operators $(id*\omega_{Jp_n, Jq_n})(W)$, which belong to $A_n(W)\cap A_n(W)^*$. So, the element $(id*\omega_{J\xi, J\eta})(W)$ belongs to $A_w(W)\cap A_w(W)^*$. 
\newline
So, we obtain that, for all $\xi$ in $D(_\alpha\gH, \nu)$, and $\eta$ in $D(\gH_{\hat{\beta}}, \nu^o)$, then the element $(id*\omega_{\xi, \eta})(W)$ belongs also to $A_w(W)\cap A_w(W)^*$. So we obtain that $A_w(W)$ is included into $A_w(W)\cap A_w(W)^*$, which gives (ii). 
\newline
In fact, we have obtained that elements in $A_n(W)\cap A_n(W)^*$ are weakly dense in $\mathcal A(W)$, which, thanks to several remarks made in \ref{AW}, gives (iii).  
\newline
If $N$ is abelian, the weight $\nu$ is a trace, and the managing operator $P$ defined in \ref{manageable} is affiliated to $\alpha (N)'\cap \hat{\beta}(N)'$. Let us write :
\[P=\int_0^{\infty}e_\lambda de_\lambda\]
and let us define $p_n=\int_{1/n}^n de_\lambda$. 
Then $p_n$ is an increasing sequence of projections, weakly converging to $1$, in $\alpha(N)'\cap \hat{\beta}(N)'$; let us take $x$ in $A_{\Phi, T}$ (with the notations of \ref{propbasic}); then the vectors $p_n\Lambda_{\Phi}(x)$ belong to $D(_\alpha H, \nu)\cap D(H_{\hat{\beta}}, \nu^o)\cap \mathcal D(P^{1/2})\cap D(P^{-1/2})$ and both $P^{1/2}p_n\Lambda_\Phi(x)$  and $P^{-1/2}p_n\Lambda_{\Phi}(x)$ belong to $D(_\alpha H, \nu)\cap D(H_{\hat{\beta}}, \nu^o)$. So, using (iii), we get that, for $x$, $y$ in $A_{\Phi, T}$, the operator $(id*\omega_{Jp_n\Lambda_\Phi(x), Jp_n\Lambda_{\Phi}(y)})(W)$ belongs to $A_n(W)\cap A_n(W)^*$. Using ([E2], 10.5), we get, taking the norm limit, that $(id*\omega_{J\Lambda_{\Phi}(x), J\Lambda_{\Phi}(y)})(W)$ belongs to $A_n(W)\cap A_n(W)^*$ for any $x$, $y$ in $A_{\Phi, T}$, which is equivalent, that $(id*\omega_{\Lambda_{\Phi}(x), \Lambda_{\Phi}(y)})(W)$ belongs to $A_n(W)\cap A_n(W)^*$. Using now \ref{basic}, we get that $(id*\omega_{\xi, \eta})(W)$ belongs to $A_n(W)\cap A_n(W)^*$, for any $\xi$ in $D(_\alpha H, \nu)$ and $\eta$ in $D(H_{\hat{\beta}}, \nu^o)$, and, therefore, we get that $A_n(W)\subset A_n(W)\cap A_n(W)^*$, which finishes the proof. \end{proof}

\section{Lesieur's measured quantum groupoids}
\label{quantum}
In this section, we give a r\'esum\'e of Lesieur's construction of "measured quantum groupoids" ([L1], [L2]) (\ref{MQG}) and "generalized measured quantum groupoids"([L3]) (\ref{gMQG}). We then obtain some technical results about the algebra generated by the left leg of the pseudo-multiplicative unitary, especially when the basis is abelian (\ref{corPhi}, \ref{Phi}, \ref{thgamma}).  
\subsection{Definitions ([L1], [L2])}
\label{LW}
Let $(N, M, \alpha, \beta, \Gamma)$ be a Hopf-bimodule, as defined in \ref{Hbimod}; a normal, semi-finite, faithful operator valued weight $T$ from $M$ to $\alpha (N)$ is said to be left-invariant if, for all $x\in \gM_T^+$, we have :
\[(id\underset{N}{_\beta*_\alpha}T)\Gamma (x)=T(x)\underset{N}{_\beta\otimes_\alpha}1\]
A normal, semi-finite, faithful operator-valued weight $T'$ from $M$ to $\beta (N)$ will be said to be right-invariant if it is left-invariant with respect to the symmetrized Hopf-bimodule, i.e., if, for all $x\in\gM_{T'}^+$, we have :
\[(T'\underset{N}{_\beta*_\alpha}id)\Gamma (x)=1\underset{N}{_\beta\otimes_\alpha}T'(x)\]
In the case of a Hopf-bimodule, with a left-invariant normal, semi-finite, faithful operator valued weight $T$ from $M$ to $\alpha (N)$, Lesieur had constructed an isometry $U$ in the following way :
let us choose a normal, semi-finite, faithful weight $\mu$ on $N$, and let us write $\Phi=\nu\circ\alpha^{-1}\circ T$, which is a normal, semi-finite, faithful weight on $M$; let us write $H=H_\Phi$, $J=J_\Phi$, $\Delta=\Delta_\Phi$ for the canonical objects of the Tomita-Takesaki theory associated to the weight $\Phi$, and let us define, for $x$ in $N$, $\hat{\beta}(x)=J\alpha(x^*)J$. Then, there exists an unique isometry $U$ from $H\underset{\mu^o}{_\alpha\otimes_{\hat{\beta}}}H$ to $H\underset{\mu}{_\beta\otimes_\alpha}H$, such that, for any $(\beta, \mu^o)$-orthogonal 
basis $(\xi_i)_{i\in I}$ of  $H_\beta$, for any $a$ in $\gN_T\cap\gN_\Phi$ and for any $v$ in $D(H_\beta, \mu^o)$, we have 
\[U(v\underset{\mu^o}{_\alpha\otimes_{\hat{\beta}}}\Lambda_\Phi (a))=\sum_{i\in I} \xi_i\underset{\mu}{_\beta\otimes_\alpha}\Lambda_{\Phi}((\omega_{v, \xi_i}\underset{\mu}{_\beta*_\alpha}id)(\Gamma(a)))\]
Then, Lesieur proved that, if there exists a right-invariant normal, semi-finite, faithful operator valued weight $T'$ from $M$ to $\beta (N)$, then this isometry is a unitary, and that $W=U^*$ is an $(\alpha, \hat{\beta}, \beta)$-pseudo-multiplicative unitary from $H\underset{\mu}{_\beta\otimes_\alpha}H$ to $H\underset{\mu^o}{_\alpha\otimes_{\hat{\beta}}}H$\vspace{5mm}. \newline
{\bf Proposition}
{\it Let $(N, M, \alpha, \beta, \Gamma)$ be a Hopf-bimodule, as defined in \ref{Hbimod}; let us suppose that there exist  a normal, semi-finite, faithful left-invariant operator valued weight $T$ from $M$ to $\alpha (N)$ and a right-invariant normal, semi-finite, faithful operator valued weight $T'$ from $M$ to $\beta (N)$; let us write $\Phi=\nu\circ\alpha^{-1}\circ T$, and let us define, for $n$ in $N$, $\hat{\beta}(n)=J_\Phi\alpha(n^*)J_\Phi,$; then the $(\alpha, \hat{\beta}, \beta)$-pseudo-multiplicative unitary from $H\underset{\mu}{_\beta\otimes_\alpha}H$ to $H\underset{\mu^o}{_\alpha\otimes_{\hat{\beta}}}H$ verifies, for any $x$, $y_1$, $y_2$ in $\gN_T\cap\gN_\Phi$ :}
\[(i*\omega_{J_\Phi\Lambda_\Phi (y_1^*y_2), \Lambda_\Phi (x)})(W)=
(id\underset{N}{_\beta*_\alpha}\omega_{J_\Phi\Lambda_\Phi(y_2), J_\Phi\Lambda_\Phi(y_1)})\Gamma (x^*)\]
\begin{proof}
This is just ([L2], 3.19). \end{proof}

\subsection{Definition ([L1], [L2])}
\label{MQG}
Lesieur defined a "measured quantum groupoid" as an $8$-uple, $(N, M, \alpha, \beta, \Gamma, T, T', \nu)$, where $(N, M, \alpha, \beta, \Gamma)$ is a Hopf-bimodule, $T$ (resp. $T'$) is a left-invariant (resp. right-invariant) normal, semi-finite, faithful operator-valued weight from $M$ to $\alpha (N)$ (resp. $\beta (N)$), and $\nu$ is a normal semi-finite faithful weight on $N$ such that, for all $t$ in $\mathbb{R}$ and $n$ in $N$ :
\[\alpha(\sigma_t^\nu (n))=\sigma_t^{T'}(\alpha (n))\]
\[\beta (\sigma_{-t}^\nu (n))=\sigma_t^T(\beta (n))\]
and this last axiom allowed him to mimick ([KV]) and to define the analog of an antipode, co-inverse, scaling automorphism group, modulus, scaling operator. More precisely, Lesieur had constructed :
\newline
- a co-inverse $R$, which is a $*$-anti-automorphism of $M$, such that $R^2=id$, $R\circ\alpha=\beta$, satisfying :
\[\Gamma\circ R=\varsigma (R\underset{N}{_\beta*_\alpha}R)\Gamma\]
- a scaling group $\tau_t$, which is a one-parameter group of automorphisms of $M$, such that $\tau_t\circ\alpha=\alpha\circ\sigma_t^{\nu}$, $\tau_t\circ\beta=\beta\circ\sigma_t^{\nu}$, and :
\[\Gamma\circ\tau_t=(\tau_t\underset{N}{_\beta*_\alpha}\tau_t)\Gamma\]
Then $S=R\circ\tau_{i/2}$ is an unbounded antipode, satifying, for all $u$, $v$ in $D(_\alpha H, \nu)\cap D(H_{\hat{\beta}}, \nu^o)$ :
\[S((i*\omega_{v, w})(W))=(i*\omega_{v, w})(W^*)\]
Moreover, we have also :
\[\Gamma\circ\tau_t=(\sigma_t^{\Phi}\underset{N}{_\beta*_\alpha}\sigma_{-t}^{\Phi\circ R})\Gamma\]
- a scaling operator $\lambda$ which is a strictly positive operator affiliated to $Z(M)\cap\alpha(N)\cap\beta(N)$, and a modulus $\delta$, which is a stricltly positive operator affiliated to $M\cap\alpha(N)'\cap\beta(N)'$, such that :
\[(D\Phi\circ R:D\Phi)_t=\lambda^{it^2/2}\delta^{it}\] 
The modulus verifies $\Gamma(\delta)=\delta\underset{\nu}{_\beta\otimes_\alpha}\delta$. \newline
Then, there exists a strictly positive operator $P$ on $H$ defined, for all $x\in \gN_\Phi$ by
\[P^{it}\Lambda_\Phi (x)=\lambda^{t/2}\Lambda_\Phi(\tau_t(x))\]
and the pseudo-multiplicative unitary $W$ is manageable, with managing operator $P$. Moreover, for any $y$ in $M$, we have $\tau_t(y)=P^{it}yP^{-it}$, and $M=A_w(W)=\mathcal A(W)$. 
\newline
Lesieur had also proved that, if $(N, M, \alpha, \beta, \Gamma, T, T', \nu)$ is a quantum groupoid, so is the symmetrized one $(N^o, M, \beta, \alpha, \varsigma_N\circ\Gamma, T', T, \nu^o)$.

\subsection{Fundamental example}
\label{fe}
Let us come back to the fundamental example described in \ref{gd} and use the notations of  \ref{gd}; Lesieur had shown ([L1], [L2]) that the 8-uple :
\[(L^{\infty}(\mathcal G ^{(0)}, \mu), L^{\infty}(\mathcal G , \nu), r_{\mathcal G}, s_{\mathcal G}, \Gamma_{\mathcal G}, T, T', \mu)\]
is a measured quantum groupoid, where $\mu$ is a quasi-invariant measure on $\mathcal G ^{(0)}$, $\nu$ the measure constructed on $\mathcal G $ using $\mu$ and the Haar system $(\lambda^u)_{u\in \mathcal G ^{(0)}}$, $r_{\mathcal G}$, $s_{\mathcal G}$, $\Gamma_{\mathcal G}$ had been defined in \ref{Hbimod}, $T$ and $T'$ are given by, for $f$ positive in $L^{\infty}(\mathcal G , \nu)$, by :
\[T(f)(u)=\int_{\mathcal G ^{(0)}}fd\lambda^u\]
\[T'(f)(u)=\int_{\mathcal G ^{(0)}}fd\lambda_u\]
where $\lambda_u$ is the image of $\lambda^u$ by the application $x\mapsto x^{-1}$.

\subsection{Definition ([L3])}
\label{gMQG}
Unfortunately, this theory, resumed in \ref{MQG}, appeared to be too restrictive : namely, the axiom on the auxilliary weight $\nu$ is too restrictive, and prevents Lesieur from obtaining a nice duality theory, except in the case when $N$ is semi-finite. So Lesieur had to construct a wider class of objects, called "generalized measured quantum groupoids"; these are $10$-uple, $(N, M, \alpha, \beta, \Gamma, T, R, \nu, \tau_t, \gamma_t)$, where $(N, M, \alpha, \beta, \Gamma)$ is a Hopf-bimodule as defined in \ref{Hbimod}, $T$ is a left-invariant normal semi-finite faithful weight, as defined in \ref{LW}, $R$ is a co-inverse as defined in \ref{MQG}, $\tau_t$ is a one parameter automorphism group of $M$, $\gamma_t$ a one parameter automorphism group of $N$, and $\nu$ a normal semi-finite faithful weight on $N$; moreover, these data are supposed to verify (where $\Phi=\nu\circ\alpha^{-1}\circ T$) :
\newline
- $\tau_t\circ\beta=\beta\circ\sigma_t^\nu$
\newline
- $(\tau_t\underset{N}{_\beta*_\alpha}\sigma_t^\Phi)\Gamma=\Gamma\circ\sigma_t^\Phi$
\newline
- $\tau_t\circ R=R\circ\tau_t$
\newline
- for any $a$, $b$ in $\gN_T\cap\gN_\Phi$, we have :
\[R((id\underset{N}{_\beta*_\alpha}\omega_{J_\Phi\Lambda_\Phi(a)})\Gamma(b^*b))=
(id\underset{N}{_\beta*_\alpha}\omega_{J_\Phi\Lambda_\Phi(b)})\Gamma(a^*a)\]
\[\tau_t((id\underset{N}{_\beta*_\alpha}\omega_{J_\Phi\Lambda_\Phi(a)})\Gamma(b^*b))=
((id\underset{N}{_\beta*_\alpha}\omega_{J_\Phi\Lambda_\Phi(\sigma_t^\Phi(a))})\Gamma(\sigma_t^\Phi(b^*b)))\]
\newline
-$\sigma_t^T\circ\beta=\beta\circ\gamma_t$
\newline
-$\nu\circ\gamma_t=\nu$
\newline
These axioms, mostly inspired from [MNW], are clearly two complicated and should be simplified. However, in that case, as $R\circ T\circ R$ is a right-invariant operator valued weight, we can apply \ref{LW} and construct a pseudo-multiplicative unitary $W$; Lesieur, in [L3], proved in that situation a nice duality theory, starting from the dual pseudo-multiplicative unitary  $\widehat{W}=\sigma_{\nu^o}W^*\sigma_{\nu^o}$. He obtained all these data for the dual object, the role of $\beta$ and $\hat{\beta}$ are interchanged, and he obtained also that $\widehat{\gamma_t}=\gamma_{-t}$. 
\newline
Using \ref{MQG}, we see that measured quantum groupoids are generalized measured quantum groupoids; there are exactly those for which $\gamma_t=\sigma_{-t}^\nu$.
\newline
Moreover, the properties listed in \ref{MQG} remains valid in that more general case, except that the modulus $\delta$ no longer commutes with $\alpha$ and $\beta$; more precisely, Lesieur had obtained that, for all $n$ in $N$ and $t$ in $\mathbb{R}$ 
\[\delta^{it}\alpha(n)\delta^{-it}=\alpha(\gamma_t\sigma_t^\nu(n))\]
\[\delta^{it}\beta(n)\delta^{-it}=\beta(\gamma_t\sigma_t^\nu(n))\]
These properties allows us to give a meaning to $\delta\underset{\nu}{_\beta\otimes_\alpha}\delta$, and Lesieur obtained again $\Gamma(\delta)=\delta\underset{\nu}{_\beta\otimes_\alpha}\delta$. 

\subsection{Proposition ([L3], 5.1)}
\label{alphabeta}
{\it Let $(N, M, \alpha, \beta, \Gamma, T, R, \nu, \tau_t, \gamma_t)$ be a generalized measured quantum groupoid in the sense of Lesieur; then :
\newline
(i) the left ideal $\gN_T\cap\gN_{R\circ T\circ R}\cap\gN_\Phi\cap\gN_{\Phi\circ R}$ is dense in $M$; the subspace $\Lambda_\Phi(\gN_T\cap\gN_{R\circ T\circ R}\cap\gN_\Phi\cap\gN_{\Phi\circ R})$ is dense in $H_\Phi$; for any $x$ in $\gN_T\cap\gN_\Phi$, there exist $x_n$ in $\gN_T\cap\gN_{R\circ T\circ R}\cap\gN_\Phi\cap\gN_{\Phi\circ R}$ such that $\Lambda_T(x_n)$ is weakly converging to $\Lambda_T(x)$; if $N$ is abelian, $\Lambda_T(x_n)$ is norm converging to $\Lambda_T(x)$. 
\newline
(ii) for all $\xi$ in $D(_\alpha H, \nu)$, there exists a sequence $\xi_n$ in $D(_\alpha H, \nu)\cap D(H_\beta, \nu^o)$, such that $\xi_n$ is converging to $\xi$, and $R^{\alpha, \nu}(\xi_n)$ is weakly converging to $R^{\alpha, \nu}(\xi)$. If $N$ is abelian, $R^{\alpha, \nu}(\xi_n)$ is norm converging to $R^{\alpha, \nu}(\xi)$.}

\subsection{Theorem}
\label{thL}
{\it Let $(N, M, \alpha, \beta, \Gamma, T, R, \nu, \tau_t, \gamma_t)$ be a generalized measured quantum groupoid in the sense of Lesieur; then, for any $\xi$, $\eta$ in $D(_\alpha H, \nu)$, for all $t$ in $\mathbb{R}$, we have :
\newline
(i) $R((i*\omega_{\xi, J_\Phi\eta})(W))=(i*\omega_{\eta, J_\Phi\xi})(W)$
\newline
(ii) $\tau_t((i*\omega_{\xi, J_\Phi\eta})(W))=(i*\omega_{P^{it}\xi, P^{it}J_\Phi\eta})(W)$
\newline
(iii)   $\sigma_t^\Phi((i*\omega_{\xi, J_\Phi\eta})(W))=
(i*\omega_{\delta^{it}J_\Phi\delta^{-it}J_\Phi\Delta_\Phi^{-it}\xi, P^{it}J_\Phi\eta})(W)$
\newline
$\sigma_t^{\Phi\circ R}((i*\omega_{\xi, J_\Phi\eta})(W))=
(i*\omega_{P^{it}\xi, \delta^{it}J_\Phi\delta^{-it}J_\Phi\Delta_\Phi^{-it}J_\Phi\eta})(W)$}
\begin{proof}
By proposition \ref{LW}, the result (i) a just the application of \ref{propbasic} to the basic prioperty of $R$ given in \ref{gMQG}, and (ii) is given by the properties of the managing operator $P$. 
\newline
Let us take $\xi=J_\Phi\Lambda_\Phi (y_1^*y_2)$, and $\eta=J_\Phi\Lambda_\Phi(x)$, with $x$, $y_1$, $y_2$ in  $\gN_T\cap\gN_\Phi$; then, using \ref{LW}, and \ref{MQG}, we get :
\[\sigma_t^\Phi((i*\omega_{J_\Phi\Lambda_\Phi (y_1^*y_2), \Lambda_\Phi(x)})(W))
=(id\underset{N}{_\beta*_\alpha}\omega_{J_\Phi\Lambda_\Phi(y_2), J_\Phi\Lambda_\Phi(y_1)}\circ\sigma_t^{\Phi\circ R})\Gamma (\tau_t(x^*))\]
which is equal to :
\begin{multline*}
(id\underset{N}{_\beta*_\alpha}\omega_{J_\Phi\Lambda_\Phi(\lambda^{t/2}\sigma_{-t}^{\Phi\circ R}(y_2), J_\Phi\Lambda_\Phi(\lambda^{t/2}\sigma_{-t}^{\Phi\circ R}(y_1)})\Gamma (\tau_t(x^*))\\
=(id\underset{N}{_\beta*_\alpha}\omega_{J_\Phi\Lambda_\Phi(\sigma_{-t}^{\Phi\circ R}(y_2), J_\Phi\Lambda_\Phi(\sigma_{-t}^{\Phi\circ R}(y_1)})\Gamma (\lambda^t\tau_t(x^*))
\end{multline*}
which, using again \ref{LW} and \ref{MQG}, is equal to :
\begin{multline*}
(i*\omega_{J_\Phi\Lambda_\Phi (\sigma_t^{\Phi\circ R}(y_1^*y_2)), \Lambda_\Phi(\lambda^t\tau_t(x))})(W)\\
=(i*\omega_{J_\Phi\delta^{-it}J_\Phi\delta^{it}J_\Phi\Delta_\Phi^{it}\Lambda_\Phi (y_1^*y_2), P^{it}\Lambda_\Phi(x)})(W)
\end{multline*}
which gives the first result of (iii), using \ref{basic}. 
\newline
By similar calculations, we obtain :
\begin{multline*}
\sigma_t^{\Phi\circ R}((i*\omega_{J_\Phi\Lambda_\Phi (y_1^*y_2), \Lambda_\Phi(x)})(W))\\
=(id\underset{N}{_\beta*_\alpha}\omega_{J_\Phi\Lambda_\Phi(y_2), J_\Phi\Lambda_\Phi(y_1)}\circ\tau_t)\Gamma (\sigma_t^{\Phi\circ R}(x^*))
\end{multline*}
which is equal to :
\begin{multline*}
(id\underset{N}{_\beta*_\alpha}\omega_{J_\Phi\Lambda_\Phi(\lambda^{t/2}\tau_t(y_2), J_\Phi\Lambda_\Phi(\lambda^{t/2}\tau_t(y_1)})\Gamma (\sigma_t^{\Phi\circ R}(x^*))\\
=(i*\omega_{J_\Phi\Lambda_\Phi (\lambda^t\tau_t(y_1^*y_2)), \Lambda_\Phi(\sigma^{\Phi\circ R}(x))})(W)
\end{multline*}
from which we obtain the second result of (iii). \end{proof}

\subsection{Corollary}
\label{corPhi}
{\it Let $(N, M, \alpha, \beta, \Gamma, T, R, \nu, \tau_t, \gamma_t)$ be a generalized measured quantum groupoid in the sense of Lesieur;  let $\xi$, $\eta$ in $D(_\alpha H, \nu)$, then :
\newline 
(i) for any $t$ in $\mathbb{R}$, we have $\sigma_t^\Phi(A_n(W))=A_n(W)$.
\newline
(ii) if $<\xi, \xi>_{\alpha, \nu}^o$ belongs to $\mathbb{C}^*(\gamma)$ and $<\eta, \eta>_{\alpha, \nu}^o$ to $\mathbb{C}^*(\sigma^\nu)$, then $(i*\omega_{\xi, J_\Phi\eta})(W)$ belongs to $\mathbb{C}^*(\sigma^\Phi)$; so, if $N$ is abelian, and if $\gamma=id$, then $A_n(W)\subset \mathbb{C}^*(\sigma^\Phi)$. 
\newline
(iii) if $<\xi, \xi>_{\alpha, \nu}^o$ belongs to $\mathbb{C}^*(\sigma^\nu)$ and $<\eta, \eta>_{\alpha, \nu}^o$ to $\mathbb{C}^*(\gamma)$, then $(i*\omega_{\xi, J_\Phi\eta})(W)$ belongs to $\mathbb{C}^*(\sigma^{\Phi\circ R})$; so, if $N$ is abelian, and if $\gamma=id$, then $A_n(W)\subset \mathbb{C}^*(\sigma^{\Phi\circ R})$.
\newline
(iv) if $<\xi, \xi>_{\alpha, \nu}^o$ and $<\eta, \eta>_{\alpha, \nu}^o$ belong to $\mathbb{C}^*(\sigma^\nu)$, then $(i*\omega_{\xi, J_\Phi\eta})(W)$ belongs to $\mathbb{C}^*(\tau)$; so, if $N$ is abelian, then $A_n(W)\subset \mathbb{C}^*(\tau)$.}
\begin{proof}
By the norm continuity of $x\mapsto \sigma_t^\Phi(x)$, the first result is a simple corollary of \ref{thL}(iii). 
Let now be $\nabla$ the self-adjoint positive operator defined on $L^2(N)$ defined, for all $n$ in $\gN_\nu$ and $t$ in $\mathbb{R}$ by :
\[\nabla^{it}\Lambda_\nu(n)=\Lambda_\nu(\gamma_t(n))\]
We have then :
\begin{eqnarray*}
R^{\alpha, \nu}(\delta^{it}J_\Phi\delta^{-it}J_\Phi\Delta_\Phi^{-it}\xi)\Lambda_\nu(n)
&=&\alpha(n)\delta^{it}J_\Phi\delta^{-it}J_\Phi\Delta_\Phi^{-it}\xi\\
&=&\sigma_{-t}^{\Phi\circ R}(\alpha(n))\xi\\
&=&\alpha(\gamma_t(n))\xi\\
&=&R^{\alpha, \nu}(\xi)\nabla^{it}\Lambda_\nu(n)
\end{eqnarray*}
from which we get that $R^{\alpha, \nu}(\delta^{it}J_\Phi\delta^{-it}J_\Phi\Delta_\Phi^{-it}\xi)=R^{\alpha, \nu}(\xi)\nabla^{it}$ and that 
\[<\delta^{it}J_\Phi\delta^{-it}J_\Phi\Delta_\Phi^{-it}\xi, \delta^{it}J_\Phi\delta^{-it}J_\Phi\Delta_\Phi^{-it}\xi>_{\alpha, \nu}^o=\gamma_{-t}(<\xi, \xi>_{\alpha, \nu}^o)\]
Therefore, if the function $t\mapsto\gamma_t(<\xi, \xi>^o)$ is norm continuous, so is the function $t\mapsto\|R^{\alpha}(\delta^{it}J_\Phi\delta^{-it}J_\Phi\Delta_\Phi^{-it}\xi)\|$. 
\newline
On the other hand, we have :
\[R^{\alpha, \nu}(P^{it}\xi)\Lambda_\nu(n)=\alpha(n)P^{it}\xi\\=P^{it}\alpha(\sigma_t^\nu(n))\xi
=P^{it}R^{\alpha, \nu}(\xi)\Delta_\nu^{it}\Lambda_\nu(n)\]
from which we get that $<P^{it}\xi, P^{it}\xi>_{\alpha, \nu}^o=\sigma_{-t}^\nu(<\xi, \xi>_{\alpha, \nu}^o)$; from these results, using \ref{thL}(iii) and \ref{thL}(ii), we get easily (ii), (iii) and (iv).  \end{proof}

\subsection{Proposition}
\label{Phi}
{\it Let $(N, M, \alpha, \beta, \Gamma, T, R, \nu, \tau_t, \gamma_t)$ be a generalized measured quantum groupoid in the sense of Lesieur; then :
\newline
(i) if $x$ is in $\gN_T\cap\gN_\Phi$, and $y$ is in $\gN_T\cap\gN_{R\circ T\circ R}\cap\gN_\Phi$, then the operator $(i*\omega_{J_\Phi\Lambda_\Phi (y^*y), \Lambda_\Phi (x^*x)})(W)$ belongs to $\gM_T^+\cap\gM_\Phi^+$ and we have :
}
\[\Phi( (i*\omega_{J_\Phi\Lambda_\Phi (y^*y), \Lambda_\Phi (x^*x)})(W))=(R\circ T\circ R(y^*y)J_\Phi\Lambda_\Phi(x)|J_\Phi\Lambda_\Phi(x))\]
\begin{multline*}
T( (i*\omega_{J_\Phi\Lambda_\Phi (y^*y), \Lambda_\Phi (x^*x)})(W))\\
=\alpha (<T\circ R(y^*y)J_\Phi\Lambda_\Phi(x),J_\Phi\Lambda_\Phi(x)>_{\alpha, \nu})\\
\leq\|T\circ R(y^*y)\|T(x^*x)
\end{multline*}
{\it (ii) if $x$, $y$ are in $\gN_T\cap\gN_{R\circ T\circ R}\cap\gN_\Phi$, then $(i*\omega_{J_\Phi\Lambda_\Phi (y^*y), \Lambda_\Phi (x^*x)})(W)$ belongs to $\gM_T^+\cap\gM_\Phi^+\cap\gM_{R\circ T\circ T}^+\cap\gM_{\Phi\circ R}^+$. 
\newline
(iii) if $x$, $y$ are in $\gM_T\cap\gM_{R\circ T\circ R}\cap\gM_\Phi\cap\gM_{\Phi\circ R}$, then $(i*\omega_{J_\Phi\Lambda_\Phi (y), \Lambda_\Phi (x)})(W)$ belongs to $\gM_T\cap\gM_\Phi\cap\gM_{R\circ T\circ R}\cap\gM_{\Phi\circ R}$.
\newline
(iv) if $N$ is abelian and $\gamma = id$, the restrctions of $\Phi$ and $\Phi\circ R$ to $A_n(W)$ are faithful lower semi-continuous densely defined KMS weights on $A_n(W)$. }
\begin{proof}
Using \ref{LW}, we obtain that the operator $(i*\omega_{J_\Phi\Lambda_\Phi (y^*y), \Lambda_\Phi (x^*x)})(W)$ is positive, and that :
\[R\circ T( (i*\omega_{J_\Phi\Lambda_\Phi (y^*y), \Lambda_\Phi (x^*x)})(W))
= R\circ T\circ R( (i*\omega_{J_\Phi\Lambda_\Phi (x^*x), \Lambda_\Phi (y^*y)})(W))\]
which, using \ref{thL}(iii) and the right-invariance of $R\circ T\circ R$, is equal to :
\begin{multline*}
R\circ T\circ R((id\underset{N}{_\beta *_\alpha}\omega_{J_\Phi\Lambda_\Phi (x)})\Gamma(y^*y))\\
=(id\underset{N}{_\beta *_\alpha}\omega_{J_\Phi\Lambda_\Phi (x)})(1\underset{N}{_\beta\otimes_\alpha}R\circ T\circ R (y^*y))\\
=\beta(<R\circ T\circ R(y^*y)J_\Phi\Lambda_\Phi(x),J_\Phi\Lambda_\Phi(x)>_{\alpha, \nu})
\end{multline*}
from which we get (i); we then get (ii) by using \ref{thL}(i), and (iii) is just given by linearity. Let us write $\tilde{\Phi}$ for the restriction of $\Phi$ to $A_n(W)$; If $N$ is abelian, using then Dini's theorem and (iii), we get that, if $x$, $y$ are in $\gN_T\cap\gN_\Phi\cap\gN_{\Phi\circ R}\cap\gN_{R\circ T\circ R}$, the operator $(i*\omega_{J_\Phi\Lambda_\Phi (y), \Lambda_\Phi (x)})(W)$ belongs to the norm closure of $\gM_{\tilde{\Phi}}\cap\gM_{\tilde{\Phi}\circ R}$; so, using \ref{alphabeta}(i) and \ref{propbasic}(iii), we get that $\tilde{\Phi}$ is densely defined in $A_n(W)$; moreover, if $\gamma =id$, using \ref{corPhi}(ii), we get that $\tilde{\Phi}$ is KMS, which finishes the proof. \end{proof}

\subsection{Lemma}
\label{lemomega}
{\it Let $\gG=(N, M, \alpha, \beta, \Gamma, T, R, \nu, \tau_t, \gamma_t)$ be a generalized measured quantum groupoid in the sense of Lesieur;  let us define $\Phi=\nu\circ\alpha^{-1}\circ T$; then, we have, for all $x$ in $\gN_\Phi$ :}
\[\omega_{J_\Phi\Lambda_{\Phi}(x)}\circ R=\omega_{J_{\Phi\circ R}\Lambda_{\Phi\circ R}(R(x^*))}\]
\begin{proof}
Let $y$ be analytical with respect to both $\Phi$ and $\Phi\circ R$; we then get that :
\[<\omega_{J_\Phi\Lambda_{\Phi}(x)}, y>=\Phi(\sigma_{i/2}^{\Phi}(y)x^*x)\]
and, therefore :
\begin{eqnarray*}
<\omega_{J_\Phi\Lambda_{\Phi}(x)}\circ R, y>
&=&\Phi(\sigma_{i/2}^{\Phi}(R(y))x^*x)\\
&=&\Phi(R(\sigma_{-i/2}^{\Phi\circ R}(y)x^*x)\\
&=&\Phi\circ R(R(x^*x)\sigma_{-i/2}^{\Phi\circ R}(y))\\
&=&\overline{\Phi\circ R(\sigma_{i/2}^{\Phi\circ R}(y^*)R(x)R(x^*))}
\end{eqnarray*}
which, by a similar calculation, is equal to $\overline{<\omega_{J_{\Phi\circ R}\Lambda_{\Phi\circ R}(R(x^*))}, y^*>}$; which gives the result. \end{proof}

\subsection{Lemma}
\label{lemA}
{\it Let $\gG=(N, M, \alpha, \beta, \Gamma, T, R, \nu, \tau_t, \gamma_t)$ be a generalized measured quantum groupoid in the sense of Lesieur, and let us suppose that the von Neuman algebra $N$ is abelian;  let us use the notations of \ref{propmanag}(iv) and consider the $\mathbb{C}^*$- subalgebra $A_n(W)$ of $M$ (\ref{AW}, \ref{propmanag}(iv)); for any $x$ in $\gN_T\cap\gN_\Phi$, there exists $x_n$ in $A_n(W)\cap\gM_T\cap\gM_\Phi$ such that $\Lambda_T(x_n)$ is norm converging to $\Lambda_T(x)$. }

\begin{proof}
As $\alpha(N)\subset M(A_n(W))$ (\ref{propmanag}(iv)), we get that $T(A_n(W)\cap\gM_T)$ is an ideal of $\alpha (N)$, which, by normality of $T$, is weakly dense in $N$; let $e_n$ be a countable approximate unit of $T(A_n(W)\cap\gM_T)$; we have $e_nT(x^*x)=T(x^*e_nx)$, which is increasing to $T(x^*x)$, and, using Dini's theorem, is therefore norm converging to $T(x^*x)$. Let $f_n$ positive in $A_n(W)$ such that $e_n=T(f_n)$; we have $e_nT(x^*x)=T(T(x^*x)^{1/2}f_nT(x^*x)^{1/2})=T(x_n^*x_n)$, where $x_n=f_n^{1/2}T(x^*x)^{1/2}$ belongs to $A_n(W)\cap\gM_T$. We then get that $\Lambda_T(x_n)$ is norm converging to $\Lambda_T(x)$. \end{proof}

\subsection{Theorem}
\label{thgamma}
{\it Let $\gG=(N, M, \alpha, \beta, \Gamma, T, R, \nu, \tau_t, \gamma_t)$ be a generalized measured quantum groupoid in the sense of Lesieur, and let us suppose that the von Neuman algebra $N$ is abelian;  let us use the notations of \ref{propmanag}(iv) and consider the $\mathbb{C}^*$- subalgebra $A_n(W)$ of $M$ (\ref{AW}, \ref{propmanag}(iv)); then, for all $x_1$, $x_2$ in $A_n(W)\cap\gN_T\cap\gN_\Phi$, $y_1$, $y_2$ in $A_n(W)\cap\gN_{R\circ T\circ R}\cap\gN_{\Phi\circ R}$, we have : 
\newline
(i) we have :
\[(id\underset{N}{_\beta *_\alpha}\omega_{J_\Phi\Lambda_\Phi(x_1), J_\Phi\Lambda_\Phi(x_2)})\Gamma(A_n(W))\subset A_n(W)\]
and the closed linear set generated by all elements of the form 
\[(id\underset{N}{_\beta *_\alpha}\omega_{J_\Phi\Lambda_\Phi(x_1), J_\Phi\Lambda_\Phi(x_2)})\Gamma(x)\]
where $x$ is in $A_n(W)$, $x_1$, $x_2$ in $A_n(W)\cap\gN_T\cap\gN_\Phi$,  is equal to $A_n(W)$. 
\newline
(ii) we have :
\[(\omega_{J_{\Phi\circ R}\Lambda_{\Phi\circ R}(y_1), \Lambda_{\Phi\circ R}(y_2)}\underset{N}{_\beta *_\alpha}id)\Gamma(A_n(W))\subset A_n(W)\]
and the closed linear set generated by all elements of the form 
\[(\omega_{J_{\Phi\circ R}\Lambda_{\Phi\circ R}(y_1), \Lambda_{\Phi\circ R}(y_2)}\underset{N}{_\beta *_\alpha}id)\Gamma(y)\]
 where $y$ is in $A_n(W)$, $y_1$, $y_2$ in $A_n(W)\cap\gN_{R\circ T\circ R}\cap\gN_{\Phi\circ R}$, is equal to $A_n(W)$. }

\begin{proof}
Let us take $x$, $x_1$, $x_2$ in $\gN_T\cap\gN_\Phi$; we have, by prop. \ref{LW} :
\[(id\underset{N}{_\beta *_\alpha}\omega_{J_\Phi\Lambda_\Phi(x_1), J_\Phi\Lambda_\Phi(x_2)})\Gamma(x^*)=(id*\omega_{J_\Phi\Lambda_\Phi(x_1^*x_2), \Lambda_\Phi(x)})(W)\]
If $x$ is in $A_n(W)\cap\gN_T\cap\gN_\Phi$; by the norm density of  $A_n(W)\cap\gN_T\cap\gN_\Phi$ into $A_n(W)$, we get that, for any $y$ in $A_n(W)$, $(id\underset{N}{_\beta *_\alpha}\omega_{J_\Phi\Lambda_\Phi(x_1), J_\Phi\Lambda_\Phi(x_2)})\Gamma(y)$ belongs to $A_n(W)$, from which we get the first result of (i). 
\newline
 Using \ref{LW}, we get that the first closed linear set contains all elements of the form 
$(id*\omega_{J_\Phi\Lambda_\Phi(x_1^*x_2), \Lambda_\Phi(x)})(W)$, where $x$, $x_1$, $x_2$ are in $A_n(W)\cap\gN_T\cap\gN_\Phi$, and, by linearity, all elements of the form $(id*\omega_{J_\Phi\Lambda_\Phi(y), \Lambda_\Phi(x)})(W)$, where $x$ is in $A_n(W)\cap\gN_T\cap\gN_\Phi$ and $y$ is in $A_n(W)\cap\gM_T\cap\gM_\Phi$; using then \ref{lemA}, we get it contains all elements of the form $(id*\omega_{J_\Phi\Lambda_\Phi(y), \Lambda_\Phi(x)})(W)$, where $x$, $y$ belong to $\gN_T\cap\gN_\Phi$; so, by prop \ref{basic}.1(iii), it contains all elements of the form $(i*\omega_{\xi, \eta})(W)$, where $\xi$ is in $D(_\alpha H, \nu)$ and $\eta$ is in $D(H_{\hat{\beta}}, \nu^o)$. So, it contains $A_n(W)$, and, by \ref{thgamma}, it is equal to $A_n(W)$, which is (i). 
\newline
We have now, using \ref{lemomega} :
\begin{multline*}
(\omega_{J_{\Phi\circ R}\Lambda_{\Phi\circ R}(x_1), \Lambda_{\Phi\circ R}(x_2)}\underset{N}{_\beta *_\alpha}id)\Gamma(x)\\
=(\omega_{J_\Phi\Lambda_\Phi(R(x_2^*), J_\Phi\Lambda_\Phi(R(x_1^*))}\circ R\underset{N}{_\beta *_\alpha}id)\Gamma(x)\\
=R((id\underset{N}{_\beta *_\alpha}\omega_{J_\Phi\Lambda_\Phi(R(x_2^*)), J_\Phi\Lambda_\Phi(R(x_1^*))})\Gamma(R(x))
\end{multline*}
which gives (ii). \end{proof}

\section{Measured continuous fields of $\mathbb{C}^*$-algebras }
\label{CXC*}
In that chapter, we prepare the study of the case when the basis is central; for that purpose, we recall the basic definitions and results about $C_0(X)-\mathbb{C}^*$-algebras and continuous fields of $\mathbb{C}^*$-algebras (\ref{defCXC*}). We then give, in that context, the definition of continuous fields  of weights (\ref{deffieldweights}) and a result (\ref{propfieldweights}) which leads to the definition of a "measured continuous field of $\mathbb{C}^*$-algebras"(\ref{defmcf}), which will be basic data for the construction of $\mathbb{C}^*$-quantum groupoids. 
\newline
In the case of two measured continuous fields of $\mathbb{C}^*$-algebras, we show that the min tensor product introduced by Blanchard can be described then with the help of Sauvageot's relative tensor product (\ref{tensor}); moreover, if we take the tensor product of a measurable continuous field of $\mathbb{C}^*$-algebras by itself, the $\mathbb{C}^*$-algebra obtained is again a continuous field of $\mathbb{C}^*$-algebras (\ref{AtensA}).

\subsection{Definitions}
\label{defCXC*}
Let $X$ be a locally compact space; following [Ka], we shall say that a $\mathbb{C}^*$-algebra $A$ is a $C_0(X)-\mathbb{C}^*$-algebra if there exists an injective non degenerate $*$-homomorphism $\alpha$ from $C_0(X)$ into $Z(M(A))$. If $x\in X$, let us write $C_x(X)$ for the ideal of $C_0(X)$ made of all functions in $C_0(X)$ with value $0$ at $x$, and let us consider $\alpha(C_x(X))A$, which is an ideal in $A$; let us write $A^x$ for the quotient $\mathbb{C}^*$-algebra $A/\alpha(C_x(X))A$. For any $a$ in $A$, let us write $a^x$ for its image in $A^x$. Then, we have ([B2], 2.8) :
\[\|a\|=sup_{x\in X}\|a^x\|\]
By definition [D], we shall say that $A$ is a continuous field over $X$ if the function $x\mapsto \| a^x\|$ is continuous. 
\newline
Let $A$ be a $C_0(X)-\mathbb{C}^*$-algebra, and $\mathcal E$ be a $C_0(X)$-Hilbert module; let $\pi$ be a $C_0(X)$-linear morphism $\pi$ from $A$ to $\mathcal L(\mathcal E)$ (which means that the specialization $\pi_x$ is a representation of $A$ on the Hilbert space $\mathcal E_x$ whose kernel contains $\alpha(C_x(X))A$). We say that $\pi$ is a continuous field of faithful representations if $Ker \pi_x=\alpha(C_x(X))A$. We may then consider $\pi_x$ as a faithful representation of the $\mathbb{C}^*$-algebra $A^x$ on $\mathcal E_x$. It is proved in ([B2], 3.3) that, if $A$ is a separable $C_0(X)-\mathbb{C}^*$-algebra, the following are equivalent :
\newline
(i) $A$ is a continuous field over $X$ of $\mathbb{C}^*$-algebras
\newline
(ii) there exists a continuous field of faithful representations of $A$\vspace{5mm}. 
\newline
Given two $C_0(X)$-algebras $A_1$ and $A_2$, Blanchard ([B1], 2.9) had defined the minimal $\mathbb{C}^*$-norm on the involutive algebra $(A_1\odot A_2)/J(A_1, A_2)$, where $A_1\odot A_2$ is the algebraic tensor product of the algebras $A_1$ and $A_2$, and $J(A_1, A_2)$ the involutive ideal in $A_1\odot A_2$ made of finite sums $\sum_{i=1}^n a_i\otimes b_i$, with $a_i\in A_1$, $b_i\in A_2$, such that $\sum_{i=1}^n a_i^x\otimes b_i^x=0$, for all $x\in X$. 
This $\mathbb{C}^*$-norm is given by :
\[\|\sum_{i=1}^n a_i\otimes b_i\|_m=sup_{x\in X}\|\sum_{i=1}^n a_i^x\otimes b_i^x\|_m\]
where, on the right hand on the formula, is taken the minimal tensor product of the $\mathbb{C}^*$-algebras $A_1^x$ and $A_2^x$. The completion with respect to that norm will be called the minimal tensor product of the $C_0(X)$-algebras $A_1$ and $A_2$, and will be denoted $A_1\otimes_{C_0(X)}^m A_2$. 
\newline
In the case of measured fields of $\mathbb{C}^*$-algebras, it is proved in ([B2], 3.21) that this $\mathbb{C}^*$-algebra, equipped with  the  morphism $f\mapsto f\otimes_{C_0(X)}^m 1=1\otimes_{C_0(X)}^m f$ is equal to the $C_0(X)-\mathbb{C}^*$-algebra $A_1\otimes_{C_0(X)}A_2$. (If $\pi_1$ (resp. $\pi_2$) is a faithful non degenerate $C_0(X)$-representation on a $C_0(X)$-Hilbert module $\mathcal E_1$ (resp. $\mathcal E_2$), $A_1\otimes_{C_0(X)}A_2$ is defined as operators on $\mathcal E_1\otimes_{C_0(X)}\mathcal E_2$ ([B2], 3.18) and does not depend on the choice of the $C_0(X)$-representations ([B2], 3.20).

\subsection{Definition}
\label{deffieldweights}
Let $X$ be a locally compact space, and $A$ a $C_0(X)-\mathbb{C}^*$-algebra. A continuous field of weights is a family $\varphi=(\varphi^x)_{x\in X}$ such that :
\newline
- for all $x$ in $X$, $\varphi^x$ is a lower semi-continuous weight on $A$, such that, for all $a$ in $\gM_{\varphi^x}$, $f$ in $C_b(X)^+$, $\alpha(f)a$ belongs to $\gM_{\varphi^x}$ and  :
\[\varphi^x(\alpha(f)a)=f(x)\varphi^x(a)\]
- for all $x\in X$, $a$ in $A^+$, $\varphi^x(a)=0$ if and only if $a\in\alpha(C_x(X))A$. 
\newline
- there exists a norm continuous one parameter group of automorphisms of $A$, $t\mapsto \sigma_t$, such that, for all $x\in X$, the weights $\varphi^x$ are KMS with respect to $\sigma_t$; 
\newline
- for all $x\in X$ and $a\in A^+$, the function $x\mapsto \varphi^x(a)$ from $X$ to $[0, +\infty]$ is Borel, and there exist a norm dense hereditary cone $F\subset\underset{x}{\cap}\gM_{\varphi^x}$ such that, for any $a$ in $F$, the function $x\mapsto \varphi^x(a)$ belongs to $C_b(X)$\vspace{5mm}. 
\newline

We then easily get that $\sigma_t$ is trivial on $\alpha(C_b(X))$, and consider $\sigma_t$ as a one parameter group of automorphisms on $A^x$, and consider then $\varphi^x$ as a faithful lower semi-continuous densely defined weight on $A^x$, which is KMS with respect to $\sigma_t$.

\subsection{Proposition}
\label{propfieldweights}
{\it Let $X$ be a locally compact space, $A$ a $C_0(X)-\mathbb{C}^*$-algebra. and $\varphi$ a continuous field of weights on $A$; then :
\newline
(i) there exists a faithful representation $\pi$ of $A$ such that $\pi\circ\alpha$ is square-integrable with respect to any positive measure $\nu$ on $X$, whose support is $X$; let $\overline{\alpha}$ be the extension of $\pi\circ\alpha$ to $L^\infty(X, \nu)$; 
\newline
(ii) there exists a normal semi-finite faithful operator-valued weight $T$ from $\pi (A)''$ to $\overline{\alpha}(L^\infty(X, \nu))$, such that, for any $a$ in $A^+$, $T(\pi (a))$ is the image by $\alpha$ of the function $x\mapsto \varphi^x(A)$; 
\newline
(iii) for all $t\in\mathbb{R}$, we have : $\pi\circ\sigma_t=\sigma_t^{T}\circ \pi$. }

\begin{proof}
Let $\nu$ be a bounded positive Radon measure on $X$, whose support is $X$; let us define a weight $\Phi$ on $A^+$ by, for all $a\in A^+$ :
\[\Phi (a)=\int_X\varphi^x(a)d\nu '(x)\]
$\Phi$ can be written as an upper limit of lower semi-continuous functions on $A^+$, and is therefore lower semi-continuous; for any $f\in\mathcal K (X)$ and $a\in F$, the function $x\mapsto \varphi^x(\alpha(f)a)$ belongs to $\mathcal K(X)$ and, therefore, to $L^1(\nu)$ and we get that $\Phi(\alpha(f)a)<\infty$; so, $\Phi$ is semi-finite; If $\Phi'(a)=0$, thanks to the support of $\nu $, we get that $\varphi^x(a)=0$ for $\nu$-almost every $x\in X$, which means that $a\in \alpha(C_x(X))A$; therefore, for $\nu$-almost every $x\in X$, we have $a^x=0$, and, by continuity of $x\mapsto \|a^x\|$, we get that $a^x=0$ for all $x\in X$; as $\|a\|=sup_x\|a^x\|$, we finally get that $a=0$, and $\Phi$ is faithful. 
\newline
It is straightforward to get that $\Phi$ is KMS with respect to $\sigma_t$; therefore, there exists an extension $\overline{\Phi}$ on $\pi_{\Phi}(A)''$ which is a normal faithful semi-finite weight, and $\pi_{\Phi}\circ\sigma_t=\sigma_t^{\overline{\Phi}}\circ\pi_{\Phi}$. 
\newline
Let $b\in\gN_{\Phi}$ such that the function $x\mapsto \varphi^x(b^*b)$ is bounded, and $f$ in $C_0(X)$; we have :
\begin{eqnarray*}
\Phi(\alpha (f)b^*b\alpha(\overline{f}))&=&\int_X|f|^2(x)\varphi^x(b^*b)d\nu '(x)\\
&\leq& sup_x\varphi^x(b^*b)\int_X\|f\|^2d\nu' (x)
\end{eqnarray*}
we get that :
\[\|\pi_{\Phi}(\alpha(f))\Lambda_\Phi(b)\|^2\leq sup_x\varphi^x(b^*b)\nu (|f|^2)\]
Therefore, $\Lambda_{\Phi}(b)$ is bounded for the representation $\pi_{\Phi}\circ\alpha$, with respect to the trace $\nu $ on $C_0(X)$; taking now $a\in \gN_{\Phi}$, we can consider the Borel set $B_n=\{x\in X; \varphi^x(a^*a)<n\}$, a compact $K_n\subset B_n$ such that $\nu(B_n-K_n)<1/n$, and a function $f_n$ in $\mathcal K(X)$, $0\leq f_n\leq 1$, with value $1$ on $K_n$; then the functions $x\mapsto \varphi^x(\alpha(f_n)^2a^*a)$ are bounded, so the vectors $\Lambda_\Phi(\alpha(f_n)a)$ are bounded for the representation $\pi_{\Phi}\circ\alpha$, with respect to the trace $\nu $ on $C_0(X)$; as $\alpha(f_n)$ is strongly converging to $1$, we get that these bounded vectors are dense. So, we infer that there exists an extension $\overline{\alpha}$ of $L^\infty (X, \nu)$ on $H_{\Phi}$. 
\newline
If we take now another positive measure $\nu'$ on $X$, with support $X$, by the equivalence of the measures $\nu$ and $\nu '$, we get (i). 
\newline
From the formula $\pi_{\Phi}\circ\sigma_t=\sigma_t^{\overline{\Phi}}\circ\pi_{\Phi}$, we obtain, by continuity, that $\sigma_t^{\overline{\Phi}}(\overline{\alpha}(f))=\overline{\alpha}(f)$ for all $f$ in $L^\infty(X, \nu)$, and, therefore, by Takesaki's theorem ([T], th. 4.2), there exists a faithful semi-finite normal operator-valued weight $T$ from $\pi_{\Phi}(A)''$ to $\overline{\alpha}(L^\infty(X, \nu))$ such that $\overline{\Phi}=\nu\circ\overline{\alpha}^{-1}\circ T$. 
\newline
Another positive measure $\nu'$ on $X$, with support $X$, will be of the form $g\nu $, and the Radon Nykodym derivative $g$ can be considered as a postive non singular element affiliated to $L^\infty(X, \nu)$; taking now the weight $\overline{\Phi'}$ on $\pi_{\Phi}(A)''$ defined by $(D\overline{\Phi'}:D\overline{\Phi})_t=\overline{\alpha}(g)^{it}$, we get, on one hand, 
$\overline{\Phi'}=\nu'\circ\overline{\alpha}^{-1}\circ T$, and, for all $a\in A^+$ :
\[\overline{\Phi'}\circ\pi_{\Phi}(a)=\int_X\varphi^x(a)d\nu'(x)\]
from which we get (ii). 
\newline
As $\overline{\alpha}(L^\infty(X, \nu))$ is central, the modular group $\sigma_t^T$ is equal to $\sigma_t^{\overline{\Phi}}$, and we have already obtained (iii). \end{proof}

\subsection{Definitions}
\label{defmcf}
Let $X$ be a locally compact space, $A$ a $C_0(X)-\mathbb{C}^*$-algebra, and $\varphi$ a continuous field of weights on $A$; let's take the notations of \ref{deffieldweights}; It is then clear that $\gN_\varphi$, considered as a right $C_b(X)$-module, and equipped with the inner product $(a,b)\mapsto T(b^*a)$, is a inner-product $C_b(X)$-module in the sense of [La]. (Sorry for being old-fashioned in writing inner products left linear). 
\newline
Using \ref{propfieldweights}, we can see its completion $\mathcal E_\varphi$ as the norm closure of the set $\{\Lambda_T(X), X\in \gN_\varphi\}$; the representation $\pi$ of $A$ defined in \ref{propfieldweights} is a $C_b(X)$-linear morphism and can be considered as a $C_b(X)$-linear morphism from $A$ into $\mathcal L(\mathcal E_\varphi)$, and we shall denote it $\pi_\varphi$. Taking the specialization at the point $x\in X$, we obtain an Hilbert space $\mathcal (E_\varphi)_x$, which is the completion of the inner product in $\gN_\varphi$ given by $(a, b)\mapsto \varphi^x(b^*a)$; from which we get that $\mathcal (E_\varphi)_x=H_{\varphi^x}$, and that the representation $(\pi_\varphi)_x$ obtained by the specialization of $\pi_\varphi$ is equal to $\pi_{\varphi^x}$. As $\varphi^x$, by definition, is faithful on $A^x$, so is $\pi_{\varphi^x}$, and, therefore, $\pi_\varphi$ is a field of faithful representations of $A$ ([B2], 2.11). 
\newline
Using now ([B2], th 3.3), we get that $A$ is a continuous field of $\mathbb{C}^*$-algebras. 
With all these notations, we shall call $(X, \alpha, A, \varphi)$ a measured continuous field over $X$ of $\mathbb{C}^*$-algebras; for all positive measure $\nu$ on $X$, the weight $\Phi$ on $A$ defined in \ref{propfieldweights} for all $a$ in $A^+$ by :
\[\Phi(a)=\int_X\varphi^x(a)d\nu(x)\]\
will be denoted $\Phi=\int_X\varphi^xd\nu(x)$\vspace{5mm}. 
\newline
Extending slightly the definition of a continuous field of states given in ([B1] 1.2), we shall say that a continuous field of positive linear forms on $A$ is a positive $C_0(X)$-linear application $\omega$ from $A$ into $M(\alpha(C_0(X)))$, such that, for any $x\in X$, the application $X\mapsto\omega(X)(x)$ factors throught a positive linear form on the $\mathbb{C}^*$-algebra $A^x$. In the situation of a measured continuous field $(X, \alpha, A, \varphi)$ over $X$ of $\mathbb{C}^*$-algebras, if $a$ belongs to $\gN_\varphi$, then the applications $\omega_a(b)=T(a^*ba)$ are continuous fields of positive linear forms on $A$; $\omega_a(b)$ is the image by $\alpha$ of the function $x\mapsto \varphi^x(a^*ba)$.

\subsection{Proposition}
\label{tensor}
{\it Let $(X, \alpha_1, A_1, \varphi_1)$ and $(X, \alpha_2, A_2, \varphi_2)$ be two measured continuous fields of $\mathbb{C}^*$-algebras. Let $\nu$ be a positive radon measure on $X$, with support $X$. Let $\Phi_1=\int_X\varphi^x_1d\nu(x)$ (resp. $\Phi_2=\int_X\varphi^x_2d\nu(x)$) be the faithful, densely defined weight on $A$ constructed in \ref{propfieldweights}. Then :
\newline
(i) the representation $\pi_{\Phi_1}\circ\alpha_1$ (resp. $\pi_{\Phi_2}\circ\alpha_2$) is a  representation of $C_0(X)$ on $H_{\Phi_1}$ (resp. $H_{\Phi_2}$), which is square-integrable with respect to $\nu$. Let us write $\overline{\alpha_1}$ (resp. $\overline{\alpha_2}$) the representation of $L^\infty(X, \nu)$ which extends $\alpha_1$ (resp. $\alpha_2$), and let us write again $\nu$ the normal semi-finite faithful trace on $L^\infty(X, \nu)$ given by the positive Radon measure $\nu$. 
\newline
(ii) the $\mathbb{C}^*$-algebra $A_1\otimes_{C_0(X)}^m A_2$ has a faithful representation $\varpi$ on the Hilbert space $H_{\Phi_1}\underset{\nu}{_{\overline{\alpha_1}}\otimes_{\overline{\alpha_2}}}H_{\Phi_2}$ such that, for all $a_1\in A_1$ and $a_2\in A_2$, we have :
\[\varpi(a_1\otimes a_2)=\pi_{\Phi_1}(a_1)\underset{N}{_{\overline{\alpha_1}}\otimes_{\overline{\alpha_2}}}\pi_{\Phi_2}(a_2)\]
(iii) for any finite sum with $a_i\in A_1$ and $b_i\in A_2$, we have :
\[\|\sum_{i=1}^n a_i\otimes b_i\|_m=\|\sum_{i=1}^n\pi_{\Phi_1}(a_i)\underset{N}{_{\overline{\alpha_1}}\otimes_{\overline{\alpha_2}}}\pi_{\Phi_2}(b_i)\|\]
We shall therefore identify the $\mathbb{C}^*$-algebra $A_1\otimes_{C_0(X)}^m A_2$ with the $\mathbb{C}^*$ generated on $H_{\Phi_1}\underset{\nu}{_{\overline{\alpha_1}}\otimes_{\overline{\alpha_2}}}H_{\Phi_2}$ by all the operators $\pi_{\Phi_1}(a_1)\underset{N}{_{\overline{\alpha_1}}\otimes_{\overline{\alpha_2}}}\pi_{\Phi_2}(a_2)$, or all $a_1\in A_1$ and $a_2\in A_2$. 
\newline
(iv) the ideal $J(A_1, A_2)$ is generated by the elements $a_1\alpha_1(f)\otimes a_2-a_1\otimes\alpha_2(f)a_2$, where $a_1$ belongs to $A_1$, $a_2$ to $A_2$ and $f$ to $C_0(X)$. }

\begin{proof}
Result (i) is clear by \ref{propfieldweights}(i). 
\newline
By ([B1], 4.1), there exists a faithful $C_0(X)$-linear representation of $A_1\otimes_{C_0(X)}^m A_2$ on the Hilbert $A_2$-module $\mathcal E_{\varphi_1}\otimes_{C_0(X)}A_2$, which sends the finite sum $\sum_{i=1}^n a_i\otimes b_i$ on the operator $\sum_{i=1}^n\pi_{\varphi_1}(a_i)\otimes_{C_0(X)}b_i$ on $\mathcal E_{\varphi_1}\otimes_{C_0(X)}A_2$. Let have a closer look at this last operator, and let's take finite families $x_j\in \gN_{\varphi_1}$, $c_j\in \gN_{\Phi_2}$ ($j=1, ..., m$). With a repeated use of Cauchy-Schwartz inequality, and with the same arguments as in ([L], 1.2), one gets that the weight $\Phi_2$ applied to :
\[<\sum_{i=1, j=1}^{i=n, j=m}\pi_{\varphi_1}(a_i)\Lambda_{T_1}(x_j)\otimes_{C_0(X)}b_ic_j,  
\sum_{i=1, j=1}^{i=n, j=m}\pi_{\varphi_1}(a_i)\Lambda_{T_1}(x_j)\otimes_{C_0(X)}b_ic_j>\]
is less than :
\[\|\sum_{i=1}^n\pi_{\varphi_1}(a_i)\otimes_{C_0(X)}b_i\|\Phi_2(<\sum_{j=1}^m\Lambda_{T_1}(x_j)\otimes_{C_0(X)}c_j, \sum_{j=1}^m\Lambda_{T_1}(x_j)\otimes_{C_0(X)}c_j>)\]
But, we easily get that :
\begin{multline*}
\Phi_2(<\sum_{j=1}^m\Lambda_{T_1}(x_j)\otimes_{C_0(X)}c_j, \sum_{j=1}^m\Lambda_{T_1}(x_j)\otimes_{C_0(X)}c_j>)=\\
\Phi_2(\sum_j\alpha_2\circ\alpha_1^{-1}T_1(x_j^*x_j)c_j^*c_j)=\\
\|\sum_{j=1}^m \Lambda_{\Phi_1}(x_j)\underset{\nu}{_{\overline{\alpha_1}}\otimes_{\overline{\alpha_2}}}\Lambda_{\Phi_2}(c_j)\|^2
\end{multline*}
and, therefore, we get that :
\[\|\sum_{i=1, j=1}^{i=n, j=m}\pi_{\Phi_1}(a_i)\Lambda_{\Phi_1}(x_j)
\underset{\nu}{_{\overline{\alpha_1}}\otimes_{\overline{\alpha_2}}}
\pi_{\Phi_2}(b_i)\Lambda_{\Phi_2}(c_j)\|^2\]
is less than :
\[\|\sum_{i=1}^n\pi_{\varphi_1}(a_i)\otimes_{C_0(X)}b_i\|^2
\|\sum_{j=1}^m \Lambda_{\Phi_1}(x_j)\underset{\nu}{_{\overline{\alpha_1}}\otimes_{\overline{\alpha_2}}}\Lambda_{\Phi_2}(c_j)\|^2\]
which, using ([B1], 4.1), is less than :
\[\|\sum_{i=1}^n a_i\otimes b_i\|_m^2
\|\Lambda_{\Phi_1}(x_j)\underset{\nu}{_{\overline{\alpha_1}}\otimes_{\overline{\alpha_2}}}\Lambda_{\Phi_2}(c_j)\|^2\]
From which we induce that :
\[\|\sum_{i=1}^n\pi_{\Phi_1}(a_i)\underset{N}{_{\overline{\alpha_1}}\otimes_{\overline{\alpha_2}}}\pi_{\Phi_2}(b_i)\|\leq\|\sum_{i=1}^n a_i\otimes b_i\|_m\]
which gives (ii). 
\newline
Let us suppose now that $\sum_{i=1}^n\pi_{\Phi_1}(a_i)\underset{N}{_{\overline{\alpha_1}}\otimes_{\overline{\alpha_2}}}\pi_{\Phi_2}(b_i)=0$; with the same calculation as above, using the faithfulness of $\Phi_2$, we get that, for any finite families $(x_j)_{j=1,..m}$ and $(c_j)_{j=1, ..m}$, we have :
\[\sum_{i=1, j=1}^{i=n, j=m}\pi_{\varphi_1}(a_i)\Lambda_{T_1}(x_j)\otimes_{C_0(X)}b_ic_j=0\]
which gives that the operator $\sum_{i=1}^n\pi_{\Phi_1}(a_i)\otimes_{C_0(X)}b_i$ on $\mathcal E_{\varphi_1}\otimes_{C_0(X)}A_2$ is equal to $0$. By the faithfulness of the representation constructed in ([B1] 4.1), we get that $\|\sum_{i=1}^n a_i\otimes b_i\|_m=0$, and, therefore, that $\sum_{i=1}^n a_i\otimes b_i$ belongs to $ J(A_1, A_2)$. But, as the semi-norm $\sum_{i=1}^n a_i\otimes b_i\mapsto \|\sum_{i=1}^n a_i\otimes b_i\|_m$ is the minimal semi-norm on $(A_1\odot A_2)/J(A_1, A_2)$ ([B1], 2.9), we get (iii). Now (iv) is given by ([B1], 3.1). \end{proof}

This last proposition extends ([GM], 3.4).

\subsection{Example}
\label{gd2}
Let $\mathcal G$ be a locally compact groupoid, as described in \ref{gd}. Then, the transpose of applications $r$ and $s$ give us mappings $\hat{r}$ and $\hat{s}$ from $C_0(\mathcal G^{(0)})$ into $C_b(\mathcal G)=M(C_0(\mathcal G))$, which gives to $C_0(\mathcal G)$ a double (via $r$, and via $s$) structure of $C_0(\mathcal G^{(0)})-\mathbb{C}^*$-algebra. 
\newline
As $r$ and $s$ are open ([R1], 2.4), $C_0(\mathcal G)$ is (in two different ways, via $r$ and via $s$) a continuous field of $\mathbb{C}^*$-algebras ([B2], prop. 3.14). Moreover, if $f\in\mathcal K(\mathcal G)$ the function $u\mapsto \lambda^u(f)$ belongs to $\mathcal K(\mathcal G^{(0)})$ ([R1] 2.3). So, it is clear that $\lambda^u$ (resp. $\lambda_u$) is then a continuous field of weights on $C_0(\mathcal G)$, via $r$ (resp. via $s$), in the sense of \ref{deffieldweights}. 
\newline
By applying \ref{tensor}(iv), we get that $C_0(\mathcal G)\underset{C_0(\mathcal G^{(0)})}{_{\hat{s}}\otimes^m_{\hat{r}}}C_0(\mathcal G)$ is the quotient of $C_0(\mathcal G)\otimes C_0(\mathcal G)$ (which can be identified with $C_0(\mathcal G^2)$) by the ideal generated by the functions $(x_1, x_2)\mapsto f(s(x_1))g(x_1, x_2)-f(r(x_2))g(x_1, x_2)$, where $x_1$, $x_2$ are in $\mathcal G$, $f$ in $C_0(\mathcal G^{(0)})$, $g$ in $C_0(\mathcal G^2)$; a non zero character on this algebra is a couple $(x_1, x_2)$ such that $s(x_1)=r(x_2)$, i.e. an element of $\mathcal G^{(2)}$; so, we can identify $C_0(\mathcal G)\underset{C_0(\mathcal G^{(0)})}{_{\hat{s}}\otimes^m_{\hat{r}}}C_0(\mathcal G)$ with $C_0(\mathcal G^{(2)})$.

\subsection{Proposition}
\label{AtensA}
{\it Let $(X, \alpha, A, \varphi)$ be a measurable continuous field of $\mathbb{C}^*$-algebras; let $\pi_\varphi$ be the representation $C_0(X)-\mathbb{C}^*$-algebra $A$ on the $\mathbb{C}^*$-module $\mathcal E_\varphi$; let us consider the $C_0(X)-\mathbb{C}^*$-algebra $A\underset{C_0(X)}{\otimes^m}A$; then, for any $x\in X$ :
\newline
(i)  for any finite family $a_i$, $b_i$ of elements of $A$, we have :
\[\|(\pi_{\varphi^x}\otimes_m\pi_{\varphi^x})(\sum_i a_i\otimes b_i)\|\leq\|\sum_i(a_i\underset{C_0(X)}{\otimes^m}b_i)\|\]
and, therefore, we can define a representation $\pi_{\varphi^x}\underset{C_0(X)}{\otimes^m}\pi_{\varphi^x}$ of $A\underset{C_0(X)}{\otimes^m}A$, which is equal to $(\pi_\varphi\otimes_{C_0(X)}\pi_\varphi)_x$; 
\newline
(ii) we have :
\[Ker(\pi_\varphi\otimes_{C_0(X)}\pi_\varphi)_x=\alpha(C_x(X))(A\underset{C_0(X)}{\otimes^m}A)\]
(iii) the $C_0(X)-\mathbb{C}^*$-algebra $A\underset{C_0(X)}{\otimes^m}A$ is a continuous field of $\mathbb{C}^*$-algebras and $(A\underset{C_0(X)}{\otimes^m}A)^x=A^x\otimes_mA^x$}

\begin{proof}
The tensor product $\mathcal E_\varphi\otimes_{C_0(X)}\mathcal E_\varphi$ is a $C_0(X)$-module defined, for $y_1$, $y_2$, $z_1$, $z_2$ in $\gN_\varphi$, by :
\begin{multline*}
<\Lambda_T(y_1)\otimes_{C_0(X)}\Lambda_T(z_1), \Lambda_T(y_2)\otimes_{C_0(X)}\Lambda_T(z_2)>\\
=<T(y_2^*y_1)\Lambda_T(z_1), \Lambda_T(z_2)>
\end{multline*}
which is equal to :
\[T(z_2^*T(y_2^*y_1)z_1)=T(T(y_2^*y_1)z_2^*z_1)=T(y_2^*y_1)T(z_2^*z_1)\]
and the $C_0(X)$-representation $\pi_\varphi\otimes_{C_0(X)}\pi_\varphi$ is therefore defined by :
\begin{multline*}
<(\pi_\varphi\otimes_{C_0(X)}\pi_\varphi)(\sum_ia_i\otimes_{C_0(X)}b_i)(\Lambda_T(y_1)\otimes_{C_0(X)}\Lambda_T(z_1)), \Lambda_T(y_2)\otimes_{C_0(X)}\Lambda_T(z_2)>\\
=\sum_iT(y_2^*a_iy_1)T(z_2^*b_iz_1)
\end{multline*}
Specializing, we obtain, because $(\mathcal E_\varphi)_x=H_{\varphi^x}$ and $(\pi_\varphi)_x=\pi_{\varphi^x}$, that  :
\[((\pi_\varphi\otimes_{C_0(X)}\pi_\varphi)_x(\sum_ia_i\otimes_{C_0(X)}b_i)(\Lambda_{\varphi^x}(y_1)\otimes\Lambda_{\varphi^x}(z_1))| \Lambda_{\varphi^x}(y_2)\otimes\Lambda_{\varphi^x}(z_2))\]
is equal to :
\[\sum_i\varphi^x(y_2^*a_iy_1)\varphi^x(z_2^*b_iz_1)=\sum_i \varphi^x(y_2^*a_iy_1) \varphi^x(z_2^*b_iz_1)\]
which is :
\[((\pi_{\varphi^x}\otimes\pi_ {\varphi^x})(\sum_ia_i\otimes b_i)(\Lambda_{\varphi^x}(y_1)\otimes\Lambda_{\varphi^x}(z_1))| \Lambda_{\varphi^x}(y_2)\otimes\Lambda_{\varphi^x}(z_2))\]
from which we get (i). 
\newline
Let now $Y$ in $Ker (\pi_\varphi\otimes_{C_0(X)}\pi_\varphi)^x$; let us write $Y$ as an infinite sum
$\sum_ia_i\underset{C_0(X)}{\otimes^m}b_i$, norm converging in $A\underset{C_0(X)}{\otimes^m}A$; by the same calculation as above, we obtain, for all $y$ and $z$ in $\gN_\varphi$, that :
\[\sum_i\varphi^x(y^*a_i^*a_iy)\varphi^x(z^*b_i^*b_iz)=0\]
which, implies, that, for all $i$, either $\pi_{\varphi^x}(a_i)$, or $\pi_{\varphi^x}(b_i)$ is equal to $0$; which, using \ref{propfieldweights}, means that, for all $i$, either $a_i$, or $b_i$ belongs to $\alpha(C_x(X))A$; and, therefore, for all $i$, the element $a_i\underset{C_0(X)}{\otimes^m}b_i$ belongs to $(\alpha(C_x(X))\underset{C_0(X)}{\otimes^m}1)(A\underset{C_0(X)}{\otimes^m}A)$, which implies (ii). 
\newline
We get (iii) by ([B2], 3.3 and 3.22). \end{proof}

\section{Measured quantum groupoids with a central basis}
\label{chapcentral}
We deal now with a generalized measured quantum groupoid $\gG=(N, M, \alpha, \beta, \Gamma, T, R, \nu, \tau_t, \gamma_t)$, such that the von Neuman algebra $\alpha(N)$ is included into the center $Z(M)$. Then, we get that it is a measured quantum groupoid (\ref{propcentral}(i)), and that the $\mathbb{C}^*$-algebra $A_n(W)$ is a $\mathbb{C}_0(X)-\mathbb{C}^*$-algebra (\ref{propcentral}(ii)), which is a measured continuous field of $\mathbb{C}^*$-algebras, in two different ways (via $\alpha$, or via $\beta$). We obtain a Plancherel-like formula for the coproduct $\Gamma$ (\ref{thcentral1}), which gives that the coproduct sends the $\mathbb{C}^*$-algebra $A_n(W)$ in the min tensor product of the $\mathbb{C}_0(X)-\mathbb{C}^*$-algebra $A_n(W)$ by itself (\ref{thcentral2}). 
We finish by considering the case of an "abelian" measured quantum groupoid (i.e. the case when the underlying von Neuman algebra itself is abelian; then we prove that it is possible to put on the spectrum of the $\mathbb{C}^*$-algebra $A_n(W)$ a structure of a locally compact groupoid, whose basis is the spectrum of $\mathbb{C}^*(\nu)$ (\ref{thgroupoid}). Starting from a measured groupoid equipped with a left-invariant Haar system, we recover Ramsay's theorem which says that this groupoid is measure-equivalent to a locally compact one (\ref{ramsay}). 

\subsection{Lemma}
\label{propcentral}

{\it Let $\gG=(N, M, \alpha, \beta, \Gamma, T, R, \nu, \tau_t, \gamma_t)$ be a generalized measured quantum groupoid, such that the von Neuman algebra $\alpha(N)$ is included into the center $Z(M)$; then :
\newline
(i) the generalized measured quantum groupoid is a measured quantum groupoid (therefore $\gamma_t=id$); moreover, 
the von Neuman algebra $\beta(N)$ is included into the center $Z(M)$ and the representation $\hat{\beta}$ is equal to $\alpha$.
\newline
(ii) Let $X$ be the spectrum of $\mathbb{C}^*(\nu)$; then the $\mathbb{C}^*$-algebra $A_n(W)$ is, in two different ways (via $\alpha$ and via $\beta$) a $\mathbb{C}_0(X)-\mathbb{C}^*$-algebra, in the sense of Kasparov-Blanchard ([Ka], [B1]).}
\begin{proof}
Let $\delta$ be the modulus of $\gG$; as $\delta^{it}$ commutes with $\alpha (n)$, for all $x$ in $N$ and $t$ in $\mathbb{R}$, we get, using \ref{gMQG} that $\gamma_t=id$, which proves that $\gG$ is a measured quantum groupoid. As $\beta=R\circ\alpha$, we get that $\beta (N)\subset Z(M)$. Moreover, the definition of $\hat{\beta}$ given in \ref{LW} finishes the proof of (i). 
\newline
By \ref{propmanag}(iv), we get that $\alpha(\mathbb{C}^*(\nu))\subset M(A_n(W))$, and with the hypothesis, we get that $\alpha (\mathbb{C}^*(\nu))\subset Z(M(A_n(W))$, which gives the result; the same holds if we take $\beta$ instead of $\alpha$. 
\end{proof}

\subsection{Theorem}
\label{continuousfield}
{\it Let $\gG=(N, M, \alpha, \beta, \Gamma, T, T', \nu)$ be a measured quantum groupoid, such that the von Neuman algebra $\alpha(N)$ is included into the center $Z(M)$, and let $X$ be the spectrum of $\mathbb{C}^*(\nu)$. As $N=M(\mathbb{C}^*(\nu))$, we shall identify $N$ with $C_b(X)$, and, for any $a$ in $\gM_T^+$ (resp. $\gM_{T'}^+$), we define $T^x(a)$ (resp. $T'^x(a)$) as $T(a)(x)$ (resp. $T'(a)(x)$); we obtain then lower semi-continuous weights $T^x$ on $\gM_T^+$ (resp. $T'^x$ on $\gM_{T'}^+$), and we shall denote $\tilde{T}^x$ (resp. $\tilde{T'}^x$) their extension to $A^+$ (\ref{C*}, [B]) :
\newline
(i) Let $\tilde{T}=(\tilde{T}^x)_{x\in X}$; it is a measured continuous field of $\mathbb{C}^*$-algebras over $X$, in the sense of \ref{defmcf}.
\newline
(ii) Let $\tilde{T'}=(\tilde{T'}^x)_{x\in X}$; it is a measured continuous field of $\mathbb{C}^*$-algebras over $X$, in the sense of \ref{defmcf}.}
\begin{proof}
For any $a\in\gM_T^+$, the function $x\mapsto \tilde{T}^x(a)$ is equal to $T(a)$, and is therefore in $N=C_b(X)$; then we get that it is a continuous field of weights, in the sense of \ref{deffieldweights}, which gives (i).  Property (ii) is proved the same way. \end{proof}

\subsection{Corollary}
\label{cortensor}
{\it Let $\gG=(N, M, \alpha, \beta, \Gamma, T, T', \nu)$ be a measured quantum groupoid, such that the von Neuman algebra $\alpha(N)$ is included into the center $Z(M)$. If $x$, $y$ are in $A_n(W)$, the application $x\otimes y\mapsto \|x\underset{N}{_\beta\otimes_\alpha}y\|$ extends to a $\mathbb{C}^*$-semi-norm on the algebraic tensor product $A_n(W)\otimes A_n(W)$ and to a $\mathbb{C}^*$-norm on the quotient of this algebraic tensor product by the ideal generated by the operators of the form \[\{x\beta(f)\otimes y-x\otimes\alpha(f)y, x,y\in A_n(W), f\in N\}\]
Therefore, the $\mathbb{C}^*$-algebra generated on $H\underset{\nu}{_\beta\otimes_\alpha}H$ by all operators $x\underset{N}{_\beta\otimes_\alpha}y$, for $x$, $y$ in $A_n(W)$ can be considered as the min tensor product of the $C_0(X)-\mathbb{C}^*$-algebra $A_n(W)$ (via $\beta$) by the $C_0(X)-\mathbb{C}^*$-algebra $A_n(W)$ (via $\alpha$). }
\begin{proof}
This is just an application of \ref{tensor}. \end{proof}

\subsection{Lemma}
\label{lemcentral}
{\it In the situation of \ref{propcentral}, 
let $(e_i)_{i\in I}$ be an $(\alpha, \nu)$-orthogonal basis of $H$; then, we have : 
\newline
(i) for all $\xi$, $\eta_2$ in $D(_\alpha H, \nu)$, and $\eta_1$ in $D(_\alpha H, \nu)\cap D(H_\beta, \nu)$ :
\[(\omega_{\eta_1, \eta_2}*id)(W)\xi=\sum_i\alpha(<id*\omega_{\xi, e_i})(W)\eta_1, \eta_2>_{\alpha, \nu})e_i\]
(ii) for all $\xi_1$ in $D(H_\beta, \nu)$, $\xi_2$ in $D(_\alpha H, \nu)\cap D(H_\beta, \nu)$ and $\eta$ in $D(_\alpha H, \nu)$ 
\[(\omega_{\xi_1, \xi_2}*id)(W)^*\eta=\sum_i\alpha(<\xi_2, (id*\omega_{e_i, \eta})(W)\xi_1>_{\beta, \nu})e_i\]}
\begin{proof}
We have :
\[W(\eta_1\underset{\nu}{_\beta\otimes_\alpha}\xi)=\sum_i(id*\omega_{\xi, e_i})(W)\eta_1\underset{\nu}{_\alpha\otimes_\alpha}e_i\]
Let now $\zeta$ be in $H$; we have then :
\begin{eqnarray*}
((\omega_{\eta_1, \eta_2}*id)(W)\xi |\zeta)&=&
\sum_i((id*\omega_{\xi, e_i})(W)\eta_1\underset{\nu}{_\alpha\otimes_\alpha}e_i|\eta_2\underset{\nu}{_\alpha\otimes_\alpha}\zeta)\\
&=&(\sum_i\alpha(<id*\omega_{\xi, e_i})(W)\eta_1, \eta_2>_{\alpha, \nu})e_i|\zeta)
\end{eqnarray*}
from which we get (i).
\newline
We have :
\[W^*(\xi_2\underset{\nu}{_\alpha\otimes_\alpha}\eta)=\sum_i(id*\omega_{e_i, \eta})(W)^*\xi_2\underset{\nu}{_\beta\otimes_\alpha}e_i\]
Let now $\zeta$ be in $H$; we have then :
\begin{eqnarray*}
((\omega_{\xi_1, \xi_2}*id)(W)^*\eta|\zeta)
&=&(\sum_i(id*\omega_{e_i, \eta})(W)^*\xi_2\underset{\nu}{_\beta\otimes_\alpha}e_i|\xi_1\underset{\nu}{_\beta\otimes_\alpha}\zeta)\\
&=&(\sum_i\alpha(<\xi_2, (id*\omega_{e_i, \eta})(W)\xi_1>_{\beta, \nu})e_i|\zeta)
\end{eqnarray*}
which finishes the proof. \end{proof}

\subsection{Lemma}
\label{lemp}
{\it In the situation of \ref{propcentral}, let $(e_i)_{i\in I}$ be an $(\alpha, \nu)$-orthogonal basis of $H$ and $J$ a finite subset of $I$; let us write $p_J=\Sigma_{i\in J}\theta^{\alpha, \nu}(e_i, e_i)$; then, for all $\Xi_1$, $\Xi_2$ in $H\underset{\nu}{_\beta\otimes_\alpha}H$, the finite sum :
\[\sum_{i\in J} ((id*\omega_{e_i, \eta})(W)\underset{N}{_\beta\otimes_\alpha}(id*\omega_{\xi, e_i})(W))\Xi_1|
\Xi_2)\]
is equal to the scalar product of the vector :
\[(\sigma_{\nu}\underset{N}{_\alpha\otimes_\alpha}1_{\gH})
(1_{\gH}\underset{N}{_\alpha\otimes_\alpha}W)
\sigma_{2\nu}
(1_{\gH}\underset{N}{_\beta\otimes_\alpha}\sigma_{\nu^o})
(1_{\gH}\underset{N}{_\beta\otimes_\alpha}(1\underset{N}{_\alpha\otimes_\alpha}p_J)W)(\Xi_1
\underset{\nu}{_\beta\otimes_\alpha}\xi)\]
with the vector $\Xi_2\underset{\nu}{_\alpha\otimes_\alpha}\eta$}

\begin{proof}
Let $\xi_1$ be in $D(H_\beta, \nu)$, $\xi_2$ in $D(_\alpha H, \nu)\cap D(H_\beta, \nu)$, $\eta_1$ in $D(_\alpha H, \nu)\cap D(H_\beta, \nu)$ and $\eta_2$ in $D(_\alpha H, \nu)$. If we take $\Xi_1=\xi_1\underset{\nu}{_\beta\otimes_\alpha}\eta_1$ and $\Xi_2=\xi_2\underset{\nu}{_\beta\otimes_\alpha}\eta_2$, the scalar product we were dealing with is equal to :
\[\sum_{i\in J}(\alpha(<(id*\omega_{e_i, \eta})(W)\xi_1, \xi_2>_{\beta, \nu})(id*\omega_{\xi, e_i})(W)\eta_1|\eta_2)\]
Using the commutativity of $N$, we consider $\alpha(<e_i, e_i>_{\alpha, \nu})(id*\omega_{\xi, e_i})(W)$, which is equal to $(id*\omega_{\xi, \alpha(<e_i, e_i>_{\alpha, \nu})e_i})(W)$, thanks to the commutation relations of $W$, and the fact that $\alpha=\hat{\beta}$. But by \ref{spatial}, we know that $\alpha(<e_i, e_i>_{\alpha, \nu})e_i=e_i$, and therefore : 
\[\alpha(<e_i, e_i>_{\alpha, \nu})(id*\omega_{\xi, e_i})(W)=(id*\omega_{\xi, e_i})(W)\]
So, this scalar product is equal to :
\[\sum_{i\in J}(\alpha(<(id*\omega_{e_i, \eta})(W)\xi_1, \xi_2>_{\beta, \nu})(id*\omega_{\xi, e_i})(W)\eta_1\underset{\nu}{_\alpha\otimes_\alpha}e_i|\eta_2\underset{\nu}{_\alpha\otimes_\alpha}e_i)\]
which is :
\[\sum_{i\in J}(\alpha(<(id*\omega_{e_i, \eta})(W)\xi_1, \xi_2>_{\beta, \nu}<(id*\omega_{\xi, e_i})(W)\eta_1, \eta_2>_{\alpha, \nu})e_i|e_i)\]
This is equal to to the scalar product of :
\[\sum_{i\in J}\alpha(<(id*\omega_{\xi, e_i})(W)\eta_1, \eta_2>_{\alpha, \nu})e_i\]
with :
\[\sum_{i\in J}\alpha(<\xi_2, (id*\omega_{e_i, \eta})(W)\xi_1>_{\beta, \nu})e_i\]
which, thanks to \ref{lemcentral}, is equal to :
\[((\omega_{\xi_1, \xi_2}*id)(W)p_J(\omega_{\eta_1, \eta_2}*id)(W)\xi|\eta)\]
Coming back to the calculations made in \ref{prop2Gamma}, we get it is equal to :
\[(W(\eta_1\underset{\nu}{_\beta\otimes_\alpha}\xi)|
\eta_2\underset{\nu}{_\alpha\otimes_\alpha}p_J(\omega_{\xi_1, \xi_2}*id)(W)^*\eta)\]
Defining now $\zeta_i$, $\zeta'_i$ as in \ref{lem1}, we get that it is equal to :
\[(W(\eta_1\underset{\nu}{_\beta\otimes_\alpha}\xi)|\eta_2\underset{\nu}{_\alpha\otimes_\alpha}\sum_i\alpha(<\zeta_i, \xi_1>_{\beta, \nu})p_J\zeta'_i)\]
Using \ref{lem2}, it is equal to the scalar product of :
\[\sigma_{2\nu}
(1_{\gH}\underset{N}{_\beta\otimes_\alpha}\sigma_{N})
(1_{\gH}\underset{N}{_\beta\otimes_\alpha}W)(\xi_1\underset{\nu}{_\beta\otimes_\alpha}\eta_1\underset{\nu}{_\beta\otimes_\alpha}\xi)\]
with :
\[\eta_2\underset{\nu}{_\alpha\otimes_\alpha}(1\underset{N}{_\alpha\otimes_\alpha}p_J)W^*(\xi_2\underset{\nu}{_\alpha\otimes_\alpha}\eta)\]
which is equal to the scalar product of :
\[\sigma_{2\nu}
(1_{\gH}\underset{N}{_\beta\otimes_\alpha}\sigma_{\nu})
(1_{\gH}\underset{N}{_\beta\otimes_\alpha}(1\underset{N}{_\alpha\otimes_\alpha}p_J)W)(\xi_1\underset{\nu}{_\beta\otimes_\alpha}\eta_1\underset{\nu}{_\beta\otimes_\alpha}\xi)\]
with :
\[\eta_2\underset{\nu}{_\alpha\otimes_\alpha}W^*(\xi_2\underset{\nu}{_\alpha\otimes_\alpha}\eta)\]
from which we get the result, by linearity, continuity and density. 
\end{proof}

\subsection{Proposition}
\label{propc}
{\it Let $(N, M, \alpha, \beta, \Gamma, T, T', \nu)$ be a quantum groupoid in the sense of Lesieur ([L1], [L2]), and let us suppose that the von Neuman algebra $\alpha(N)$ is included into the center $Z(M)$;
let $(e_i)_{i\in I}$ be an $(\alpha, \nu)$-orthogonal basis of $H$ and $J$ a finite subset of $I$; let us write $p_J=\Sigma_{i\in J}\theta^{\alpha, \tau}(e_i, e_i)$; let $k_1$, $k_2$ in $\mathcal{K}_{\alpha, \nu}$, $\xi$ in $D(_\alpha H, \nu)$, $\eta$ in $H$; then, we have :
\[lim_J\|(k_2\underset{N}{_\alpha\otimes_\alpha}(1-p_J))W(k_1\eta\underset{\nu}{_\beta\otimes_\alpha}\xi)\|=0\]}

\begin{proof}
Let $\eta_1$ in $D(_\alpha H, \nu)\cap D(H_\beta, \nu)$ and $\eta_2$ in $D(_\alpha H, \nu)$; we have :
\[R^{\alpha, \nu}(p_J(\omega_{\eta_1, \eta_2}*id)(W)\xi)=p_J(\omega_{\eta_1, \eta_2}*id)(W)R^{\alpha, \nu}(\xi)\]
and, therefore :
\begin{multline*}
<p_J(\omega_{\eta_1, \eta_2}*id)(W)\xi, p_J(\omega_{\eta_1, \eta_2}*id)(W)\xi>_{\alpha, \nu}=\\
R^{\alpha, \nu}(\xi)^*(\omega_{\eta_1, \eta_2}*id)(W)^*p_J(\omega_{\eta_1, \eta_2}*id)(W)R^{\alpha, \nu}(\xi)
\end{multline*}
which is increasing with $J$ towards 
\[<(\omega_{\eta_1, \eta_2}*id)(W)\xi, (\omega_{\eta_1, \eta_2}*id)(W)\xi>_{\alpha, \nu}\]
Let $X$ be the spectrum of $\mathbb{C}^*(\nu)$, and let us identify $\mathbb{C}^*(\nu)$ to $C_0(X)$;  using then Dini's theorem, we get it is norm converging, from which we infer that :
\[lim_J\|R^{\alpha, \nu}((1-p_J)(\omega_{\eta_1, \eta_2}*id)(W)\xi)\|=0\]
But, by \ref{lemcentral}(i), we have :
\[(1-p_J)(\omega_{\eta_1, \eta_2}*id)(W)\xi=\Sigma_{i\notin J}\alpha(<(id*\omega_{\xi, e_i})(W)\eta_1, \eta_2>_{\alpha, \nu})e_i\]
and, therefore :
\[R^{\alpha, \nu}((1-p_J)(\omega_{\eta_1, \eta_2}*id)(W)\xi)=\Sigma_{i\notin J}R^{\alpha, \nu}(e_i)<(id*\omega_{\xi, e_i})(W)\eta_1, \eta_2>_{\alpha, \nu}\]
and :
\begin{multline*}
\|R^{\alpha, \nu}((1-p_J)(\omega_{\eta_1, \eta_2}*id)(W)\xi)\|^2=\\
\|\Sigma_{i\notin J}<(id*\omega_{\xi, e_i})(W)\eta_1, \eta_2>_{\alpha, \nu}^*<(id*\omega_{\xi, e_i})(W)\eta_1, \eta_2>_{\alpha, \nu}\|
\end{multline*}
We have :
\begin{multline*}
<(id*\omega_{\xi, e_i})(W)\eta_1, \eta_2>_{\alpha, \nu}=\\R^{\alpha, \nu}(\eta_2)^*(\rho^{\alpha, \alpha}_{e_i})^*W\rho^{\beta, \alpha}_\xi R^{\alpha, \nu}(\eta_1)=\\
(\rho^{\alpha, \alpha}_{e_i})^*(R^{\alpha, \nu}(\eta_2)^*\underset{N}{_\alpha\otimes_\alpha}1)W\rho_\xi^{\beta, \alpha} R^{\alpha, \nu}(\eta_1)
\end{multline*}
and, therefore :
\begin{multline*}
\Sigma_{i\notin J}<(id*\omega_{\xi, e_i})(W)\eta_1, \eta_2>_{\alpha, \nu}^*<(id*\omega_{\xi, e_i})(W)\eta_1, \eta_2>_{\alpha, \nu}=\\
R^{\alpha, \nu}(\eta_1)^*(\rho_\xi^{\beta, \alpha})^*W^*(\theta^{\alpha, \nu}(\eta_2, \eta_2)\underset{N}{_\alpha\otimes_\alpha}(1-p_J))W\rho_\xi^{\beta, \alpha} R^{\alpha, \nu}(\eta_1)
\end{multline*}
and its norm is equal to :
\[\|(\theta^{\alpha, \nu}(\eta_2, \eta_2)\underset{N}{_\alpha\otimes_\alpha}(1-p_J))W\rho_\xi^{\beta, \alpha} R^{\alpha, \nu}(\eta_1)\|^2\]
So, we have that :
\[lim_J\|(\theta^{\alpha, \nu}(\eta_2, \eta_2)\underset{N}{_\alpha\otimes_\alpha}(1-p_J))W\rho_\xi^{\beta, \alpha} R^{\alpha, \nu}(\eta_1)\|=0\]
and, therefore :
\[lim_J\|(\theta^{\alpha, \nu}(\eta_2, \eta_2)\underset{N}{_\alpha\otimes_\alpha}(1-p_J))W\rho_\xi^{\beta, \alpha} \theta^{\alpha, \nu}(\eta_1, \eta_1)\|=0\]
and we get :
\[lim_J\|(\theta^{\alpha, \nu}(\eta_2, \eta_2)\underset{N}{_\alpha\otimes_\alpha}(1-p_J))W( \theta^{\alpha, \nu}(\eta_1, \eta_1)\eta\underset{\nu}{_\beta\otimes_\alpha}\xi)\|=0\]
from which we get the result. \end{proof}

\subsection{Theorem}
\label{thcentral1}
{\it Let $(N, M, \alpha, \beta, \Gamma, T, T', \nu)$ be a measured quantum grou-poid in the sense of Lesieur ([L1], [L2]), and let us suppose that the von Neuman algebra $\alpha(N)$ is included into the center $Z(M)$;
let $(e_i)_{i\in I}$ be an $(\alpha, \nu)$-orthogonal basis of $H$; then, we have, for all $\xi$, $\eta$ in $D(_\alpha H, \nu)$ 
\[\Gamma ((id*\omega_{\xi, \eta})(W))=\sum_i (id*\omega_{e_i, \eta})(W)\underset{N}{_\beta\otimes_\alpha}(id*\omega_{\xi, e_i})(W)\]
where the sum is weakly and strictly convergent. }
\begin{proof}
Let $\xi$, $\eta$ in $D(_\alpha H, \nu)$. Using \ref{lemp}, we get that, for all finite $J\subset I$, we have :
\[\|\sum_{i\in J} ((id*\omega_{e_i, \eta})(W)\underset{N}{_\beta\otimes_\alpha}(id*\omega_{\xi, e_i})(W))\|
\leq \|R^{\alpha, \nu}(\xi)\|\|R^{\alpha, \nu}(\eta)\|\]
Let $\xi_1$ be in $D(H_\beta, \nu)$, $\xi_2$ in $D(_\alpha H, \nu)\cap D(H_\beta, \nu)$, $\eta_1$ in $D(_\alpha H, \nu)\cap D(H_\beta, \nu)$ and $\eta_2$ in $D(_\alpha H, \nu)$; using \ref{prop2Gamma}, we get that the scalar product 
\[(\Gamma((id*\omega_{\xi, \eta})(W))(\xi_1\underset{\nu}{_\beta\otimes_\alpha}\eta_1)|
\xi_2\underset{\nu}{_\beta\otimes_\alpha}\eta_2)\]
is equal to :
\[((\omega_{\xi_1, \xi_2}*id)(W)(\omega_{\eta_1, \eta_2}*id)(W)\xi|\eta)\]
which, thanks to the calculations made in \ref{lemp}, is equal to :
\[\sum_i ((id*\omega_{e_i, \eta})(W)\xi_1\underset{\nu}{_\beta\otimes_\alpha}(id*\omega_{\xi, e_i})(W)\eta_1|
\xi_2\underset{\nu}{_\beta\otimes_\alpha}\eta_2)\]
Using then \ref{alphabeta}, we get the weak convergence.
\newline
Moreover, we get, using \ref{lemp} that, for any $k_1$, $k_2$ in $\mathcal K_{\alpha, \nu}$ :
\begin{multline*}
|(\Sigma_{i\notin J}(id*\omega_{e_i, \eta})(W)\xi_1\underset{\nu}{_\beta\otimes_\alpha}k_1^*(id*\omega_{\xi, e_i})(W)k_2\eta_1|\xi_2\underset{\nu}{_\beta\otimes_\alpha}\eta_2)|\leq\\
\|(k_2\underset{N}{_\alpha\otimes_\alpha}(1-p_J))W(k_1\eta\underset{\nu}{_\beta\otimes_\alpha}\xi)\|
\|\xi_1\underset{\nu}{_\beta\otimes_\alpha}\eta_1\|\|\xi_2\underset{\nu}{_\beta\otimes_\alpha}\eta_2\|
\end{multline*}
from which, thanks to \ref{propc}, we get that :
\[lim_J\|(\Sigma_{i\notin J}(id*\omega_{e_i, \eta})(W)\underset{N}{_\beta\otimes_\alpha}k_1^*(id*\omega_{\xi, e_i})(W)k_2\|=0\]
which gives the result.  \end{proof}

\subsection{Theorem}
\label{thcentral2}
{\it Let $(N, M, \alpha, \beta, \Gamma, T, T', \nu)$ be a measured quantum grou-poid in the sense of Lesieur ([L1], [L2]), and let us suppose that the von Neuman algebra $\alpha(N)$ is included into the center $Z(M)$;
then, for all $x$ in the $\mathbb{C}^*$-algebra $A_n(W)$, $\Gamma (x)$ belongs to the multipliers of the  $\mathbb{C}^*$-algebra generated by all elements of the form $a\underset{N}{_\beta\otimes_\alpha}b$, where $a$, $b$ belong to $A_n(W)$; using \ref{cortensor}, we get that $\Gamma (x)\in M(A_n(W)\underset{\mathbb{C}^*(\nu)}{_\beta\otimes_\alpha^m}A_n(W))$. }

\begin{proof}
Let $\xi$, $\eta$ in $D(_\alpha H, \nu)$; using \ref{thcentral1}, the operator $\Gamma ((id*\omega_{\xi, \eta})(W))$ is a strict limit of elements in $A_n(W)\underset{N}{_\beta\otimes_\alpha}A_n(W)$, and therefore belongs to $M(A_n(W)\underset{N}{_\beta\otimes_\alpha}A_n(W))$, from which we get the result, by definition of $A_n(W)$. \end{proof}

\subsection{Proposition}
\label{thgroupoid}
{\it Let $(N, M, \alpha, \beta, \Gamma, T, T', \nu)$ be a quantum groupoid in the sense of Lesieur ([L1], [L2]), and let us suppose that the von Neuman algebra $M$ is abelian; let us write $\mathcal G$ for the spectrum of the $\mathbb{C}^*$-algebra $A_n(W)$, and $\mathcal G^{(0)}$ for the spectrum of the $\mathbb{C}^*$-algebra $\mathbb{C}^*(\nu)$. Then :
\newline
(i) there exists a continuous open application $r$ from $\mathcal G$ onto $\mathcal G^{(0)}$, such that, for all $f\in C_0(\mathcal G^{(0)})=A$, we have $\alpha (f)=f\circ r$; there exists a continous open application $s$ from $\mathcal G$ onto $\mathcal G^{(0)}$, such that, for all $f\in C_0(\mathcal G^{(0)})=A$, we have $\beta (f)=f\circ s$. 
\newline
(ii) then, there exists a partially defined multiplication on $\mathcal G$, which gives to $\mathcal G$ a structure of locally compact groupoid, with $\mathcal G ^{(0)}$ as set of units. 
\newline
(iii) The application defined for all $F$ continuous, positive, with compact support in $\mathcal G$, by $F\mapsto \alpha^{-1}(T(F))(u)$, defines a positive Radon measure $\lambda^u$ on $\mathcal G$, whose support is $\mathcal G^u$. The measures $(\lambda^u)_{u\in \mathcal G ^{(0)}}$ are a Haar system on $\mathcal G$.
\newline
(iv) the trace $\nu$ on $\mathbb{C}^*(\nu)$ leads to a quasi-invariant measure (denoted again by $\nu$) on $\mathcal G ^{(0)}$. If we write $\mu$ for the measure on $\mathcal G$ constructed from $\nu$ and the Haar system, we have then : 
\[(N, M, \alpha, \beta, \Gamma)=(L^{\infty}(\mathcal G ^{(0)}, \nu), L^{\infty}(\mathcal G , \mu), r_{\mathcal G}, s_{\mathcal G}, \Gamma_{\mathcal G})\]
where $r_{\mathcal G}$, $s_{\mathcal G}$, $\Gamma_{\mathcal G}$ have been defined in \ref{Hbimod}. Moreover, then, the operator-valued weights are given, for any positive $F$ in $L^{\infty}(\mathcal G , \nu)$ by :
\[T(F)(u)=\int_{\mathcal G}Fd\lambda^u\]
\[T'(F)(u)=\int_{\mathcal G}Fd\lambda_u\]
where $\lambda_u$ is the image of $\lambda^u$ under the application $(x\mapsto x^{-1})$.}

\begin{proof}
As $\alpha(N)\subset M(A_n(W))$, we can construct by restriction a continuous application $r$ from $\mathcal G$ into $\mathcal G^{(0)}$ such that, for all $f\in C_0(\mathcal G^{(0)})=\mathbb{C}^*(\nu)$, we have $\alpha (f)=f\circ r$; we can construct the same way a continuous application $s$ from $\mathcal G$ into $\mathcal G^{(0)}$ such that, for all $f\in C_0(\mathcal G^{(0)})=\mathbb{C}^*(\nu)$, we have $\beta (f)=f\circ r$. The applications $r$ and $s$ are open by ([B2], 3.14), which gives (i).
\newline
The application $R$ from $A_n(W)$ into itself leads to an involutive application in $\mathcal G$, we shall write $x\mapsto x^{-1}$, and, using that $R\circ\alpha=\beta$, we get that $r(x^{-1})=s(x)$ and $s(x^{-1})=r(x)$.
\newline
Thanks to \ref{continuousfield}, we may apply \ref{tensor} to $A_n(W)$, which we identify to $C_0(\mathcal G)$, and we obtain that the commutative $\mathbb{C}^*$-algebra $A_n(W)\underset{\mathbb{C}^*(\nu)}{_\beta\otimes^m_\alpha}A_n(W)$ is the quotient of $A_n(W)\otimes A_n(W)$ (identified with $C_0(\mathcal G^2)$) by the ideal generated by all the functions $(x_1, x_2)\mapsto f(s(x_1))g(x_1, x_2)-f(r(x_2))g(x_1, x_2)$, where $x_1$, $x_2$ are in $\mathcal G$, $f$ in $\mathbb{C}^*(\nu)$ (identified with $C_0(\mathcal G^{(0)})$), and $g$ in $C_0(\mathcal G^2)$. So, a non zero character on $A_n(W)\underset{\mathbb{C}^*(\nu)}{_\beta\otimes^m_\alpha}A_n(W)$ is a couple $(x_1, x_2)$ in $\mathcal G^2$ such that $s(x_1)=r(x_2)$; let us write $\mathcal G^{(2)}$ for the subset of such elements of $\mathcal G^2$. 
\newline
So, we shall identify $A_n(W)\underset{\mathbb{C}^*(\nu)}{_\beta\otimes^m_\alpha}A_n(W)$ to $C_0(\mathcal G^{(2)})$. Therefore, we see that the restriction of $\Gamma$ to $A_n(W)$ leads to a continuous application from $\mathcal G^{(2)}$ into $\mathcal G$, which gives to $\mathcal G$ a structure of locally compact groupoid, which is (ii). 
\newline
As $A_n(W)\cap\gM_T$ is a dense ideal in $A_n(W)$, it contains the ideal $\mathcal K(\mathcal G)$ of continuous functions on $\mathcal G$, with compact support; for all $F$ in $\mathcal K(\mathcal G)$, $\alpha^{-1}(T(F))$ belongs to $C_b(\mathcal G^{(0)})$, and, for all $u\in\mathcal G^{(0)}$, $F\mapsto\alpha^{-1}(T(F))(u)$ defines a non zero positive Radon measure $\lambda^u$ on $\mathcal G$; it is now straightforward to get, from the left invariance of $T$, that $(\lambda^u)_{u\in \mathcal G^{(0)}}$ is a Haar system on the groupoid. Starting from $R\circ T\circ R$, we obtain measures $\lambda_u$, which are the images of $\lambda^u$ by the inverse. 
\newline
The modulus $\delta$ of the measured quantum groupoid gives that the trace $\nu$ on $\mathbb{C}^*(\nu)$ leads to a quasi-invariant measure $\mu$ on $\mathcal G^{(0)}$. 
\newline
Now, by density reasons, we shall identify $N$ with $L^{\infty}(\mathcal G^{(0)}, \mu)$, $M$ with $L^{\infty}(\mathcal G, \mu)$, where $\mu$ is the measure on $\mathcal G$ constructed from $\mu$ and the Haar system, $\alpha$ with $r_{\mathcal G}$, $\beta$ with $s_{\mathcal G}$, $\Gamma$ with $\Gamma_{\mathcal G}$, and we obtain the required formulae for the left and right Haar systems.
\end{proof}

\subsection{Ramsay's theorem}
\label{ramsay}
{\it Let $\mathcal G$ be a measured groupo\"{\i}d, with $\mathcal G^{(0)}$ as space
of units, and $r$ and $s$ the range and source functions from $\mathcal G$ to $\mathcal G^{(0)}$, with a Haar system $(\lambda^u)_{u\in \mathcal G^{(0)}}$ and a quasi-invariant measure $\mu$ on $\mathcal G^{(0)}$. Let us write $\nu$ the associated measure on $\mathcal G$. Let $\Gamma_{\mathcal G}$, $r_{\mathcal G}$, $s_{\mathcal G}$ be the morphisms associated in \ref{gd}. Then, there exists a locally compact groupoid $\tilde{\mathcal G}$, with set of units $\tilde{\mathcal G}^{(0)}$, with a Haar system $(\tilde{\lambda}^u)_{u\in \tilde{\mathcal G}^{(0)}}$, and a quasi-invariant measure $\tilde{\mu}$ on $\tilde{\mathcal G}^{(0)}$, such that, if $\tilde{\nu}$ is the measure on $\tilde{\mathcal G}$ constructed from $\tilde{\mu}$ and $(\tilde{\lambda}^u)_{u\in \tilde{\mathcal G}^{(0)}}$, we get that $L^{\infty}(\mathcal G^{(0)}, \mu)$ is isomorphic to $L^{\infty}(\tilde{\mathcal G}^{(0)}, \tilde{\mu})$, $L^{\infty}(\mathcal G, \nu)$ is isomorphic to $L^{\infty}(\tilde{\mathcal G}, \tilde{\nu})$, these isomorphisms sends $\mu$ on $\tilde{\mu}$, $\nu$ on $\tilde{\nu}$, $r_{\mathcal G}$ and $s_{\mathcal G}$ on $r_{\tilde{\mathcal G}}$ and $s_{\tilde{\mathcal G}}$, $\Gamma_{\mathcal G}$ on $\Gamma_{\tilde{\mathcal G}}$, and the operator-valued weight constructed with $(\lambda^u)_{u\in \mathcal G^{(0)}}$ on the operator-valued weight constructed with $(\tilde{\lambda}^u)_{u\in \tilde{\mathcal G}^{(0)}}$.}
\begin{proof}
Let us apply \ref{thgroupoid} to the commutative measured quantum groupoid constructed from the measured groupoid $\mathcal G$. Then, we get the result. \end{proof}

\section{$\mathbb{C}^*$-measured quantum groupoids (when the basis is central)}
\label{C*central}

In this chapter, we define $\mathbb{C}^*$ measured quantum groupoids, when the basis is central. In that case, we can use the notion of measured continuous fields of $\mathbb{C}^*$-algebras, defined in \ref{defmcf}. We need then to define appropriate notion of coproduct (\ref{defCXHopf}) to define $C_0(X)$-Hopf algebras; we define also left (resp. right) invariant $\mathbb{C}^*$-weight (\ref{defCXHopf}), and we give a definition for $\mathbb{C}^*$ quantum groupoids (\ref{defC*quantum}); starting from a measured quantum groupoid, with a central basis, we obtain that the $\mathbb{C}^*$-algebra $A_n(W)$ has a natural structure of a $\mathbb{C}^*$ quantum groupoids (\ref{AC*quantum}); more precisely, we obtain that it is, in a certain sense, the greatest one (\ref{thU2}(v)); conversely, starting from a $\mathbb{C}^*$ quantum groupoid, we obtain a measured quantum groupoid (\ref{thU2}). 

\subsection{Definitions and notations}
\label{defnot}
Let $(X, \alpha, A, \varphi)$ and $(X, \beta, A, \psi)$ be two measured continuous fields of $\mathbb{C}^*$-algebras, in the sense of \ref{defmcf}; we have seen in chapter \ref{CXC*} that it is possible to consider the minimal tensor product of the $C_0(X)$-algebras $A$ (via $\beta$) with $A$ (via $\alpha$), which will be denoted $A\underset{C_0(X)}{_\beta\otimes^m_\alpha}A$. 
\newline
This $\mathbb{C}^*$ algebra is clearly, via the morphisms $f\mapsto \alpha(f)\underset{C_0(X)}{_\beta\otimes^m_\alpha}1$ and $f\mapsto 1\underset{C_0(X)}{_\beta\otimes^m_\alpha}\beta(f)$, a $C_0(X)$-algebra, in two different ways, and we can consider the minimal tensor product $(A\underset{C_0(X)}{_\beta\otimes^m_\alpha}A)\underset{C_0(X)}{_\beta\otimes^m_\alpha}A$, which is, as $A$ is a continuous field of $\mathbb{C}^*$-algebras, using ([B1], 3.1), given by the minimal $\mathbb{C}^*$-norm on the algebra obtained by taking the quotient of the algebraic tensor product of $A\underset{C_0(X)}{_\beta\otimes^m_\alpha}A$ with $A$, by the ideal generated by the elements of the form \[\{X(1\underset{C_0(X)}{_\beta\otimes^m_\alpha}\beta(f))\otimes a-X\otimes\alpha(f)a, X\in A\underset{C_0(X)}{_\beta\otimes^m_\alpha}A, a\in A, f\in C_0(X)\}\]
We can also consider the $\mathbb{C}^*$-algebra $A\underset{C_0(X)}{_\beta\otimes^m_\alpha}(A\underset{C_0(X)}{_\beta\otimes^m_\alpha}A)$, and, using ([B1], end of chapter 4), we get that these two $\mathbb{C}^*$-algebras are equal, and shall be denoted by $A\underset{C_0(X)}{_\beta\otimes^m_\alpha}A\underset{C_0(X)}{_\beta\otimes^m_\alpha}A$. 
\newline
Let now consider an injective $*$-homomorphism $\Gamma$ from $A$ into the multiplier algebra $M(A\underset{C_0(X)}{_\beta\otimes^m_\alpha}A)$, and let us denote again $\Gamma$ its extension to $M(A)$. Let us suppose that, for any $f$ in $C_0(X)$, we have :
\[\Gamma(\alpha(f))= \alpha(f)\underset{C_0(X)}{_\beta\otimes^m_\alpha}1\]
\[\Gamma(\beta(f))=1\underset{C_0(X)}{_\beta\otimes^m_\alpha}\beta(f)\]
and let us consider the $*$-homomorphism $\Gamma\otimes id$ from the algebraic tensor product $A\otimes A$ in the $\mathbb{C}^*$-algebra $M(A\underset{C_0(X)}{_\beta\otimes^m_\alpha}A\underset{C_0(X)}{_\beta\otimes^m_\alpha}A)$, which contains $M(A\underset{C_0(X)}{_\beta\otimes^m_\alpha}A)\underset{C_0(X)}{_\beta\otimes^m_\alpha}A$.  We have, for any $a_1$, $a_2$ in $A$ :
\begin{multline*}
(\Gamma\otimes id)(\beta(f)a_1\otimes a_2-a_1\otimes\alpha(f)a_2)=\\
(1\underset{C_0(X)}{_\beta\otimes^m_\alpha}\beta(f))\Gamma(a_1)\underset{C_0(X)}{_\beta\otimes^m_\alpha} a_2-\Gamma(a_1)\underset{C_0(X)}{_\beta\otimes^m_\alpha}\alpha(f)a_2=0
\end{multline*}
and, therefore, this application goes to the quotient, and, using the min property, to the closure $A\underset{C_0(X)}{_\beta\otimes^m_\alpha}A$, and we shall denote it then $\Gamma\underset{C_0(X)}{_\beta\otimes^m_\alpha}id$.
\newline
The definition of $id\underset{C_0(X)}{_\beta\otimes^m_\alpha}\Gamma$ is done the same way. 
\newline
For any $a$ in $\gN_{\varphi}$, let now (for $b\in A$) $\omega_a(b)$ be the image by $\alpha$ of the function $x\mapsto\varphi(a^*ba)$;  $\omega_a$ is a continuous field of positive linear forms on the $C_0(X)$-algebra $A$ (via $\alpha$) defined in \ref{defmcf}. Let $\nu$ be a poitive measure on $X$, whose support is $X$; let us denote $\Phi=\int_X\varphi^xd\nu(x)$ be the lower semi-continuous semi-finite KMS faithful weight on $A$ constructed in \ref{defmcf}; let $\overline{\alpha}$ be the representation of $L^\infty(X, \nu)$ on $H_\Phi$ and $T$ be the normal faithful semi-finite operator valued weight from $\pi_\phi(A)''$ to $\overline{\alpha}(L^\infty(X,\nu))$ as constructed in \ref{propfieldweights}; as $\omega_a(b)=\Lambda_T(a)^*\pi_{\varphi}(b)\Lambda_{T}(a)$, it is possible, using  \ref{propfieldweights} and ([B1], 4.1), to define $(id\underset{C_0(X)}{_\beta\otimes^m_\alpha}\omega_a)$ from $A\underset{C_0(X)}{_\beta\otimes^m_\alpha}A$ on $A$, and, by extension, from $M(A\underset{C_0(X)}{_\beta\otimes^m_\alpha}A)$ on $M(A)$. 
\newline
For any $b$ in $\gN_{\psi}$, we define as well $\omega_b(a)$ as the image by $\beta$ of the function $x\mapsto\psi^x(b^*ab)$ and $(\omega_b\underset{C_0(X)}{_\beta\otimes^m_\alpha}id)$ from  $A\underset{C_0(X)}{_\beta\otimes^m_\alpha}A$ to $A$, or from $M(A\underset{C_0(X)}{_\beta\otimes^m_\alpha}A)$ to $M(A)$.

\subsection{Definition}
\label{defCXHopf}
Let $X$ be a locally compact space, $A$ a $\mathbb{C}^*$-algebra, $\alpha$ (resp. $\beta$), a non degenerate $*$-homomorphism from $C_0(X)$ into $Z(M(A))$, $\varphi$ (resp. $\psi$) such that $(X, \alpha, A, \varphi)$ (resp. $(X, \beta, A, \psi)$) is a measured continuous field of $\mathbb{C}^*$-algebra, as defined in \ref{defmcf}. Let us consider the minimal tensor product $A\underset{C_0(X)}{_\beta\otimes^m_\alpha}A$ as constructed in \ref{defnot}; let now $\Gamma$ be a non degenerate $*$-homomorphism from $A$ into $M(A\underset{C_0(X)}{_\beta\otimes^m_\alpha}A)$, such that, for any $f$ in $C_0(X)$, we have :
\[\Gamma(\alpha(f))= \alpha(f)\underset{C_0(X)}{_\beta\otimes^m_\alpha}1\]
\[\Gamma(\beta(f))=1\underset{C_0(X)}{_\beta\otimes^m_\alpha}\beta(f)\]
and let us suppose that :
\[(\Gamma\underset{C_0(X)}{_\beta\otimes^m_\alpha}id)\Gamma=
(id\underset{C_0(X)}{_\beta\otimes^m_\alpha}\Gamma)\Gamma\]
where $\Gamma\underset{C_0(X)}{_\beta\otimes^m_\alpha}id$ and $id\underset{C_0(X)}{_\beta\otimes^m_\alpha}\Gamma$ had been defined in \ref{defnot}. Such an application will be called a $C_0(X)$-coproduct. 
\newline
Moreover, we shall suppose that, for all $c$ in $A$, for all $a$ in $\gN_{\varphi}$, and $b$ in $\gN_{\psi}$, the elements $(id\underset{C_0(X)}{_\beta\otimes^m_\alpha}\omega_a)\Gamma (c)$
and $(\omega_b\underset{C_0(X)}{_\beta\otimes^m_\alpha}id)(\Gamma (c))$ belong to $A$ (and not only $M(A)$). 
\newline
Such a $7$-uple $(X, \alpha, \beta, A, \Gamma, \varphi, \psi)$ will be called a $C_0(X)$-Hopf algebra.

\subsection{Example}
\label{gd3}
Let $\mathcal G$ be a locally compact groupoid, as defined in \ref{gd}. We have seen in \ref{gd2} that $C_0(\mathcal G)$ is, in two ways (via $r$ and via $s$) a continuous field of $\mathbb{C}^*$-algebras over $\mathcal G^{(0)}$, and that the tensor product $C_0(\mathcal G)\underset{C_0(\mathcal G^{(0)})}{_{\hat{s}}\otimes^m_{\hat{r}}}C_0(\mathcal G)$ can be identified with $C_0(\mathcal G^{(2)})$. So, the product of the groupoid leads to a coproduct $\Gamma_{\mathcal G}$ in the sense of \ref{defCXHopf}. 
\newline
If $f$ belongs to $C_0(\mathcal G)$, and if $g$ belongs to $C_0(\mathcal G)\cap (\underset{u\in\mathcal G^{(0)}}{\cap}\L^2(\mathcal G, \lambda^u))$, let us consider the function on $\mathcal G$ :
\[s\mapsto \int_{\mathcal G^{(0)}}d\mu(u)\int_{G^u}f(st)|g|^2(t)d\lambda^u(t)\]
If both $f$ and $g$ belongs to $\mathcal K(\mathcal G)$, this function belongs also to $\mathcal K(\mathcal G)$, and, therefore, by continuity, it belongs, in general, to $C_0(\mathcal G)$. So, we get that $(\mathcal G^{(0)}, \hat{r}, \hat{s}, C_0(\mathcal G), \Gamma_{\mathcal G}, \lambda^u, \lambda_u)$ is a $C_0(\mathcal G^{(0)})$-Hopf algebra\vspace{2mm}.  
 \newline
 Conversely, if $A$ is an abelian $\mathbb{C}^*$-algebra, and $(X, \alpha, \beta, A, \Gamma, \varphi, \psi)$ is a $C_0(X)$-Hopf algebra, let us identify $A$ with $C_0(Y)$, where $Y$ is the spectrum of $A$; then, by transposition, $\alpha$ (resp. $\beta$) gives a continuous application from $Y$ to $X$, which is open by ([B2], 3.14); writing $Y^{(2)}$ the subset of $Y^2$ made of couples $(x, y)$ such that $s(x)=r(y)$, we may identify $A\underset{C_0(X)}{_\beta\otimes^m_\alpha}A$ with $C_0(Y^{(2)})$, and, by  transposition of $\Gamma$, we get a continuous product from $Y^{(2)}$ to $Y$ which gives to $Y$ a structure of locally compact groupoid. 
 
 \label{inv}
 \subsection{Definition}
Let $(X, \alpha, \beta, A, \Gamma, \varphi, \psi)$ be a $C_0(X)$-Hopf algebra. Then, $\varphi$ will be said left-invariant with respect to $\Gamma$ if, for all $a\in \gM_{\varphi}$, all $b\in\gN_{\psi}$, we have, for all $x\in X$ :
\[\varphi^x((\omega_b\underset{C_0(X)}{_\beta\otimes^m_\alpha}id)(\Gamma (a))=\psi^x(b^*b)\varphi^x(a)\]
and $\psi$ will be said right-invariant with respect to $\Gamma$ if, for all $b\in\gM_{\psi}$, all $a\in\gN_{\varphi}$, we have :
\[\psi^x((id\underset{C_0(X)}{_\beta\otimes^m_\alpha}\omega_a)\Gamma (b))=\varphi^x(a^*a)\psi^x(b)\]

Moreover, $\Gamma$ will be said simplifiable if :
\newline
(i) $A$ is equal to the norm closure of the subspace generated by all the elements of the form $(id\underset{C_0(X)}{_\beta\otimes^m_\alpha}\omega_a)\Gamma (x)$, where $x$ is in $A$ and $a$ is in $\gN_{\varphi}$; 
\newline
(ii) $A$ is equal to the norm closure of the subspace generated by all the elements of the form $(\omega_b\underset{C_0(X)}{_\beta\otimes^m_\alpha}id)(\Gamma (x))$, where $x$ is in $A$ and $b$ is in $\gN_{\psi}$\vspace{2mm}. 
\newline
Coming back to the case of a locally compact groupoid (\ref{gd3}), it is straightforward to verify that the Haar system $\lambda^u$ (which gives a continuous field of weights by \ref{gd2}) is left-invariant, and that $\lambda_u$ is right-invariant. 
Taking now an approximate unit ([R1], 1.9), one can prove that any function $f$ in $\mathcal K(\mathcal G)$ belongs to norm closure of the subspace generated by all elements of the form $(id\underset{C_0(\mathcal G^{(0)})}{_{\hat{s}}\otimes^m_{\hat{r}}}\omega_g)\Gamma_{\mathcal G}(f)$, where $g$ is in $\mathcal K(\mathcal G)$; from which we deduce that $\Gamma_{\mathcal G}$ is simplifiable. 

\subsection{Lemma}
\label{lemU}
{\it Let $(X, \alpha, \beta, A, \Gamma, \varphi, \psi)$ be a $C_0(X)$-Hopf algebra, and let us suppose that $\varphi$ is left-invariant. Let $\nu$ be a positive Radon measure on $X$, with support $X$. Let us write $\Phi=\int_X\varphi^xd\nu(x)$, $\Psi=\int_X\psi^xd\nu(x)$, $\overline{\alpha}$ (resp. $\overline{\beta}$) the representation of $L^\infty (X, \nu)$ on $H_{\Phi}$ (resp. $H_{\Psi}$) which extends $\alpha$ (resp. $\beta$). Let us write $\overline{\Phi}$ (resp. $\overline{\Psi}$) for the normal semi-finite faithful weight on $\pi_\Phi(A)''$ (resp. $\pi_\Psi(A)''$) which extends $\Phi$ (resp. $\Psi$), and $T$ (resp. $S$) for the normal semi-finite faithful operator-valued weight from $\pi_\Phi (A)''$ to $\overline{\alpha}(L^\infty(X, \nu))$ (resp. from $\pi_\Psi(A)''$ to $\overline{\beta}(L^\infty(X, \nu))$) which extends $\varphi$ (resp. $\psi$). We shall consider $A\underset{C_0(X)}{_\beta\otimes^m_\alpha}A$ as a $\mathbb{C}^*$-algebra on $H_{\Psi}\underset{\nu}{_{\overline{\beta}}\otimes_{\overline{\alpha}}}H_{\Phi}$, as explained in \ref{tensor}. Then : 
\newline
(i) for any $a$ in $\gN_{\varphi}\cap\gN_{\Phi}$, $x$ in $\gN_{\psi}\cap\gN_{\Psi}$, $\xi$ in $D((H_{\Psi})_{\overline{\beta}}, \nu^o)$, the operator 
\[(\omega_{J_{\Psi}\Lambda_{\Psi}(x),\xi}\underset{N}{_{\overline{\beta}}*_{\overline{\alpha}}}id)\Gamma (a)\]
 belongs to $\gN_{T}\cap\gN_{\overline{\Phi}}$. 
 \newline
 (ii) for any $a$ in $\gN_{\varphi}\cap\gN_{\Phi}$, $x$ in $\gN_{\psi}\cap\gN_{\Psi}$, and $(\xi_i)_{i\in I}$ an $(\overline{\beta}, \nu)$-orthogonal basis of $H_{\Psi}$, the sum :
 \[\sum_{i\in I}\xi_i\underset{\nu}{_{\overline{\beta}}\otimes_{\overline{\alpha}}}\Lambda_{\overline{\Phi_1}}((\omega_{J_{\Phi_2}\Lambda_{\Phi_2}(x),\xi_i}\underset{N}{_{\overline{\beta}}*_{\overline{\alpha}}}id)\Gamma (a))\]
 converges in $H_{\Psi}\underset{\nu^o}{_{\overline{\beta}}\otimes_{\overline{\alpha}}}H_{\Phi}$, and it does not depend on the $(\overline{\beta}, \nu^o)$-orthogonal basis of $H_{\Psi}$. }
 \begin{proof}
 The operator $X=(\omega_{J_{\Psi}\Lambda_{\Psi}(x),\xi}\underset{N}{_{\overline{\beta}}*_{\overline{\alpha}}}id)\Gamma (a)$ belongs to $\pi_{\Phi}(A)''$. By definition, we have :
\begin{eqnarray*}
X^*X
&=&(\lambda^{\overline{\beta}, \overline{\alpha}}_{J_{\Psi}\Lambda_{\Psi}(x)})^*\Gamma (a^*)
\lambda^{\overline{\beta}, \overline{\alpha}}_{\xi}(\lambda^{\overline{\beta}, \overline{\alpha}}_{\xi})^*\Gamma(a)\lambda^{\overline{\beta}, \overline{\alpha}}_{J_{\Psi}\Lambda_{\Psi}(x)}\\
&\leq&\|R^{\overline{\beta},\nu}(\xi)\|^2
(\omega_{J_{\Psi}\Lambda_{\Psi}(x)}\underset{N}{_{\overline{\beta}}*_{\overline{\alpha}}}id)\Gamma (a^*a)\\
&=&\|R^{\overline{\beta},\nu}(\xi)\|^2
(\Lambda_{S}(x)\underset{\nu}{\otimes_{\overline{\alpha}}}id)^*\Gamma(a^*a)(\Lambda_{S}(x)\underset{\nu}{\otimes_{\overline{\alpha}}}id)\\
&=&\|R^{\overline{\beta}, \nu}(\xi)\|^2(\omega_x\underset{C_0(X)}{_\beta\otimes^m_\alpha}id)\Gamma (a^*a)
\end{eqnarray*}
and, therefore :
\begin{eqnarray*}
T(X^*X)
&\leq&\|R^{\overline{\beta}, \nu}(\xi)\|^2
T((\omega_x\underset{C_0(X)}{_\beta\otimes^m_\alpha}id)\Gamma (a^*a))\\
&=&\|R^{\overline{\beta}, \nu}(\xi)\|^2T(a^*a)\alpha\circ\beta^{-1}(S(b^*b))
\end{eqnarray*}
and $\overline{\Phi}(X^*X)\leq \|R^{\overline{\beta}, \nu}(\xi)\|^2
\int_X\varphi^x(a^*a)\psi^x(b^*b)d\nu(x)$, which proves (i).
\newline 
Using (i), part (ii) is proved the same way as ([L2], 3.13).
\end{proof}
 
\subsection{Proposition}
\label{propU}
{\it Let $(X, \alpha, \beta, A, \Gamma, \varphi, \psi)$ be a $C_0(X)$-Hopf algebra, and let us suppose that $\varphi$ is left-invariant. Let $\nu$ be a positive Radon measure on $X$, with support $X$. Let us write $\Phi=\int_X\varphi^xd\nu(x)$, $\Psi=\int_X\psi^xd\nu(x)$, $\overline{\alpha}$ (resp. $\overline{\beta}$) the representation of $L^\infty (X, \nu)$ on $H_{\Phi}$ (resp. $H_{\Psi}$) which extends $\alpha$ (resp. $\beta$). Let us write $\overline{\Phi}$ (resp. $\overline{\Psi}$) for the normal semi-finite faithful weight on $\pi_\Phi(A)''$ (resp. $\pi_\Psi(A)''$) which extends $\Phi$ (resp. $\Psi$), and $T$ (resp. $S$) for the normal semi-finite faithful operator-valued weight from $\pi_\Phi (A)''$ to $\overline{\alpha}(L^\infty(X, \nu))$ (resp. from $\pi_\Psi(A)''$ to $\overline{\beta}(L^\infty(X, \nu))$) which extends $\varphi$ (resp. $\psi$). We shall consider $A\underset{C_0(X)}{_\beta\otimes^m_\alpha}A$ as a $\mathbb{C}^*$-algebra on $H_{\Psi}\underset{\nu}{_{\overline{\beta}}\otimes_{\overline{\alpha}}}H_{\Phi}$, as explained in \ref{tensor}. Let us suppose, moreover, that the representation $\pi_{\Psi}\circ\alpha$ of $C_0(X)$ is square-integrable with respect to $\nu$, and let us write $\overline{\alpha_2}$ the representation of $L^\infty(X, \nu)$ on $H_{\Psi}$ constructed then. On the other hand, let us write, for all $f\in L^\infty(X, \nu)$ 
\[\hat{\beta}(f)=J_{\overline{\Phi}}\overline{\alpha}(f^*)J_{\overline{\Phi}}\]
Then :
\newline
(i) there exists a unique isometry $U$ from $H_{\Psi}\underset{\nu}{_{\overline{\alpha_2}}\otimes_{\hat{\beta}}}H_{\Phi}$ into $H_{\Psi}\underset{\nu}{_{\overline{\beta}}\otimes_{\overline{\alpha}}}H_{\Phi}$ such that, for any $x\in\gN_{\psi}\cap\gN_{\Psi}$, $a\in\gN_{\varphi}\cap\gN_{\Phi}$, and $(\xi_i)_{i\in I}$ orthonormal $(\overline{\beta}, \nu)$-basis of $H_{\Psi}$, we have :
\[U(J_{\overline{\Psi}}\Lambda_{\Psi}(x)\underset{\nu}{_{\overline{\alpha_2}}\otimes_{\hat{\beta}}}\Lambda_{\Phi}(a))=
\sum_{i\in I}\xi_i\underset{\nu}{_{\overline{\beta}}\otimes_{\overline{\alpha}}}
\Lambda_{\overline{\Phi}}((\omega_{J_{\Psi}\Lambda_{\Psi}(x),\xi_i}\underset{L^\infty(X, \nu)}{_{\overline{\beta}}*_{\overline{\alpha}}}id)\Gamma (a))\]
(ii)  for all $a$ in $\gN_{\varphi}\cap\gN_{\Phi}$, $x$ in $\gN_{\psi}\cap\gN_{\Psi}$, $\xi$ in $D((H_{\Psi})_{\overline{\beta}}, \nu)$, we have 
\[\Lambda_{\overline{\Phi}}((\omega_{J_{\overline{\Psi}}\Lambda_{\Psi}(x),\xi}\underset{L^\infty(X, \nu)}{_{\overline{\beta}}*_{\overline{\alpha}}}id)\Gamma (a))=(\lambda^{\overline{\beta}, \overline{\alpha}}_\xi)^*U(J_{\overline{\Psi}}\Lambda_{\Psi}(x)\underset{\nu}{_{\overline{\alpha_2}}\otimes_{\hat{\beta}}}\Lambda_{\Phi}(a))\]
(iii) for all $a_1$, $a_2$ in $\gN_{\varphi}\cap\gN_{\Phi}$, we have :
\[(1\underset{L^\infty(X, \nu)}{_{\overline{\beta}}\otimes_{\overline{\alpha}}}J_{\overline{\Phi}}\pi_{\Phi}(a_2)J_{\overline{\Phi}})U\rho^{\overline{\alpha_2}, \hat{\beta}}_{\Lambda_{\Phi}(a_1)}=
\Gamma(a_1)\rho^{\overline{\beta}, \overline{\alpha}}_{J_{\overline{\Phi}}\Lambda_{\Phi}(a_2)}\]
and, for any $x$ in $A$, we have :
\[\Gamma (x)U=U(1\underset{L^\infty(X, \nu)}{_{\overline{\alpha_2}}\otimes_{\hat{\beta}}}\pi_{\Phi}(x))\] 
(iv) for all $f\in C_0(X)$, we have :
\[U(1 \underset{L^\infty(X, \nu)}{_{\overline{\alpha_2}}\otimes_{\hat{\beta}}}\pi_{\Phi}\circ\alpha (f))=
(\pi_{\Psi}\circ\alpha (f)\underset{L^\infty(X, \nu)}{_{\overline{\beta}}\otimes_{\overline{\alpha}}}1)U\]
\[U(1 \underset{L^\infty(X, \nu)}{_{\overline{\alpha_2}}\otimes_{\hat{\beta}}}\pi_{\Phi}\circ\beta (f))=
(1\underset{L^\infty(X, \nu)}{_{\overline{\beta}}\otimes_{\overline{\alpha}}}\pi_{\Phi}\circ\beta (f))U\]
\[U(\pi_{\Psi}\circ\beta (f)\underset{L^\infty(X, \nu)}{_{\overline{\alpha_2}}\otimes_{\hat{\beta}}}1)=
((1\underset{L^\infty(X, \nu)}{_{\overline{\beta}}\otimes_{\overline{\alpha}}}\hat{\beta}(f))U\]

 }

\begin{proof}
Using \ref{lemU}(ii), result (i)  is proved the same way as ([L2], 3.14); then result (ii) is straightforward, and, using (i) and (ii), results (iii), and (iv) are proved as ([L2], 3.16), ([L2], 3.25) and ([L2], 3.22).

\end{proof}

\subsection{Theorem}
\label{thU}
{\it Let $(X, \alpha, \beta, A, \Gamma, \varphi, \psi)$ be a $C_0(X)$-Hopf algebra, and let us suppose that $\varphi$ is left-invariant. Let $\nu$ be a positive Radon measure on $X$, with support $X$. Let us write $\Phi=\int_X\varphi^xd\nu(x)$, $\Psi=\int_X\psi^xd\nu(x)$, $\overline{\alpha}$ (resp. $\overline{\beta}$) the representation of $L^\infty (X, \nu)$ on $H_{\Phi}$ (resp. $H_{\Psi}$) which extends $\alpha$ (resp. $\beta$). Let us write $\overline{\Phi}$ (resp. $\overline{\Psi}$) for the normal semi-finite faithful weight on $\pi_\Phi(A)''$ (resp. $\pi_\Psi(A)''$) which extends $\Phi$ (resp. $\Psi$), and $T$ (resp. $S$) for the normal semi-finite faithful operator-valued weight from $\pi_\Phi (A)''$ to $\overline{\alpha}(L^\infty(X, \nu))$ (resp. from $\pi_\Psi(A)''$ to $\overline{\beta}(L^\infty(X, \nu))$) which extends $\varphi$ (resp. $\psi$). We shall consider $A\underset{C_0(X)}{_\beta\otimes^m_\alpha}A$ as a $\mathbb{C}^*$-algebra on $H_{\Psi}\underset{\nu}{_{\overline{\beta}}\otimes_{\overline{\alpha}}}H_{\Phi}$, as explained in \ref{tensor}. Let us suppose, moreover, that the representation $\pi_{\Psi}\circ\alpha$ of $C_0(X)$ is square-integrable with respect to $\nu$, and let us write $\overline{\alpha_2}$ the representation of $L^\infty(X, \nu)$ on $H_{\Psi}$ constructed then. On the other hand, let us write, for all $f\in L^\infty(X, \nu)$ 
\[\hat{\beta}(f)=J_{\overline{\Phi}}\overline{\alpha}(f^*)J_{\overline{\Phi}}\]
Then :
\newline
(i) $\Gamma$ is square-integrable with respect to $\Phi$, and there exists an injective normal $*$-homomorphism $\overline{\Gamma_1}$ from $\pi_{\Phi}(A)''$ to $\pi_{\Psi}(A)''\underset{L^\infty (X, \nu)}{_{\overline{\beta}}*_{\overline{\alpha}}}\pi_{\Phi}(A)''$ such that, for all $a$ in $A$, we have $\overline{\Gamma_1}(\pi_{\Phi}(a))=\Gamma (a)$. 
\newline
(ii) the normal representation $\overline{\alpha_2}$ of $L^\infty (X, \nu)$ which extends $\pi_{\Psi}\circ\alpha$ is injective.  }

\begin{proof}
Let $x$, $e$ in $\gN_{\varphi}\cap\gN_{\Phi}$, $v$ in $H_{\Psi}$. Using \ref{propU}(iii), we get that :
\[\Gamma (x)(v\underset{\nu}{_{\overline{\beta}}\otimes_{\overline{\alpha}}}J_{\overline{\Phi}}\Lambda_{\Phi}(e))=
(1\underset{L^\infty (X, \nu)}{_{\overline{\beta}}\otimes_{\overline{\alpha}}}J_{\overline{\Phi}}\pi_{\Phi}(e)J_{\overline{\Phi}})U(v \underset{\nu}{_{\overline{\alpha_2}}\otimes_{\hat{\beta}}}\Lambda_{\Phi}(x))\]
and, therefore, if $v$ belongs to $D(_{\overline{\alpha_2}}(H_{\Psi}), \nu)$, we have :
\[\|\Gamma (x)(v\underset{\overline{\nu}}{_{\overline{\beta}}\otimes_{\overline{\alpha}}}J_{\overline{\Phi}}\Lambda_{\Phi}(e))\|\leq\|e\|\|R^{\overline{\alpha_2}, \nu}(v)\|\|\Lambda_{\Phi}(x)\|\]
and all such vectors $v\underset{\overline{\nu}}{_{\overline{\beta}}\otimes_{\overline{\alpha}}}J_{\overline{\Phi}}\Lambda_{\Phi}(e)$ in $_\Gamma(H_{\Psi}\underset{\nu}{_{\overline{\beta}}\otimes_{\overline{\alpha}}}H_{\Phi})$ are bounded with respect to $\Phi$. Using \ref{basic}, we get that the set of bounded vectors in $H_{\Psi}\underset{\nu^o}{_{\overline{\beta}}\otimes_{\overline{\alpha}}}H_{\Phi}$ is dense, and, therefore, there exists an unique normal representation $\overline{\Gamma_1}$ of $\pi_{\Phi}(A)''$ on this Hilbert space such that $\overline{\Gamma_1}(\pi_{\Phi}(a))=\Gamma (a)$, for all $a$ in $A$. Using now \ref{propU}(iii) again, we see that, for all $a$ in $A$, we have :
\[\overline{\Gamma_1}(\pi_{\Phi}(a))U=U(1\underset{L^\infty (X, \nu)}{_{\overline{\alpha_2}}\otimes_{\hat{\beta}}}\pi_{\Phi}(a))\] 
and, by continuity, we get, for any $x$ in $\pi_{\Phi}(A)''$ :
\[\overline{\Gamma_1}(x)U=U(1\underset{L^\infty (X, \nu)}{_{\overline{\alpha_2}}\otimes_{\hat{\beta}}}x)\] 
from which, $U$ being an isometry, we easily get that $\overline{\Gamma_1}$ is injective, which is (i). 
\newline
Moreover, we get, by continuity, for any $f\in L^\infty (X, \nu)$ :
\[\overline{\Gamma_1}(\overline{\alpha}(f))=\overline{\alpha_2}(f)\underset{L^\infty (X, \nu)}{_{\overline{\beta}}\otimes_{\overline{\alpha}}}1\]
from which we get the injectivity of $\overline{\alpha_2}$. 
\end{proof}

 \subsection{Definition}
 \label{defC*quantum}
Let $(X, \alpha, \beta, A, \Gamma,\varphi, \psi)$ be a $C_0(X)$-Hopf algebra, and let us suppose that $\varphi$ is left-invariant and $\psi$ right-invariant. Let $\nu$ be a positive Radon measure on $X$, with support $X$. Let us write $\Phi=\int_X\varphi^xd\nu(x)$, $\Psi=\int_X\psi^xd\nu(x)$; we shall say that  $(X, \alpha, \beta, A, \Gamma,\varphi,\psi, \nu)$ is a $\mathbb{C}^*$-quantum groupoid if, moreover :
\newline
(a) the representation $\pi_{\Psi}\circ\alpha$ of $C_0(X)$ is square-integrable with respect to $\nu$; 
\newline
(b) the representation $\pi_{\Phi}\circ\beta$ of $C_0(X)$ is square-integrable with respect to $\nu$.

\subsection{Theorem}
\label{AC*quantum}
{\it Let $(N, M, \alpha, \beta, \Gamma, T_1, T_2, \nu)$ be a measured quantum groupoid in the sense of Lesieur, and let us suppose that $\alpha(N)$ is central in $M$; let $W$ be the pseudo-multiplicative unitary of this measured quantum groupoid; let us consider the $\mathbb{C}^*$-algebra $A_n(W)$ (\ref{AW}, \ref{propmanag}(iv)); let us consider the $\mathbb{C}^*$-algebra $\mathbb{C}^*(\nu)$ (which is then nothing but the norm closure of $\gM_\nu$) and let $X$ be the spectrum of $\mathbb{C}^*(\nu)$; we shall denote by $\tilde{\nu}$ the restriction of $\nu$ to $\mathbb{C}^*(\nu)$, $\tilde{\alpha}$ (resp. $\tilde{\beta}$) the restriction of $\alpha$ (resp. $\beta$) to 
$\mathbb{C}^*(\nu)$. We shall denote $\tilde{T_1}$ (resp. $\tilde{T_2}$) the restrction of $T_1$ (resp. $T_2$) to $A_n(W)\cap\gM_{T_1}^+$ (resp. $A_n(W)\cap\gM_{T_2}^+$). Then :
\newline
(i) $(X, \tilde{\alpha}, A_n(W), \tilde{T_1})$ is a measured continuous field of $\mathbb{C}^*$-algebras; 
\newline
(ii) $(X, \tilde{\beta}, A_n(W), \tilde{T_2})$ is a measured continuous field of $\mathbb{C}^*$-algebras; 
\newline
(iii) $\Gamma_{|A_n(W)}$ is a $C_0(X)$-coproduct of $A_n(W)$;
\newline
(iv) for any $a$ in $\gN_{\tilde{T_1}}$, $b$ in $\gN_{\tilde{T_2}}$ and $x$ in $A_n(W)$, we have :
\[(id\underset{C_0(X)}{_{\tilde{\beta}}\otimes_{\tilde{\alpha}}^m}\omega_a)\Gamma(x)\in A_n(W)\]
\[(\omega_b\underset{C_0(X)}{_{\tilde{\beta}}\otimes_{\tilde{\alpha}}^m}id)\Gamma(x)\in A_n(W)\]
(v) $\tilde{T_1}$ is left-invariant, and $\tilde{T_2}$ is right-invariant with respect to $\Gamma$. 
\newline
(vi) $(X, \tilde{\alpha}, \tilde{\beta}, A_n(W), \Gamma_{|A_n(W)}, \tilde{\nu})$ is a $\mathbb{C}^*$-quantum groupoid, and $\Gamma_{|A_n(W)}$ is simplifiable. }

\begin{proof}
(i) and (ii) had been proved in \ref{continuousfield}, and (iii) in \ref{thcentral2}. If $a$ belongs to $A_n(W)\cap\gN_{T_1}\cap\gN_{\Phi_1}$ (where $\Phi_1=\nu\circ\alpha^{-1}\circ T_1$), we have :
\[(id\underset{C_0(X)}{_{\tilde{\beta}}\otimes_{\tilde{\alpha}}^m}\omega_a)\Gamma(x)=
(id\underset{N}{_\beta*_\alpha}\omega_{J_{\Phi_1}\Lambda_{\Phi_1}(a)})\Gamma(x)\]
which belongs to $A_n(W)$ by \ref{thgamma}. 
\newline
Moreover, if $a$ belongs to $A_n(W)\cap\gN_{T_1}$, and $f_n$ is an increasing family in $C_0(X)$ converging to $1$, we have :
  \begin{multline*}
 \|\Lambda_{T_1}(a(1-\alpha(f_n)))\|^2=\|(1-\alpha(f_n))T_1(a^*a)(1-\alpha(f_n))\|\\=\|T_1(a^*a)^{1/2}((1-\alpha(f_n))^2T_1(a^*a)^{1/2}\|
  \end{multline*}
 As $(1-\alpha(f_n))^2$ is decreasing to $0$, we get, by Dini's theorem, that $\Lambda_{T_1}(a\alpha(f_n))$ is norm converging to $\Lambda_{T_1}(a)$; therefore, $\omega_{a\alpha(f_n)}$ is norm converging to $\omega_a$, and we get that $(id\underset{C_0(X)}{\otimes^m}\omega_a)\Gamma (x)$ belongs to $A_n(W)$ for any $x$ in $A$ and $a$ in $A_n(W)\cap\gN_{T_1}$. 
 \newline
 The second part of (iv) is proved the same way. 
 \newline
 From the right-invariance of $T_2$, we get that, for all $x$ in $A_n(W)\cap\gM_{T_2}^+$, and $a$ in $A_n(W)\cap\gN_{T_1}\cap\gN_{\Phi_1}$ :
 \[T_2((id\underset{C_0(X)}{_\beta\otimes^m_\alpha}\omega_a)\Gamma (x))=\beta\circ\alpha^{-1}(\omega_a(T_2(x)))\]
If now $a$ belongs to $A_n(W)\cap\gN_{T_1}$,  using the same approximation of $a$ by $a\alpha(f_n)$, we get that 
 \[T_2((id\underset{C_0(X)}{_\beta\otimes^m_\alpha}\omega_{a\alpha(f_n)})\Gamma (x))\uparrow_nT_2((id\underset{C_0(X)}{_\beta\otimes^m_\alpha}\omega_a)\Gamma (x))\]
and, on the other hand, that $\beta\circ\alpha^{-1}(\omega_{a\alpha(f_n)}(T_2(x)))$ is converging to $\beta\circ\alpha^{-1}(\omega_a(T_2(x)))$, from which, we get, by linearity, for all $b$ in $A_n(W)\cap\gM_{T_2}^+$ and $a$ in $A_n(W)\cap \gN_{T_1}$ : 
 \[\tilde{T_2}((id\underset{C_0(X)}{_\beta\otimes^m_\alpha}\omega_a)\Gamma (b))=\beta\circ\alpha^{-1}(\omega_a(\tilde{T_2}(b)))\]
 and, taking the evaluation at the point $x\in X$, we get :
 \[\tilde{T_2}^x((id\underset{C_0(X)}{_\beta\otimes^m_\alpha}\omega_a)\Gamma (b))=\tilde{T_2}^x(b)\tilde{T_1}^x(a^*a)\]
 which gives that $\tilde{T_2}$ is right-invariant. Left-invariance of $\tilde{T_1}$ is proved the same way. 
\newline
Now it is trivial to prove that $(X, \tilde{\alpha}, \tilde{\beta}, A_n(W), \Gamma_{|A_n(W)}, \tilde{\nu})$ is a $\mathbb{C}^*$-quantum groupoid, and the fact that $\Gamma_{|A_n(W)}$ is simplifiable have been proved in \ref{thgamma}.  \end{proof}

\subsection{Theorem}
\label{thU2}
{\it Let $(X, \alpha, \beta, A, \Gamma,\varphi, \psi, \nu)$ be a $\mathbb{C}^*$-quantum groupoid (\ref{defC*quantum}). Let us write $\Phi=\int_X\varphi^xd\nu(x)$, $\Psi=\int_X\psi^xd\nu(x)$; let $\overline{\alpha}$ (resp. $\overline{\beta}$) be the representation of $L^\infty(X, \nu)$ on $H_\Phi$ (resp. $H_\Psi$) which extends $\pi_\Phi\circ\alpha$ (resp. $\pi_\Psi\circ\beta$); let $T$ (resp. $S$) be the normal semi-finite faithful operator-valued weight from $\pi_\Phi(A)''$ to $\overline{\alpha}(L^\infty(X, \nu))$ (resp. from $\pi_\Psi(A)''$ to $\overline{\beta}(L^\infty(X, \nu))$) which extends $\varphi$ (resp. $\psi$); let us write $\overline{\alpha_2}$ the representation of $L^\infty(X, \nu)$ which extends $\pi_{\Psi}\circ\alpha$, and $\overline{\beta_1}$ the representation of $L^\infty(X, \nu)$ which extends $\pi_{\Phi}\circ\beta$. Then :
\newline
(i) there exists an isomorphism between $\pi_{\Phi}(A)''$ and $\pi_{\Psi}(A)''$ which makes the representations $\pi_{\Phi}$ and $\pi_{\Psi}$ equivalent, the representations $\overline{\alpha}$ and $\overline{\alpha_2}$ equivalent, and the representations $\overline{\beta}$ and $\overline{\beta_1}$ equivalent; 
\newline
(ii) there exists a coproduct $\overline{\Gamma}$ from $M=\pi_{\Phi_1}(A)''$ into $M\underset{L^\infty(X, \nu)}{_{\overline{\beta_1}}*_{\overline{\alpha}}}M$ which makes $(L^\infty(X, \nu), M, \overline{\alpha}, \overline{\beta_1}, \overline{\Gamma})$ be a Hopf-bimodule;
\newline
(iii) the normal semi-finite faithful operator-valued weight $T$ from $M$ to $\overline{\alpha}(L^\infty(X, \nu))$ is a left-invariant operator-valed weight, and the isomorphism obtained in (i) sends $S$ on a right-invariant weight from $M$ to $\overline{\beta_1}(L^\infty(X, \nu))$, we shall denote $T'$. 
\newline
(iv) $(L^\infty(X, \nu), M, \overline{\alpha}, \overline{\beta_1}, \overline{\Gamma}, T, T')$ is a measured quantum groupoid in the sense of Lesieur ([L1], [L2]). 
\newline
(v) if $\Gamma$ is simplifiable, then $A\subset A_n(W)$, where $W$ is the pseudo-multiplicative unitary of the measured quantum groupoid constructed in (iv). }

\begin{proof}
Using twice \ref{thU}, we can construct an injective $*$-homomorphism $\overline{\Gamma_1}$ (resp. $\overline{\Gamma_2}$) from $\pi_{\Phi}(A)''$ (resp. $\pi_{\Psi}(A)''$) to $\pi_{\Psi}(A)''\underset{L^\infty(X, \nu)}{_{\overline{\beta}}*_{\overline{\alpha}}}\pi_{\Phi}(A)''$ such that $\overline{\Gamma_1}(\pi_{\Phi}(a))=\Gamma (a)$ (resp. $\overline{\Gamma_2}(\pi_{\Psi}(a))=\Gamma (a)$), for all $a\in A$. So, we get that  $\overline{\Gamma_1}(\pi_{\Phi}(A)'')=\overline{\Gamma_2}(\pi_{\Psi}(A)'')$, and $\overline{\Gamma_2}^{-1}\overline{\Gamma_1}$ gives the isomorphism required and (i). With the help of that isomorphism, one gets then (ii). If now $x$ is in $\gN_{S}\cap\gN_{\Psi}$, and $a$ in $\gN_{T}$, we get, by the left-invariance of $\varphi$, that :
\[T((\omega_{J_{\overline{\Psi}}\Lambda_{\Psi}(x)}\underset{L^\infty(X, \nu)}{_{\overline{\beta}}*_{\overline{\alpha}}}id)\Gamma (a^*a))
=\alpha(<T(a^*a)J_{\overline{\Psi}}\Lambda_{\Psi}(x),J_{\overline{\Psi}}\Lambda_{\Psi}(x)>_{\overline{\beta}, \nu})\]
which can be written, by density :
\[(id \underset{L^\infty(X, \nu)}{_{\overline{\beta}}*_{\overline{\alpha}}}T)\Gamma (a^*a)=
T(a^*a)\underset{L^\infty(X, \nu)}{_{\overline{\beta}}\otimes_{\overline{\alpha}}}1\]
Using now (i), we get that :
\[(id \underset{L^\infty(X, \nu)}{_{\overline{\beta_1}}*_{\overline{\alpha}}}T)\overline{\Gamma} (a^*a)=
T(a^*a)\underset{L^\infty(X, \nu)}{_{\overline{\beta_1}}\otimes_{\overline{\alpha}}}1\]
and, by normality, we get, for all $a\in\gN_T$ :
\[(id \underset{L^\infty(X, \nu)}{_{\overline{\beta_1}}*_{\overline{\alpha}}}T)\overline{\Gamma} (a^*a)=
T(a^*a)\underset{L^\infty(X, \nu)}{_{\overline{\beta_1}}\otimes_{\overline{\alpha}}}1\]
which proves that $T$ is left-invariant. 
\newline
The proof of the right-invariantness of $S$ is the same, which finishes the proof of (iii).
\newline
 As $\alpha (L^\infty(X, \nu))\subset Z(M)$ (resp. $\beta_1(L^\infty(X, \nu))\subset Z(M)$), we get that $\sigma_t^{\overline{T_2}}(\alpha (f))=\alpha (f)$ (resp. $\sigma_t^{\overline{T_1}}(\beta_1(f))=\beta_1(f)$) for all $f\in L^\infty(X, \nu)$, which finishes the proof of (iv). 
 \newline
 If $a$ belongs to $\gN_{\varphi}^*\cap\gN_{\Phi}^*$, and $b$ to $\gN_{\varphi}\cap\gN_{\Phi}$, we have :
 \begin{multline*}
 (id\underset{C_0(X)}{\otimes^m}\omega_b)\Gamma (a)=(id\underset{L^\infty(X, \nu)}{_{\overline{\beta}}\otimes_{\overline{\alpha}}}\omega_{J_{\overline{\Phi}}\Lambda_{\Phi}(b)})\overline{\Gamma}(a)\\=(i*\omega_{J_{\overline{\Phi}}\Lambda_{\Phi}((b^*b), \Lambda_{\Phi}(a^*)})(W)
 \end{multline*}
 which belongs to $A_n(W)$. 
 \newline
 By continuity, this remains true for any $a$ in $A$; moreover, if $b$ belongs to $\gN_{\varphi}$, and $f_n$ is an increasing family in $C_0(X)$ converging to $1$, we have :
  \begin{multline*}
 \|\Lambda_T(b(1-\alpha(f_n)))\|^2=\|(1-\alpha(f_n))T(b^*b)(1-\alpha(f_n))\|\\=\|T(b^*b)^{1/2}((1-\alpha(f_n))^2T(b^*b)^{1/2}\|
  \end{multline*}
 As $(1-\alpha(f_n))^2$ is decreasing to $0$, we get, by Dini's theorem, that $\Lambda_T(b\alpha(f_n))$ is norm converging to $\Lambda_T(b)$; therefore, $\omega_{b\alpha(f_n)}$ is norm converging to $\omega_b$, and we get that $(id\underset{C_0(X)}{\otimes^m}\omega_b)\Gamma (a)$ belongs to $A_n(W)$ for any $a$ in $A$ and $b$ in $\gN_{\varphi}$; so, by the simplifiability of $\Gamma$, we obtain (v).   \end{proof} 
 
\subsection{Example}
\label{gd4}
Let $\mathcal G$ be a locally compact groupoid, as described in \ref{gd}. Then $(\mathcal G^{(0)}, \hat{r}, \hat{s}, C_0(\mathcal G), \Gamma_{\mathcal G}, \lambda^u, \lambda_u, \mu)$ is a $\mathbb{C}^*$-quantum groupoid, and $\Gamma_{\mathcal G}$ is simplifiable. Applying then \ref{thU2}, we recover the fundamental example of a measured quantum groupoid described in \ref{fe}, and, in particular, we get again that $C_0(\mathcal G)\subset A_n(W_{\mathcal G})$, which was remarked in \ref{gd}. 

\section{Continuous field of locally compact quantum groups}
\label{fieldqg}
In that last chapter, we define a notion of continuous field of locally compact quantum groups (\ref{deffieldqg}), which was underlying in [B2], and which gives examples of the $\mathbb{C}^*$-quantum groups introduced in \ref{C*central}. We show that these are exactly Lesieur's measured quantum groupoids with central basis, such that the dual object is of the same kind (\ref{bicentral}). We finish by recalling concrete examples (\ref{SU}, \ref{ax+b}, \ref{Emu}) given by Blanchard, which give, examples of measured quantum groupoids ([L2], 6.8.3).

\subsection{Definition}
\label{deffieldqg}
A 6-uple $(X, \alpha, A, \Gamma^x, \varphi^x, \psi^x)$ will be called a continuous field of locally compact quantum groups if :
\newline
(i) $(X, \alpha, A, \varphi^x)$ and $(X, \alpha, A, \psi^x)$ are measurable continuous fields of $\mathbb{C}^*$-algebras; 
\newline
(ii) for any $x\in X$, there exists a simplifiable coproduct $\Gamma^x$ from $A^x$ to $M(A^x\otimes_m A^x)$ such that $(A^x, \Gamma^x, \varphi^x, \psi^x)$ is a locally compact quantum group, in the sense of [KV1]. 
\newline
Let us recall that simplifiable means that the closed linear set generated by $\Gamma^x(A^x)(A^x\otimes_m 1)$ (resp. $\Gamma^x(A^x)(1\otimes_m A^x)$) is equal to $A^x\otimes_m A^x$. 

\subsection{Theorem}
\label{thfieldqg}
{\it Let $(X, \alpha, A, \Gamma^x, \varphi^x, \psi^x)$ be a continuous field of locally compact quantum groups; then :
\newline
(i) there exists a simplifiable $C_0(X)$-coproduct $\Gamma$  on $A$ such that, for any $x\in X$, and $a\in A$, we have $\Gamma(a)^x=\Gamma^x(a^x)$; 
\newline
(ii) the continuous field of weights $\varphi=(\varphi^x)_{x\in X}$ (resp. $\psi=(\psi^x)_{x\in X}$) is left-invariant (resp. right-invariant).  
\newline
(iii) $(X, \alpha, \alpha, A, \Gamma, \varphi, \psi, \nu)$ is a $\mathbb{C}^*$-quantum groupoid, for any Radon measure $\nu$ on $X$, whose support is $X$. }

\begin{proof}
We get the existence of $\Gamma$ by \ref{AtensA}; we know that the elements $(\omega_{\Lambda_{\varphi^x}(b)}\otimes id)\Gamma^x(a)$, for all $a\in A^x$ and $b\in \gN_{\varphi^x}$ belongs to $A^x$, and, as $\Gamma^x$ is simplifiable, we know that $A^x$ is the norm closure of the linear set generated by such elements. From which we get that the elements $(\omega_b\underset{C_0(X)}\otimes^m id)\Gamma (a)$, for all $a\in A$ and $b\in \gN_{\varphi}$ belong to $A$, and, moreover, that $A$ is the norm closure of the linear set generated by such elements. 
\newline
Let $R^x$ be the unitary co-inverse of the locally compact quantum group  $(A^x, \Gamma^x, \varphi^x, \psi^x)$; we easily get that $\varphi\circ R=(\varphi^x\circ R^x)_{x\in X}$ is a continuous field of weights, and that, for all $x\in X$, $\psi^x$ is proportional to $\varphi^x\circ R^x$;  using it, we get the same way that the elements $(id\underset{C_0(X)}\otimes^m \omega_a)\Gamma (b)$ belong to $A$ for any $a$ in $\gN_{\varphi\circ R}$ and $b$ in $A$, and that $\Gamma$ is simplifiable. 
\newline
The left-invariance of $\varphi^x$ gives that, for all $a\in\gM_{\varphi^x}$ and $b\in \gN_{\varphi^x\circ R^x}$, we have :
\[\varphi^x((\omega_{\Lambda_{\varphi^x\circ R^x}(b)}\otimes id)\Gamma^x(a))=\varphi^x(a)\varphi^x\circ R^x(b^*b)\]
from which we get that $\varphi$ is left-invariant; we obtain the same way that $\psi$ is right-invariant. Now, using \ref{defC*quantum}, we finish the proof. \end{proof}

\subsection{Theorem}
\label{bicentral}
{\it Let $\gG=(N, M, \alpha, \beta, \Gamma, T, R, \nu, \tau_t, \gamma_t)$ be a generalized measured quantum 
groupoid in the sense of Lesieur, and let us denote $\widehat{\gG}=(N, \widehat{M}, \alpha, \widehat{\beta}, \widehat{\Gamma}, \widehat{T}, \widehat{R}, \nu, \widehat{\tau_t}, \gamma_{-t})$ the dual generalized quantum groupoid; let $W$ and $\widehat{W}$ be the pseudo-multiplicative associated; let us suppose that $\alpha (N)$ is central in both $M$ and $\widehat{M}$; then :
\newline
(i) we have $\alpha=\beta=\hat{\beta}$, $\gamma_t=id$, and $(N, M, \alpha, \alpha, \Gamma, T, T\circ R, \nu)$ and $(N, \widehat{M}, \alpha, \alpha, \widehat{\Gamma}, \widehat{T}, \widehat{T}\circ\widehat{R}, \nu)$ are both measured quantum groupoids in the sense of Lesieur. 
\newline
(ii) let $X$ be the spectrum of the $\mathbb{C}^*$algebra $\mathbb{C}^*(\nu)$ and $\tilde{\alpha}$ the restriction of $\alpha$ to $\mathbb{C}^*(\nu)=C_0(X)$; for all $x\in X$ there exists on the $\mathbb{C}^*$-algebra $A_n(W)^x$ a coproduct $\Gamma^x$ and a weight $\varphi^x$ such that $(X, \tilde{\alpha}, A_n(W), \Gamma^x, \varphi^x)$ is a continuous field of locally compact quantum groups in the sense of \ref{deffieldqg}. Let $W^x$ be the multiplicative unitary associated to the locally quantul group $(A_n(W)^x,\Gamma^x, \varphi^x)$. 
\newline
(iii) for all $x\in X$ there exists on the $\mathbb{C}^*$-algebra $A_n(\widehat{W})^x$ a coproduct $\widehat{\Gamma}^x$ and a weight $\widehat{\varphi}^x$ such that $(X, \tilde{\alpha}, A_n(\widehat{W}), \widehat{\Gamma}^x, \widehat{\varphi}^x)$ is a continuous field of locally compact quantum groups in the sense of \ref{deffieldqg}. Let $\widehat{W}^x$ be the multiplicative unitary associated to the locally quantum group $(A_n(\widehat{W})^x=\widehat{\Gamma}^x, {\varphi}^x)$. Then, we get that $\widehat{W}^x,\widehat{W^x}$ and that the locally compact quantum group $(A_n(\widehat{W})^x,\widehat{\Gamma}^x, {\varphi}^x)$ is the dual of $(A_n(W)^x,\Gamma^x, \varphi^x)$. }

\begin{proof}
As $\alpha (N)$ is central in $M$, we get $\alpha=\hat{\beta}$, and applying this result to $\gG$, we get that $\alpha=\beta$, which is (i). 
\newline
Let us write $\varphi$ for the continuous field of weights which is the restriction of $T$ to $A_n(W)\cap\gM_T^+$. 
Let us consider now the $\mathbb{C}^*$-quantum groupoid 
\[(X, \tilde{\alpha}, \tilde{\alpha}, A_n(W), \Gamma_{|A_n(W)}, \varphi, \varphi\circ R)\]
constructed in \ref{AC*quantum}. 
\newline
Using \ref{AtensA}, we get that the coproduct $\Gamma_{|A_n(W)}$ leads, for each $x\in X$,  to an involutive homomorphism $\Gamma^x$ from $A_n(W)^x$ to $A_n(W)^x\otimes_m A_n(W)^x$; the left-invariance of $T$ leads to the formula $(id\otimes\varphi^x)\Gamma^x(a)=\varphi^x(a)$, for all $a$ in $\gM_{\varphi^x}^+$, and we obtain 
that $(X, \tilde{\alpha}, A_n(W), \Gamma^x, \varphi^x)$ is a continuous field of locally compact quantum groups in the sense of \ref{deffieldqg}, which is (ii). 
\newline
We obtain the same way that there exists a continuous field of locally compact quantum groups $(X, \tilde{\alpha}, A_n(\widehat{W}), \widehat{\Gamma}^x, \widehat{\varphi}^x)$; moreover, as $W$, in that context, appears as a continuous field of multiplicative unitaries in $\mathcal L(\mathcal E_{\varphi}\otimes_{C_0(X)}E_{\varphi})$, in the sense of ([B2] 4.1), and, for all $x\in X$, we get that $W^x$ is the multiplicative unitary of the locally compact quantum group $(A_n(W)^x, \Gamma^x, \varphi^x)$. Then, the duality between $W$ and $\widehat{W}$ leads to the duality between $(A_n(W)^x, \Gamma^x, \varphi^x)$ and $(A_n(\widehat{W}), \widehat{\Gamma}^x, \widehat{\varphi}^x)$, which is (iii). 
\end{proof}

\subsection{Example}
\label{SU}
As in ([B1], 7.1), let us consider the $\mathbb{C}^*$-algebra $A$ whose generators $\alpha$, $\gamma$ and $f$ verify :
\newline
(i) $f$ commutes with $\alpha$ and $\gamma$;
\newline
(ii) the spectrum of $f$ is $[0,1]$;
\newline
(iii) the matrix 
$\left(\begin{array}{cc}
\alpha&-f\gamma\\
\gamma & \alpha^*
\end{array}
\right)$ is unitary in $M_2(A)$. 
Then, using the sub $\mathbb{C}^*$-algebra generated by $f$, $A$ is a $C([0,1])$-algebra; let us consider now $A$ as a $C_0(]0,1])$-algebra.  Then, Blanchard had proved ([B2] 7.1) that $A$ is a continuous field over $]0,1]$ of $\mathbb{C}^*$-algebras, and that, for all $q\in]0,1]$, we have $A^q=SU_q(2)$, where the $SU_q(2)$ are the compact quantum groups constructed by Woronowicz and $A^1=C(SU(2))$. 
\newline
Moreover, using the coproducts $\Gamma^q$ defined by Woronowicz as 
\[\Gamma^q(\alpha)=\alpha\otimes\alpha-q\gamma^*\otimes\gamma\]
\[\Gamma^q(\gamma)=\gamma\otimes\alpha+\alpha^*\otimes\gamma\]
and the (left and right-invariant) Haar state $\varphi^q$, which verifies : 
\newline
$\varphi^q(\alpha^k\gamma^{*m}\gamma^n)=0$, for all $k\geq 0$, and $m\not=n$, 
\newline
$\varphi^q(\alpha^{*|k|}\gamma^{*m}\gamma^n)=0$, for all $k<0$, and $m\not=n$, 
\newline
and $\varphi^q((\gamma^*\gamma)^m)=\frac{1-q^2}{1-q^{2m+2}}$, 
\newline
we obtain this way a continuous field of compact quantum groups; this leads to put on $A$ a structure of $\mathbb{C}^*$ quantum groupoid (of compact type, because $1\in A$). 
\newline
This structure is given by a coproduct $\Gamma$ which is $C_0([0,1])$-linear from $A$ to $A\underset{C_0(]0,1])}{\otimes^m}A$, and given by :
\[\Gamma(\alpha)=\alpha\underset{C_0(]0,1])}{\otimes^m}\alpha-f\gamma^*\underset{C_0(]0,1])}{\otimes^m}\gamma\]
\[\Gamma(\gamma)=\gamma\underset{C_0(]0,1])}{\otimes^m}\alpha+\alpha^*\underset{C_0(]0,1])}{\otimes^m}\gamma\]
and by a conditional expectation $E$ from $A$ on $M(C_0(]0,1]))$ given by :
\newline
$E(\alpha^k\gamma^{*m}\gamma^n)=0$, for all $k\geq 0$, and $m\not=n$
\newline
$E(\alpha^{*|k|}\gamma^{*m}\gamma^n)=0$, for all $k<0$, and $m\not=n$
\newline
$E((\gamma^*\gamma)^m)$ is the bounded function $x\mapsto\frac{1-q^2}{1-q^{2m+2}}$. 
Then $E$ is both left and right-invariant with respect to $\Gamma$. 

\subsection{Example}
\label{ax+b}
One can find in [B2] another example of a continuous field of locally compact quantum group. Namely, in ([B2], 7.2), Blanchard constructs a $\mathbb{C}^*$-algebra $A$ which is a continuous field of $\mathbb{C}^*$-algebras over $X$, where $X$ is a compact included in $]0, 1]$, with $1\in X$. For any $q\in X$, $q\not=1$, we have $A^q=SU_q(2)$, and $A^1=\mathbb{C}^*_r(G)$, where $G$ is the "$ax+b$" group. ([B2], 7.6). 
\newline
Moreover, he constructs a coproduct (denoted $\delta$) ([B2] 7.7(c)), which is simplifiable ([B2] 7.8 (b)), and a continuous field of weights (denoted $\Phi$) ([B2] 7.2.3), which is left-invariant (end of remark after [B2] 7.2.3). 
\newline
Finally, he constructs a unitary $U$ in $\mathcal L(\mathcal E_\Phi)$ ([B2] 7.10), with which it is possible to construct a co-inverse $R$ of $(A, \delta)$, which leads to the fact that $\Phi\circ R$ is right-invariant. 

\subsection{Example}
\label{Emu}
Let us finish by quoting a last example given by Blanchard in ([B2], 7.4) : for $X$ compact in $[1, \infty[$, with $1\in X$, he constructs a $\mathbb{C}^*$-algebra which is a continuous field over $X$ of $\mathbb{C}^*$ algebras, whose fibers, for $\mu\in X$, are $A^\mu=E_\mu(2)$, with a coproduct $\delta$ and a continuous field of weights $\Phi$, which is left-invariant. As in \ref{ax+b}, he then constructs a unitary $U$ on $\mathcal L(\mathcal E_\Phi)$, which will lead to a co-inverse, and, therefore, to a right-invariant $\mathbb{C}^*$-weight. 

\subsection{Example ([L2], 6.8.1)}
Let us return to \ref{deffieldqg}; let $I$ be a (discrete) set, and, for all $i$ in $I$, let $(A_i, \Gamma_i, \varphi_i, \psi_i)$ be a locally compact quantum group; then the sum $\oplus_i A_i$ is a continuous field of locally compact quantum groups, and can be given a natural structure of $\mathbb{C}^*$-quantum groupoid, described in ([L2], 6.8.1). 


\section{Bibliography}
[AR] C. Anantharaman-Delaroche et J. Renault : Amenable Groupo-ids; Monographies de l'Enseignement Math\'ematique, {\bf 36}. L'Enseigne-ment
Math\'ematique, Gen\`eve, 2000. 196 pp
\newline\indent
[B] S. Baaj : Prolongement d'un poids, {\it C.R. Acad. Sci. Paris}, {\bf 288} (1979), 1013-1015.
\newline\indent
[BS] S. Baaj et G. Skandalis : Unitaires multiplicatifs et dualit\'{e} pour les produits  crois\'{e}s
de
$\mathbb{C}^*$-alg\`{e}bres, {\it Ann. Sci. ENS}, {\bf 26} (1993), 425-488.
\newline\indent
[B1] E. Blanchard : Tensor products of $\mathbb{C}(X)$-algebras over $\mathbb{C}(X)$, {\it Ast\'{e}trisque} {\bf 232} (1995), 81-92. 
\newline\indent
[B2] E. Blanchard : D\'{e}formations de $\mathbb{C}^*$-alg\`{e}bres de Hopf, {\it Bull. Soc. Math. France}, {\bf 24} (1996), 141-215.
\newline\indent
[BSz1] G. B\"{o}hm and K. Szlach\'{a}nyi : A Coassociative $\mathbb{C}^*$-Quantum group with Non
Integral Dimensions, {\it Lett. Math. Phys.}, {\bf 38} (1996), 437-456.
\newline\indent
[BSz2] G. B\"{o}hm and K. Szlach\'{a}nyi : Weak $\mathbb{C}^*$-Hopf Algebras : the coassociative
symmetry of non-integral dimensions, in Quantum Groups and Quantum spaces {\it Banach Center
Publications}, {\bf 40} (1997), 9-19.
\newline\indent
[Co1] F. Combes : Poids sur une $\mathbb{C}^*$-algbre, {\it J. math. pures et appl.}, {\bf 47} (1968), 57-100.
\newline\indent
[Co2] F. Combes : Syst\`emes Hilbertiens \`a gauche et repr\'esentation de Gelfand-Segal, in {\it Operator Algebras and Group representations vol 1; proceedings of the International Conference held in Neptun, sept 80 (Romania)}, Monographs and Studies in mathematics 17, 71-107
\newline\indent
[C1] A. Connes: On the spatial theory of von Neumann algebras, 
{\it J. Funct. Analysis}, {\bf 35} (1980), 153-164.
\newline\indent
[C2] A. Connes: Non commutative Geometry, Academic Press, 1994
\newline\indent
[E1] M. Enock : Inclusions irr\'eductibles de facteurs et
unitaires multiplicatifs II, {\it J. Funct. Analysis}, {\bf 154} (1998), 67-109.
\newline\indent
[E2] M. Enock : Quantum groupoids of compact type, {\it J. Inst. Math. Jussieu}, {\bf 4} (2005), 29-133. 
\newline\indent
[E3] M. Enock : Inclusions of von Neumann algebras and quantum groupo\"{\i}ds III, preprint 369 Institut de Math\'ematiques de Jussieu, math.OA/0405480, to be published in {\it J. Funct. Analysis}. 
\newline\indent
[EN] M. Enock, R. Nest : Inclusions of factors,
multiplicative unitaries and Kac algebras, {\it J. Funct. Analysis}, {\bf137} (1996),
466-543.
\newline\indent
[ES] M. Enock, J.-M. Schwartz : Kac algebras and Duality of locally compact Groups,
Springer-Verlag, Berlin, 1989.
\newline\indent
[EV] M. Enock, J.-M. Vallin : Inclusions of von Neumann algebras and quantum groupo\"{\i}ds,
{\it J. Funct. Analalysis}, {\bf 172} (2000), 249-300.
\newline\indent
[GM] T. Giordano and J.A. Mingo : Tensor products of $\mathbb{C}^*$-algebras over abelian subalgebras, {\it J. London Math. soc. }, {\bf 55} (1997), 170-180.
\newline\indent
[J] V. Jones: Index for subfactors, {\it Invent. Math.},  {\bf 72} (1983), 1-25.
\newline\indent
[Ka] G.G. Kasparov : Equivariant KK-theory and the Novikov conjecture, {\it Invent. Math.}, {\bf 91} (1988), 147-201.
\newline\indent
[KV1] J. Kustermans and S. Vaes : Locally compact quantum groups, {\it Ann. Sci. ENS}, {\bf 33} (2000), 837-934.
\newline\indent
[KV2] J. Kustermans and S. Vaes, Locally compact quantum groups in the von Neumann algebraic setting, {\it Math. Scand.}, {\bf 92} (2003), 68-92. 
\newline\indent
[La] E.C. Lance : Hilbert $\mathbb{C}^*$-modules, a toolkit for operator algebraists, London Mat. Soc. lecture Note Series 210, Cambridge University Press, 1995. 
\newline\indent
[L1] F. Lesieur : thesis, University of Orleans, available at :
\newline
http://tel.ccsd.cnrs.fr/documents/archives0/00/00/55/05
\newline\indent
[L2] F. Lesieur : Measured Quantum groupoids, math.OA/0409380
\newline\indent
[L3] F. Lesieur : Measured Quantum Groupoids- Duality, \newline
math.OA/0504104
\newline\indent
[MN] T. Masuda and Y. Nakagami : A von Neumann Algebra framework for
the duality of the quantum groups; {\it Publ. RIMS Kyoto}, {\bf 30} (1994), 799-850.
\newline\indent
[MNW] T. Masuda, Y. Nakagami and S.L. Woronowicz, A  $\mathbb{C}^*$-algebraic framework for quantum groups, {\it Internat. J. Math.}, {\bf 14} (2003), 903-1001. 
\newline\indent
[NV1] D. Nikshych, L. Va\u{\i}nerman : Algebraic versions of a finite dimensional quantum
groupoid, in Lecture Notes in Pure and Applied Mathematics, Marcel Dekker, 2000.
\newline\indent
[NV2] D. Nikshych, L. Va\u{\i}nerman : A characterization of depth 2 subfactors of $II_1$
factors, {\it J. Funct. Analysis}, {\bf 171} (2000), 278-307.
\newline\indent
[P] G.K. Pedersen : $\mathbb{C}^*$-Algebras and their Automorphism Groups, London Mathematical Society Monographs 14, Academic Press, 1979.
\newline\indent
[R1] J. Renault : A Groupoid Approach to $\mathbb{C}^*$-Algebras, {\it Lecture Notes in
Math.} {\bf 793}, Springer-Verlag
\newline\indent
[R2] J. Renault : The Fourier algebra of a measured groupoid and its multipliers, {\it J. Funct. Analysis}, {\bf 145} (1997), 455-490.
\newline\indent
[S1] J.-L. Sauvageot : Produit tensoriel de $Z$-modules et applications, in Operator
Algebras and their Connections with Topology and Ergodic Theory, Proceedings Bu\c{s}teni, Romania,
1983, {\it Lecture Notes in
Math.} {\bf 1132}, Springer-Verlag, 468-485.
\newline\indent
[S2] J.-L. Sauvageot : Sur le produit tensoriel
relatif d'espaces de Hilbert,  {\it J. Operator Theory}, {\bf 9} (1983), 237-352.
\newline\indent
[Sz] K. Szlach\'{a}nyi : Weak Hopf algebras, in Operators Algebras and Quantum Field Theory,
S. Doplicher, R. Longo, J.E. Roberts, L. Zsido editors, International Press, 1996.
\newline\indent
[T] M. Takesaki : Theory of Operator Algebras II, Springer, Berlin, 2003. 
\newline\indent
[Val1] J.-M. Vallin : $\mathbb{C}^*$-alg\`ebres de Hopf et $\mathbb{C}^*$-alg\`ebres de Kac, {\it Proc. London Math. Soc.}, {\bf 50} (1985), 131-174.
\newline\indent
[Val2] J.-M. Vallin : Bimodules de Hopf
et Poids op\'eratoriels de Haar, {\it J. Operator theory}, {\bf 35} (1996), 39-65
\newline\indent
[Val3] J.-M. Vallin : Unitaire pseudo-multiplicatif associ\'e \`a un groupo\"{\i}de;
applications \`a la moyennabilit\'e,  {\it J. Operator theory}, {\bf 44} (2000), 347-368.
\newline\indent
[Val4] J.-M. Vallin : Groupo\"{\i}des quantiques finis, {\it J. Algebra}, {\bf 239} (2001), 215-261.
\newline\indent
[Val5] J.-M. Vallin : Multiplicative partial isometries and finite quantum groupoids, in Locally Compact Quantum Groups and Groupoids, IRMA Lectures in Mathematics and Theoretical Physics 2, V. Turaev, L. Vainerman editors, de Gruyter, 2002. 
\newline\indent
[W1] S.L. Woronowicz : Tannaka-Krein duality
for compact matrix pseudogroups. Twisted
$SU(N)$ group. {\it Invent. Math.}, {\bf 93} (1988), 35-76.
\newline\indent
[W2] S.L. Woronowicz : Compact quantum group, in "Sym\'{e}tries quantiques" (Les Houches, 1995), North-Holland, Amsterdam (1998),
845-884.
\newline\indent
[W3] S.L. Woronowicz : From multiplicative unitaries to quantum groups, {\it Int. J. Math.}, {\bf 7}
(1996), 127-149. 
\newline\indent
[Y1] T. Yamanouchi : Duality for actions and coactions of
measured Groupoids on von Neumann Algebras, {\it Memoirs of the A.M.S.}, {\bf 101} (1993), 1-109.
\newline\indent
[Y2] T. Yamanouchi : Duality for generalized Kac algebras and a characterization on
finite groupoids algebras, {\it J. Algebra}, {\bf 163} (1994), 9-50.
\newline\indent

\end{document}